\documentclass[11pt]{article}
\usepackage{amsmath,amssymb,type1cm,bm,color,amsfonts}
\usepackage[normalem]{ulem} 

\usepackage{mathrsfs}
\usepackage{longtable,tabu}
\usepackage[flushleft]{threeparttable}
\usepackage[numbers]{natbib}
\usepackage{hyperref}
\usepackage{booktabs}
\usepackage{multirow}
\usepackage{mathtools}
\usepackage{subfig}
\usepackage{array}
\usepackage{caption}
\usepackage{color}
\usepackage{graphicx,psfrag,epsf}
\usepackage{anysize}

\marginsize{1.1in}{1.1in}{0.55in}{0.8in} %

\newcommand{\ep}{\varepsilon}

\DeclareMathOperator\diag{diag}

\DeclareMathOperator\supp{supp}

\DeclareMathOperator\E{\mathbb{E}}

\DeclareMathOperator\Pro{\mathbb{P}}

\DeclareMathOperator\sgn{sgn}
\DeclareMathOperator\vect{vec}
\DeclareMathOperator\vech{vech}
\DeclareMathOperator\FDR{FDR}
\DeclareMathOperator\FDP{FDP}

\def\bA{\mathbf{A}}
\def\bB{\mathbf{B}}
\def\bC{\mathbf{C}}
\def\bD{\mathbf{D}}
\def\bE{\mathbf{E}}
\def\bF{\mathbf{F}}

\def\bH{\mathbf{H}}
\def\bI{\mathbf{I}}

\def\bM{\mathbf{M}}

\def\bT{\mathbf{T}}
\def\bU{\mathbf{U}}
\def\bV{\mathbf{V}}

\def\bX{\mathbf{X}}

\def\ba{\mathbf{a}}
\def\bb{\mathbf{b}}

\def\be{\mathbf{e}}
\def\bff{\mathbf{f}}

\def\bt{\mathbf{t}}

\def\bv{\mathbf{v}}

\def\bx{\mathbf{x}}
\def\by{\mathbf{y}}
\def\bz{\mathbf{z}}
\def\cA{\mathcal{A}}
\def\cB{\mathcal{B}}

\def\cS{\mathcal{S}}

\def\bgamma{\boldsymbol{\gamma}}
\def\bbeta{\boldsymbol{\beta}}

\def\bdelta{\boldsymbol{\delta}}
\def\bLambda{\boldsymbol{\Lambda}}
\def\bTheta{\boldsymbol{\Theta}}
\def\btheta{\boldsymbol{\theta}}
\def\bbbeta{\boldsymbol{\eta}}

\def\bSigma{\boldsymbol{\Sigma}}
\def\blambda{\boldsymbol{\lambda}}
\def\bxi{\boldsymbol{\xi}}

\def\ep{\varepsilon}

\def\bep{\boldsymbol{\varepsilon}}
\def\bzero{\mathbf{0}}

\def\what{\widehat}
\def\wtilde{\widetilde}

\newtheorem{thm}{Theorem}
\newtheorem{lem}{Lemma}
\newtheorem{con}{Condition}
\newtheorem{defi}{Definition}
\newtheorem{res}{Result}
\newtheorem{proc}{Procedure}

\begin{document}
	
\title{IPAD: Stable Interpretable Forecasting with Knockoffs Inference %
\thanks{
Yingying Fan is Dean's Associate Professor in Business Administration, Data Sciences and Operations Department, Marshall School of Business, University of Southern California, Los Angeles, CA 90089, USA (E-mail: \textit{fanyingy@marshall.usc.edu}). %
Jinchi Lv is McAlister Associate Professor in Business Administration, Data Sciences and Operations Department, Marshall School of Business, University of Southern California, Los Angeles, CA 90089, USA (E-mail: \textit{jinchilv@marshall.usc.edu}). %
Mahrad Sharifvaghefi is Ph.D. candidate, Department of Economics, University of Southern California, Los Angeles, CA 90089, USA (E-mail: \textit{sharifva@usc.edu}). %
Yoshimasa Uematsu is Assistant Professor, Department of Economics and Management, Tohoku University, Sendai 980-8576, Japan (E-mail: \textit{yoshimasa.uematsu.e7@tohoku.ac.jp}). Most of this work was completed while Uematsu visited USC Marshall as JSPS Overseas Research Fellow and Postdoctoral Scholar. %
This work was supported by NIH Grant 1R01GM131407-01, NSF CAREER Award DMS-1150318, a grant from the Simons Foundation, Adobe Data Science Research Award, and a Grant-in-Aid for JSPS Overseas Research Fellowship 29-60. We are grateful for comments from participants of seminars at USC, the 2018 Institute of Mathematical Statistics Asia Pacific Rim Meeting in Singapore, the 2018 International Symposium on Financial Engineering and Risk Management in Shanghai, and the 2018 ICSA China Statistics Conference in Qingdao. %
}
} %
	
\date{\today}
\author{Yingying Fan$^1$, Jinchi Lv$^1$, Mahrad Sharifvaghefi$^1$ and Yoshimasa  Uematsu$^2$%
\\
University of Southern California$^1$ and Tohoku University$^2$\\
}
	
\maketitle
	
\begin{abstract}
Interpretability and stability are two important features that are desired in many contemporary big data applications arising in economics and finance. While the former is enjoyed to some extent by many existing forecasting approaches, the latter in the sense of controlling the fraction of wrongly discovered features which can enhance greatly the interpretability is still largely underdeveloped in the econometric settings. To this end, in this paper we exploit the general framework of model-X knockoffs introduced recently in Cand\`{e}s, Fan, Janson and Lv (2018), which is nonconventional for reproducible large-scale inference in that the framework is completely free of the use of p-values for significance testing, and suggest a new method of intertwined probabilistic factors decoupling (IPAD) for stable interpretable forecasting with knockoffs inference in high-dimensional models. The recipe of the method is constructing the knockoff variables by assuming a latent factor model that is exploited widely in economics and finance for the association structure of covariates. Our method and work are distinct from the existing literature in that we estimate the covariate distribution from data instead of assuming that it is known when constructing the knockoff variables, our procedure does not require any sample splitting, we provide theoretical justifications on the asymptotic false discovery rate control, and the theory for the power analysis is also established. Several simulation examples and the real data analysis further demonstrate that the newly suggested method has appealing finite-sample performance with desired interpretability and stability compared to some popularly used forecasting methods.
\end{abstract}

\textit{Running title}: IPAD


\section{Introduction} \label{sec:intro}
Forecasting is a fundamental problem that arises in economics and finance. With the availability of big data, many machine learning algorithms such as the Lasso and random forest can be resorted to for such a purpose by exploring a large pool of potential features. Many of these existing procedures provide a certain measure of feature importance which can then be utilized to judge the relative importance of selected features for the goal of interpretability. Yet the issue of stability in the sense of controlling the fraction of wrongly discovered features is still largely underdeveloped in the econometric settings. As argued in \cite{DeMol2008}, it is difficult to obtain interpretability and stability simultaneously even in simple Lasso forecasting. A natural question is how to ensure both interpretability and stability for flexible forecasting.

Naturally stability is related to statistical inference. The recent years have witnessed a growing body of work on high-dimensional inference in the econometrics and statistics literature. For example, \cite{WooldridgeZhu2018} proposed a simple procedure for inference of the average partial effects based on a debiased $\ell_1$-regularized method in approximately sparse panel probit models. 
\cite{StuckyGeer2018} used the de-sparsified estimator for constructing pointwise and group confidence sets. \cite{ZhangCheng2017} conducted simultaneous inference for high-dimensional sparse linear models based on a bootstrap and desparsifying Lasso estimator. They also applied their procedure for the family-wise error rate control. \cite{CCDDHNR2018} provided a double/debiased machine learning (DML) method for estimation and inference of treatment effects, which utilizes the Neyman orthogonal scores and cross-fitting. \cite{CNR2018} then extended this idea to linear functionals. 
\cite{CHHW2018} considered debiased simultaneous inference in a system of high-dimensional regression equations with temporal and cross-sectional dependency based on a uniform robust post-selection procedure.  \cite{ShahBuhlmann2018} proposed Lasso residual-based tests for checking goodness-of-fit in (low- and) high-dimensional linear models.  \cite{GuoEtAl2018} presented a method for estimating the effect of the treatment on the outcome by using instrumental variables where the instruments are not necessarily valid. 


Most existing work on high-dimensional inference for interpretable models has focused primarily on the aspects of post-selection inference known as selective inference and debiasing for regularization and machine learning methods. In real applications, one is often interested in conducting \textit{global} inference relative to the full model as opposed to \textit{local} inference conditional on the selected model. Moreover, many statistical inferences are based on p-values form significance testing. However, oftentimes obtaining valid p-values even for the Lasso in relatively complicated high-dimensional nonlinear models also remains largely unresolved, not to mention for the case of more complicated model fitting procedures such as random forest. Indeed high-dimensional inference is intrinsically challenging even in the parametric settings \cite{NUP}.

The desired property of stability for interpretable forecasting in this paper concentrates on \textit{global} inference by controlling precisely the fraction of wrongly discovered features in high-dimensional models, which is also known as reproducible large-scale inference. Such a problem involves testing the joint significance of a large number of features simultaneously, which is known widely as the problem of multiple testing in statistical inference. For this problem, the null hypothesis for each feature states that the feature is unimportant in the joint model which can be understood as the property that this individual feature and the response are \textit{independent} conditional on all the remaining features, while the corresponding alternative hypothesis states the opposite. Conventionally p-values from the hypothesis testing are used to decide whether or not to reject each null hypothesis with a significance level to control the probability of false discovery in a single hypothesis test, meaning rejecting the null hypothesis when it is true. When performing multiple hyothesis tests, the probability of making at least one false discovery which is known as the family-wise error rate can be inflated compared to that for the case of a single hypothesis test. The work on controlling such an error rate for multiple testing dates back to \cite{bonferroni1935calcolo}, where a simple, useful idea is lowering the significance level for each individual test as the target level divided by the total number of tests to be performed. The Bonferroni correction procedure is, however, well known to be conservative with relatively low power. Later on, \cite{holm1979simple} proposed a step-down procedure which is less conservative than the Bonferroni procedure. More recently, \cite{romano2005exact} suggested a procedure in which the critical values of individual tests are constructed sequentially. 

A more powerful and extremely popular approach to multiple testing is the Benjamini--Hochberg (BH) procedure for controlling the false discovery rate (FDR) which was originated in \cite{benjamini1995controlling}, where the FDR is defined as the expectation of the fraction of falsely rejected null hypotheses known as the false discovery proportion. Given the p-values from the multiple hypothesis tests, this procedure sorts the p-values from low to high and chooses a simple, intuitive cutoff point, which can be viewed as an adaptive extension of the Bonferroni correction for multiple comparisons, of the p-values for rejecting the null hypotheses. The BH procedure was shown to be capable of controlling the FDR at the desired level for independent test statistics in \cite{benjamini1995controlling} 
and for positive regression dependency among the test statistics in \cite{benjamini2001control}, where it was shown that a simple modification of the procedure can control the FDR under other forms of dependency but such a modification is generally conservative. There is a huge literature on the theory, applications, and various extensions of the original BH procedure for FDR control. See, for instance,  \cite{benjamini2010discovering} for a review of related developments, \cite{FanHanGu2012} for a factor model approach to FDR control under arbitrary covariance dependence, and  \cite{BCCHK2018} for a review of key results on estimation and inference including multiple testing with FDR control in high-dimensional models.

It is worth mentioning that \cite{OCMT} recently introduced a one covariate a time, multiple testing procedure for high-dimensional variable selection in linear regression models. In particular, their method was shown to have asymptotic FDR equal to the ratio of the number of pseudo signals and the total number of pseudo signals and true signals, where the true signals have nonzero regression coefficients and the pseudo signals have zero regression coefficients but nonzero marginal correlations with the response. Unlike \cite{OCMT}, the main interest of our paper is the FDR control with respect to only the set of true signals.

The aforementioned econometric and statistical inference methods including the BH-type procedures for FDR control are all rooted on the availability and validity of computable p-values for evaluating variable importance. As mentioned before, such a prerequisite can become a luxury that is largely unclear how to obtain in high dimensions even for the case of Lasso in general nonlinear models and random forest. 
In contrast, \cite{barber2015controlling} proposed a novel procedure named the knockoff filter for FDR control that bypasses the use of p-values in Gaussian linear model with deterministic design matrix, where the dimensionality is no larger than the sample size, and \cite{barber2016knockoff} generalized the method to high-dimensional linear models as a two-step procedure based on sample splitting, where a feature screening approach is used to reduce the dimensionality to below sample size (see, e.g., \cite{FanFan2008} and \cite{FanLv2008}) and then the knockoff filter is applied to the set of selected features after the screening step for selective inference. The key ingredient of the  knockoff filter is constructing the so-called knockoff variables in a geometrical way that mimic perfectly the correlation structure among the original covariates and can be used as control variables to evaluate the importance of original variables. Recently, \cite{candes2016panning} extended the work of \cite{barber2015controlling} by introducing the framework of model-X knockoffs for FDR control in general high-dimensional nonlinear models. A crucial distinction is that the knockoff variables are constructed in a probabilistic fashion such that the joint dependency structure of the original variables and their knockoff copies is invariant to the swapping of any set of original variables and their knockoff counterparts, which enables us to go beyond linear models and handle high dimensionality. As a result, model-X knockoffs enjoys exact finite-sample FDR control at the target level. However, a major assumption in \cite{candes2016panning} is that the joint distribution of all the covariates needs to be \textit{known} for the valid FDR control.

Motivated by applications in economics and finance, in this paper we model the association structure of the covariates using the latent factor model, which reduces effectively the dimensionality and enables reliable estimation of the \textit{unknown} joint distribution of all the covariates. By taking into account the latent factor model structure, we first estimate the association structure of covariates and then construct \textit{empirical} knockoffs matrix using the estimated dependency structure.  Our empirical knockoffs matrix can be regarded as an  approximation to the \textit{oracle} knockoffs matrix in \cite{candes2016panning} that requires the knowledge of the true covariate distribution. Exploiting the general framework of model-X knockoffs in \cite{candes2016panning}, we suggest the new method of intertwined probabilistic factors decoupling (IPAD) for stable interpretable forecasting with knockoffs inference in high-dimensional models. The innovations of our method and work are fourfold. First, we estimate the covariate distribution from data instead of assuming that it is known when constructing the knockoff variables. Second, our procedure does not require any sample splitting and is thus more practical when the sample size is limited. Third, we provide theoretical justifications on the asymptotic false discovery rate control when the estimated dependency structure is employed. Fourth, the theory for power analysis is also established which reveals that there can be asymptotically no power loss   in applying the knockoffs procedure compared to the underlying variable selection method. Therefore, FDR control by knockoffs can be a pure gain. Compared to earlier work, an additional challenge of our study is that knowing the true underlying distribution does \textit{not} lead to the most efficient construction of the oracle knockoffs matrix due to the presence of latent factors. The appealing interpretability and stability of our new method compared to some popularly used forecasting methods are confirmed with several simulation and real data examples.

The rest of the paper is organized as follows. Section \ref{sec:metho} introduces the model setting with a review of the model-X knockoffs inference framework and presents the new IPAD procedure. We establish the asymptotic properties of IPAD in Section \ref{sec:theory}. Sections \ref{sec:simu} and \ref{sec:rdata} present several simulation and real data examples to showcase the finite-sample performance and the advantages of our newly suggested procedure compared to some popularly used ones. We discuss some implications and extensions of our work in Section \ref{Sec6}. The proofs of the main results and additional technical details are relegated to the Appendix.

\section{Intertwined probabilistic factors decoupling} \label{sec:metho} 
To facilitate the technical presentation, we will introduce the model setting for the high-dimensional FDR control problem in Section \ref{sec:model} with a review of the model-X knockoffs inference framework in Section \ref{subsec:background}, and present the new IPAD procedure in Section \ref{subsec:newmethod}.

\subsection{Model setting} \label{sec:model}
Consider the high-dimensional linear regression model 
\begin{align}\label{model1}
\by = \bX\bbeta + \bep,
\end{align}
where $\by\in\mathbb{R}^n$ is the response vector, $\bX\in\mathbb{R}^{n\times p}$ is the random matrix of a large number of potential regressors, $\bbeta = (\beta_1,\cdots, \beta_p)'\in \mathbb{R}^p$ is the regression coefficient vector, $\bep \in\mathbb{R}^n$ is the vector of model errors, and $n$ and $p$ denote the sample size and dimensionality, respectively. Here without loss of generality, we assume that both the response and the covariates are centered with mean zero and thus there is no intercept. Motivated by many applications in economics and finance, we further assume that the design matrix $\bX$ follows the \textit{exact} factor model 
\begin{align}\label{factmodel}
\bX = \bF^0{\bLambda^0}' +\bE = \bC^0 + \bE,
\end{align}
where $\bF^0=(\bff_1^0,\dots,\bff_n^0)'\in\mathbb{R}^{n\times r}$ is a random matrix of latent factors, $\bLambda^0=(\blambda_1^0,\dots,\blambda_p^0)'\in\mathbb{R}^{p\times r}$ is a matrix of deterministic factor loadings, and $\bE\in\mathbb{R}^{n\times p}$ captures the remaining variation that cannot be explained by these latent factors. We assume that the number of factors $r$ is fixed but \textit{unknown} and  the components of $\bE$ are 
independent and identically distributed (i.i.d.) from some unknown parametric distribution with cumulative distribution function $G(\cdot; \bbbeta^0)$, where  $\bbbeta^0 \in \mathbb{R}^m$ is a finite-dimensional parameter vector. 
For simplicity, models \eqref{model1} and \eqref{factmodel} are assumed to have no endogeneity. 

In this paper, we focus on the high-dimensional scenario when the dimensionality $p$ can be much larger than sample size $n$. Therefore, to ensure model identifiability  we impose the sparsity assumption that the true regression coefficient vector $\bbeta$ has only a small portion of nonzeros; specifically,  $\bbeta$ takes nonzero values only on some (unknown) index set $\cS^0\subset \{1,\dots,p\}$ and  $\beta_j = 0$ for all $j \in \cS^1:=\{1,\dots,p\}\backslash \cS^0$. Denote by $s = |\cS^0|$ the size of $\cS^0$. We assume that $s = o(n)$ throughout the paper.  

We are interested in identifying the index set $\cS^0$ with a theoretically guaranteed error rate.  
To be more precise, we try to select variables in $\cS^0$  while keeping the false discovery rate (FDR)  under some prespecified desired level $q\in(0,1)$, where the FDR is defined as
\begin{align}\label{def:FDR}
\text{FDR} := \E \left[ \text{FDP} \right] \quad \text{with} \quad \text{FDP} := \frac{|\what{\cS}\cap \cS^1|}{|\what{\cS}| \vee 1}.
\end{align}
Here the FDP stands for the false discovery proportion and $\widehat \cS$ represents the set of variables selected by some procedure using observed data $(\bX, \by)$. A slightly modified version of FDR is defined as
\begin{align}\label{def:mFDR}
\text{mFDR} := \E\left[\frac{|\what{\cS} \cap \cS^1|}{|\what{\cS}| + q^{-1}} \right].
\end{align}
Clearly, FDR is more conservative than mFDR in that the latter is always under control if the former is.

It is easy to see that FDR is a measurement of type I error for variable selection. The other important aspect of variable selection is power, which is defined as 
\begin{align}\label{eq: power}
\text{Power} := \E \left[ \frac{|\what{\cS} \cap \cS^0|}{|\cS^0|} \right] =  \E \left[ \frac{|\what{\cS} \cap \cS^0|}{s} \right].
\end{align}
It is well known that FDR and power are two sides of the same coin. We aim at developing a variable selection procedure with theoretically guaranteed  FDR control and meanwhile achieving high power. 

\subsection{Review of model-X knockoffs framework} \label{subsec:background}  
The key idea of the model-X knockoffs framework is to construct the so-called model-X knockoff variables, which were introduced originally in \cite{candes2016panning} and whose definition  is stated formally as follows for completeness.

\begin{defi}[Model-X knockoff variables \cite{candes2016panning}]\label{def:knockoff}\normalfont
	For a set of random variables $ \bx = (X_1,\dots,X_p) $, a new set of random variables $ \wtilde{\bx} = (\wtilde{X}_1,\cdots,\wtilde{X}_p) $ is called a set of model-X knockoff variables if it satisfies the following properties:
	\begin{enumerate}
		\item[1)] For any subset $ \cS \subset \{1,\dots,p\} $, we have 
		$[\bx,\wtilde{\bx}]_{\textsf{swap}(\cS)} = [\bx,\wtilde{\bx}]$ in distribution, 
		where the vector $ [\bx,\wtilde{\bx}]_{\textsf{swap}(\cS)} $ is obtained by swapping $ X_j $ and $ \wtilde{X}_j $ for each $ j \in \cS $. 
		
		\item[2)] Conditional on $ \bx $, the knockoffs vector $ \wtilde{\bx} $ is independent of response $Y$.
	\end{enumerate}
\end{defi} 
An important consequence is that the null regressors $\{X_j:j \in \cS^1\}$ can be swapped with their knockoffs without changing the joint distribution of the original variables $ \bx$, their knockoffs $\wtilde{\bx}$, and response $ Y $. That is, we can obtain for any $\cS \subset \mathcal{S}^1$,
\begin{align} \label{exchange}
([\bx,\wtilde{\bx}]_{\textsf{swap}(\cS)},Y) \overset{d}{=} ([\bx,\wtilde{\bx}],Y),
\end{align}
where $\overset{d}{=}$ denotes equal in distribution. 
Such a property is known as the \textit{exchangeability property} using the terminology in \cite{candes2016panning}. For more details, see Lemma 3.2 therein. Following \cite{candes2016panning}, one can obtain a knockoffs matrix $\wtilde{\bX} \in \mathbb{R}^{n \times p}$ given observed design matrix $\bX$. 

Using the augmented design matrix $ [\bX,\wtilde{\bX}] $ and response vector $ \by $ constructed by  stacking the $n$ observations, \cite{candes2016panning} suggested constructing knockoff statistics $W_j = w_j([\bX,\wtilde{\bX}],\by)$, $j\in \{1,\dots, p\}$, for measuring the importance of the $j$th variable, where $w_j$ is some function that  satisfies the property that swapping $\bx_j \in \mathbb{R}^n$ with its corresponding knockoff variable $\wtilde{\bx}_j \in \mathbb{R}^n$ changes the sign of $ W_j $; that is,  
\begin{equation}\label{flipsign}
w_j([\bX,\wtilde{\bX}]_{\textsf{swap}(\cS)},\by) = 
\begin{cases}
w_j([\bX,\wtilde{\bX}],\by), & j \notin \cS, \\
-w_j([\bX,\wtilde{\bX}],\by), & j \in \cS.
\end{cases}
\end{equation}
The  knockoff statistics constructed above $W_j=w_j([\bX,\wtilde{\bX}],\by)$ satisfy the so-called sign-flip property; that is, 
conditional on $|W_j|$'s the signs of the null $W_j$'s with $j\not\in \cS^0$ are i.i.d.\ coin flips (with equal chance $1/2$). For the examples on valid constructions of knockoff statistics, see  \cite{candes2016panning}. 

Let $ t>0 $ be a fixed threshold and define $\what{\cS}=\{j: W_j \geq t \}$ as the set of discovered variables. 
Then intuitively, the sign-flip property entails
\begin{align*}
\left|\widehat{\cS}\cap\cS^1\right| \overset{d}{=} \left|\{j:W_j \leq -t\}\cap\cS^1\right| \leq \left|\{j:W_j \leq -t\}\right|.
\end{align*}
Therefore, the FDP function can be estimated (conservatively) as 
\begin{align*}
\text{FDP} = \frac{|\widehat{\cS}\cap \cS^1|}{|\widehat{\cS}|\vee 1} 
\leq \frac{\left|\{j:W_j \leq -t\}\right|}{|\widehat{\cS}|\vee 1} =: \widehat{\text{FDP}}
\end{align*}
for each $t$. 
In light of this observation,  \cite{candes2016panning} proposed to choose the threshold by resorting to the above $\widehat{\text{FDP}}$. Their results are summarized formally as follows. 

\begin{res}[\cite{candes2016panning}] \label{thm:candes}
	Let $ q \in(0,1)$ denote the target FDR level. Assume that we choose a threshold $ T_1 > 0 $ such that
	\begin{equation*}
	T_1= \min \left\{ t > 0:\frac{|\{j:W_j \le -t \}|}{| \{ j: W_j \ge t\}| \vee 1} \le q \right\}
	\end{equation*}
	or $ T_1=+\infty $ if the set is empty. Then the procedure selecting the variables $\what{\cS} = \{j:W_j \geq T_1\}$ controls the mFDR in \eqref{def:mFDR} to no larger than $q$. 
	Moreover, assume that we choose a slightly more conservative threshold $ T_2 > 0 $ such that
	\begin{equation*}
	T_2= \min \left\{ t > 0:\frac{1+|\{j:W_j \le-t \}|}{| \{ j: W_j \ge t\}| \vee 1} \le q \right\}
	\end{equation*}
	or $ T_2=+\infty $ if the set is empty. Then the procedure selecting the variables $ \what{\cS} = \{j:W_j \geq T_2\} $ controls the FDR in \eqref{def:FDR} to no larger than $q$.
\end{res}

It is worth noting that Result \ref{thm:candes} was derived under the assumption that the joint distribution of the $p$ covariates is known. In our  model setting \eqref{model1} and \eqref{factmodel}, however there exist unknown parameters that need to be estimated from data. In such case, it is natural to construct the knockoff variables and knockoff statistics with estimated distribution of the $p$ covariates. Such a plug-in principle usually leads to breakdown of the exchangeability property in Definition \ref{def:knockoff}, preventing us from using directly Result \ref{thm:candes}. To address this challenging issue, we will introduce our new method in the next section and provide detailed theoretical analysis for it. 

It is also worth mentioning that recently, \cite{barber2018Rknockoff} provided an elegant new line of theory which ensures FDR control of model-X knockoffs procedure under the approximate exchangeability assumption, which is weaker than the exact exchangeability condition required in Definition \ref{def:knockoff}. However, the conditions they need on estimation error of the joint distribution of $\bx$ is difficult to be satisfied in high dimensions.  \cite{RANK} investigated the robustness of model-X knockoffs procedure with respect to unknown covariate distribution when covariates $\bx$ follow a joint Gaussian distribution. Their procedure needs data splitting and their proofs rely heavily on the Gaussian distribution assumption, and thus their development may not be suitable for economic data with limited sample size and heavy-tailed distribution.   For these reasons, our results complement substantially those in \cite{candes2016panning}, \cite{RANK}, and \cite{barber2018Rknockoff}.   

\subsection{IPAD} \label{subsec:newmethod}
It has been seen from the previous section that the key for the model-X knockoffs framework is the construction of valid knockoff variables.  We begin with the ideal situation where the the factor model structure \eqref{factmodel} is fully available to us; that is, we know the realization $\bC^0$ and the distribution $G(\cdot; \bbbeta^0)$ for the error matrix $\bE$.  In such case,
the oracle knockoffs matrix $\wtilde{\bX}(\btheta^0)$ can be constructed as 
\begin{align}\label{eq: ideal-knockoff}
\wtilde{\bX}(\btheta^0) = \bC^0 + \bE_{\bbbeta^0},
\end{align}
where $\bE_{\bbbeta^0}$ is an i.i.d.\ copy of $\bE$ and $\btheta^0 = (\bC^0, \bbbeta^0)$ is the augmented parameter vector. Note that $\bE_{\bbbeta_0}$ itself is not a function of $\bbbeta_0$, but we slightly abuse the notation to emphasize the dependence of the distribution function on parameter $\bbbeta_0$.  It is easy to check that $\wtilde{\bX}(\btheta^0)$ constructed above is a valid knockoffs matrix and satisfies the properties in Definition \ref{def:knockoff}. Although $\wtilde{\bX}(\btheta^0)$ is generally unavailable to us, it plays an important role in our theoretical developments. 

We remark that in the construction above, we slightly misuse the concept and call $\bC^0$ a parameter. This is because although $\bC^0$ is a random matrix, for the construction of valid knockoff variables it is the particular realization $\bC^0$ leading to the observed data matrix $\bX$ that matters.  In other words, a valid construction of knockoff variables requires the knowledge of the specific realization $\bC^0$  instead of the distribution of $\bC^0$.   To understand this, consider the scenario where the underlying parameter $\bbbeta^0$ and the exact distribution of $\bC^0$ are fully available to us. If we independently generate random variables from this known distribution and form a new data matrix $\bX_1$,  because of the independence between $\bX_1$ and $\bX$, the exchangeability assumption in Definition \ref{def:knockoff} will be violated and thus $\bX_1$ cannot be a valid knockoffs matrix.   On the other hand, as long as we know the realization $\bC^0$ and parameter $\bbbeta^0$,  a valid knockoffs matrix $\widetilde{\bX}(\btheta^0)$ can be constructed using \eqref{eq: ideal-knockoff} regardless of whether the exact distribution of $\bC^0$ is available to us or not.    

In practice, however  $\btheta^0$  is unavailable to us and consequently, $\widetilde{\bX}(\btheta^0)$ is inaccessible. To overcome this difficulty, we next introduce our new method IPAD.  We start with introducing the \textit{knockoff generating function} -- for each given parameter vector $\btheta = (\bC, \bbbeta)$, define
\begin{align}\label{funcXk}
\wtilde{\bX}(\btheta) = \bC + \bE_{\bbbeta},
\end{align}
where $\bE_{\bbbeta}$ is a matrix composed of i.i.d.\ random samples from the distribution  $G(\cdot;\bbbeta)$. 
Letting $\hat{\btheta}$ denote an estimator  (obtained using data $\bX$) of $\btheta^0$, we name $\wtilde{\bX}(\hat{\btheta})$ as the \textit{empirical} knockoffs matrix while $\wtilde{\bX}(\btheta^0)$ as the \textit{oracle (ideal)} knockoffs matrix. 

With the aid of empirical knockoffs matrix, we suggest the following IPAD procedure for FDR control with knockoffs inference. 
\begin{proc}[IPAD]\label{proc:1}\normalfont
	\begin{enumerate}
		\item[1)] (Estimation of parameters) Estimate the unknown parameters in $\btheta^0$ using the design matrix $\bX$. Denote by $\hat\btheta = (\widehat\bC, \hat\bbbeta)$ the resulting estimated parameter vector.   
		
		\item[2)] (Construction of empirical knockoffs matrix) Construct the empirical knockoffs matrix by applying the knockoff generating function in \eqref{funcXk} to the estimated parameter $\hat\btheta$; that is,   
		\begin{align}\label{Xk}
		\wtilde{\bX}(\hat{\btheta}) = \what{\bC} + \bE_{\hat{\bbbeta}},
		\end{align}
		where $\bE_{\hat{\bbbeta}} \in \mathbb{R}^{n\times p}$ is a matrix composed of i.i.d.\ random variables from $G(\cdot;\hat{\bbbeta})$, and is independent of $(\bX, \by)$ conditional on $\hat\bbbeta$. 
		\item[3)] (Application of knockoffs inference) Calculate  knockoff statistics  $W_j(\hat{\btheta})$  using data $([\bX,\wtilde{\bX}(\hat{\btheta})], \by)$ and then construct  $\what{\cS}$ by applying knockoffs inference to $W_j(\hat{\btheta})$. 
	\end{enumerate}
\end{proc}

Intuitively, the accuracy of the estimator $\hat\btheta$ in Step 1 will affect the performance of our IPAD procedure.  In fact, as shown later in our Theorem \ref{thm:FDRcontrol}, the consistency rate of $\hat\btheta$ is indeed reflected in the asymptotic  FDR control of the IPAD procedure. For the specific case when the error distribution is $N(0,\sigma^2)$, the parameter $\sigma^2$ can be estimated naturally as  $(np)^{-1}\sum_{i,j}\hat{e}_{ij}^2$, where $\hat{e}_{ij}$'s are the entries of $\what{\bE} = \bX - \what{\bC}$. 
In Step 3, various methods can be used to construct knockoff statistics. For the illustration purpose, we use the  \textit{Lasso coefficient difference} (LCD) statistic as in \cite{candes2016panning}. Specifically,
with $\by$ the response vector and $ ([\bX, \wtilde\bX(\hat\btheta)])$ the augmented design matrix we consider the variable selection procedure Lasso \cite{Tibshirani1996} which solves the following optimization problem 
\begin{equation}\label{ex:lasso}
\hat{\bbeta}^\textsf{aug}(\hat\btheta;\lambda) = \arg\min_{\bb\in\mathbb{R}^{2p}}\left\{\| \by - [\bX, \widetilde{\bX}(\hat\btheta)] \bb\|_2^2 + \lambda \|\bb\|_1\right\}, 
\end{equation}
where $\lambda\geq 0$ is the regularization parameter and $\|\cdot\|_m$ with $m \geq 1$ denotes the vector $\ell_m$-norm. 
Then for each variable $ \bx_j $, the knockoff statistic can be constructed as 
\begin{equation}\label{def: LCD}
W_j(\hat\btheta; \lambda)=|\hat{\beta}_{j}^\textsf{aug}(\hat\btheta;\lambda)|-|\hat{\beta}_{p+j}^\textsf{aug}(\hat\btheta;\lambda)|,
\end{equation} 
where $\hat{\beta}_{\ell}^\textsf{aug}(\hat\btheta;\lambda)$ is the $\ell$th component of the Lasso regression coefficient vector $\hat{\bbeta}^\textsf{aug}(\hat\btheta;\lambda)$. 
It is seen that intuitively the LCD knockoff statistics evaluate the relative importance of the $j$th original variable by comparing its Lasso coefficient $\hat{\beta}_{j}^\textsf{aug}(\hat\btheta; \lambda)$  with that of its knockoff copy $\hat{\beta}_{j+p}^\textsf{aug}(\hat\btheta; \lambda)$.  
In the ideal case when the oracle knockoffs matrix $\widetilde{\bX}(\btheta^0)$ is used instead of $\widetilde{\bX}(\hat\btheta)$ in \eqref{ex:lasso}, it is easy to verify that the LCD is a valid construction of knockoff statistics and satisfies the sign-flip property in \eqref{flipsign}. Consequently, the general theory in \cite{candes2016panning} can be applied to show that the FDR is controlled in finite sample. We next show that even with the empirical knockoffs matrix employed in \eqref{ex:lasso}, the FDR can still be asymptotically controlled with delicate technical analyses.

\section{Asymptotic properties of IPAD} \label{sec:theory}
We now provide theoretical justifications for our IPAD procedure suggested in Section \ref{sec:metho} with the LCD knockoff statistics  $W_j(\hat\btheta; {\lambda})=w_j([\bX,\wtilde{\bX}(\hat\btheta)],\by; {\lambda})$ defined in \eqref{def: LCD}. 
We will first present some technical conditions in Section \ref{Sec3.1}, then prove in Section \ref{subsec:FDR} that the FDR is asymptotically under control at desired target level $q$, and finally in Section \ref{subsec:Power} show that asymptotically IPAD has no power loss compared to the Lasso under some regularity conditions. 

\subsection{Technical conditions} \label{Sec3.1}
We first introduce some notation and definitions which will be used later on.  We use $X\sim \mbox{subG}(C_x^2)$ to denote that $X$ is a sub-Gaussian  random variable with \textit{variance proxy} $C_x^2>0$ if $\E [X]=0$ and its tail probability satisfies $\Pro(|X|>u) \leq 2\exp(u^2/C_x^2)$ for each $u\geq 0$. 
In all technical assumptions below, we use $M>1$ to denote a large enough generic constant. Throughout the paper, for any vector $\bv = (v_i)$ let us denote by $\|\bv\|_1$, $\|\bv\|_2$, and $\|\bv\|_{\max}$ the $\ell_1$-norm, $\ell_2$-norm, and max-norm defined as $\|\bv\|_1=\sum_{i}|v_{i}|$, $\|\bv\|_2=(\sum_iv_i^2)^{1/2}$, and  $\|\bv\|_{\max}=\max_{i}|v_{i}|$, respectively. 
For any matrix $\bM=(m_{ij})$, we denote by $\|\bM\|_F$, $\|\bM\|_1$, $\|\bM\|_2$, and $\|\bM\|_{\max}$ the Frobenius norm, entrywise $\ell_1$-norm, spectral norm, and  entrywise $\ell_{\infty}$-norm defined as  $\|\bM\|_F=\|\vect(\bM)\|_2$, $\|\bM\|_1=\|\vect(\bM)\|_1$, $\|\bM\|_2=\sup_{\bv\not= \bzero}\|\bM\bv\|_2/\|\bv\|_2$, and  $\|\bM\|_{\max}=\|\vect(\bM)\|_{\max}$, respectively, where $\vect(\bM)$ represents the vectorization of matrix $\bM$. For a symmetric matrix $\bM$, $\vech(\bM)$ stands for the vectorization of the lower triangular part of $\bM$.

\begin{con}[Regression errors]\normalfont \label{ass:regerr}
	The model error vector $\bep$ has i.i.d.\ components from subG($C_\ep^2$).
\end{con}

\begin{con}[Latent factors]\normalfont \label{ass:fac}
	The rows of $\bF^0$ consist of mean zero i.i.d.\ random vectors $\bff_i^0\in\mathbb{R}^r$ such that $\|\bF^0\|_{\max} \leq M$ almost surely (a.s.) and $\|\bSigma_f\|_2+\|\bSigma_f^{-1}\|_2 \leq M$, where $\bSigma_f:=\E[\bff_i^0{\bff_i^0}']$. 
\end{con}

\begin{con}[Factor loadings]\normalfont \label{ass:floa}
	The rows of $\bLambda^0$ consist of deterministic vectors $\blambda_j^0\in\mathbb{R}^r$ such that $\|\bLambda^0\|_{\max} \leq M$ and $\|p^{-1}{\bLambda^0}'\bLambda^0\|_2+\|(p^{-1}{\bLambda^0}'\bLambda^0)^{-1}\|_2 \leq M$. 
\end{con}

\begin{con}[Factor errors]\normalfont \label{ass:err_new}
	The entries of matrix $\bE_{\bbbeta^0}$ are i.i.d.\ copies of $e_{\bbbeta^0}\sim \mbox{subG}(C_e^2)$ with continuous distribution function $G(\cdot;\bbbeta^0)$. For  each $1\leq \ell \leq m$, the $\ell$th element of $\bbbeta^0$  is specified as
	$\eta_\ell^0 = h_\ell (\E[e_{\bbbeta^0}],\dots, \E[e_{\bbbeta^0}^{m}])$ with $h_\ell:\mathbb{R}^{m}\to\mathbb{R}$ some local Lipschitz continuous function in the sense that 
	\begin{align*}
	\left|h_\ell(t_1,\dots,t_{m})-h_\ell(\E[e_{\bbbeta^0}],\dots, \E[e_{\bbbeta^0}^{m}])\right| \leq M\max_{k\in\{1,\dots,m\}}\left|t_k-\E[e_{\bbbeta^0}^k]\right|
	\end{align*}
	for each $t_k \in \{t: |t-\E[e_{\bbbeta^0}^k]| \leq Mc_{np} \}$ and $1 \leq k \leq m$, where $c_{np}:= (p^{-1}\log n)^{1/2} + (n^{-1}\log p)^{1/2}$. 
	Moreover, there exists some stochastic process $(e_{\bbbeta})_{\bbbeta}$ such that 
	\begin{enumerate}
		\item[i)]   for each $\bbbeta \in \{\bbbeta \in \mathbb{R}^m: \|\bbbeta-\bbbeta^0\|_{\max}\leq Mc_{np}\}$, the entries of $\bE_{\bbbeta}$  in \eqref{funcXk} have identical distribution to $e_{\bbbeta}$,
		\item[ii)] 	for some sub-Gaussian random variable $Z \sim \mbox{subG}\left(c_e^2\right)$ with some positive constant $c_e$, \begin{align}\label{ineq:subGprocess}
		\sup_{\bbbeta :\, \|\bbbeta-\bbbeta^0\|_{\max}\leq Mc_{np} }|e_{\bbbeta}-e_{\bbbeta^0}| \leq  M^{1/2}c_{np}^{1/2}|Z|.
		\end{align}	 
	\end{enumerate}
\end{con}

\begin{con}[Eigenseparation]\normalfont \label{ass:eigensepa}
	The $r$ eigenvalues of $p^{-1}{\bLambda^0}'\bLambda^0 \bSigma_f$  are distinct for all $p$. 
\end{con}

The number of factors $r$ is assumed to be known for developing the theory with simplification, but in practice it can be estimated consistently using methods such as information criteria \citep{bai2002determining} and test statistics  \citep{ahn2013}. The sub-Gaussian assumptions in Conditions \ref{ass:regerr} and \ref{ass:err_new}  can be replaced with some other tail conditions as long as similar concentration inequalities hold. Condition \ref{ass:floa} is standard in the analysis of factor models. Stochastic loadings can be assumed in Condition \ref{ass:floa} with some appropriate distributional assumption, such as sub-Gaussianity, at the cost of much more tedious technical arguments. The boundedness of the eigenvalues of $\bSigma_f$ in Condition \ref{ass:fac} is standard while the i.i.d.\ assumption and boundedness of $\bff_i^0$ are stronger compared to the existing literature (e.g., \cite{bai2002determining} and \cite{bai2003}). However, these conditions are imposed mostly for technical simplicity. In fact, the boundedness condition on $\bff_i^0$ can be replaced with (unbounded) sub-Gaussian or other heavier-tail assumption whenever concentration inequalities are available at the cost of slower convergence rates and stronger sample size requirement. Our theory on FDR control is based on that in \cite{candes2016panning}, which applies only to the case of i.i.d.\ rows of design matrix $\bX$. This is the main reason for imposing the  i.i.d.\ assumption on $\ep_i$ and $\bff_i$ in Conditions \ref{ass:regerr} and \ref{ass:fac}. However, we conjecture
that similar results can also hold  in the presence of some sufficiently weak serial dependence in $\ep_i$ and $\bff_i$. Condition \ref{ass:err_new} introduces a \textit{sub-Gaussian process}  $e_{\bbbeta}$ with respect to $\bbbeta$. The norm in \eqref{ineq:subGprocess} can be replaced with any other norm since $\bbbeta$ is finite dimensional. In the specific case when the components of  $\bE$ have Gaussian distribution such that $\bbbeta$ is a scalar parameter representing variance, by the the reflection principle for the Wiener process (\cite{Billingsley1995}, p.511),  $e_{\bbbeta}$  can be constructed as a Wiener process and the  inequality \eqref{ineq:subGprocess} can be satisfied. For more information on sub-Gaussian processes, see, e.g., \cite{Vizcarra2007}. Condition \ref{ass:eigensepa} guarantees that $\hat{\bF}'\bF^0/n$ is asymptotically nonsingular, which has been proved in \cite{bai2003} and is used in the proof of Lemma \ref{lem:HVbound} in Appendix. 

Recall that in the IPAD procedure, we first obtain the augmented Lasso estimator $\hat{\bbeta}^\textsf{aug}(\btheta; \lambda)\in\mathbb{R}^{2p}$ by regressing  $\by$ on $[\bX,\wtilde{\bX}(\btheta)]$. 
Denote by  $\cA^\textsf{aug}(\btheta; \lambda)=\supp ( \hat{\bbeta}^\textsf{aug}(\btheta; \lambda) ) \subset \{1,\dots,2p\}$ the active set of the augmented Lasso regression coefficient vector. Throughout this section, we content ourselves with sparse estimates satisfying
\begin{align}\label{def: k}
\left|	\cA^\textsf{aug}(\btheta; \lambda)\right| \leq k/2
\end{align} 
for some positive integer $k$ which may diverge with $n$ at an order to be specified later; see, e.g., \cite{FanLv2013} and \cite{Lv2013} for a similar constraint and justifications therein. This can always be achieved since users have the freedom to choose the size of the Lasso  model. 

\subsection{FDR control} \label{subsec:FDR} 
To develop the theory for IPAD, we consider the principle component estimator $\what\bC$ for the realization $\bC^0$.  More specifically, we first conduct the singular value decomposition (SVD)  $ \bX = \bU\bD\bV'$ with $\bU$ and $\bV$ the left and right singular matrices and $\bD$ a diagonal matrix of singular values, and then threshold the diagonal matrix $\bD$ by setting the smallest $n- r$ singular values to zero. Let us denote the thresholded matrix as $\bD_{r}$. Then matrix $\bC^0$ can be estimated as $\what\bC = \bU\bD_{r}\bV'$. Denote by $\what \bE = (\hat{ e}_{ij}) = \bX - \what\bC$. 
The estimator $\hat\bbbeta = (\hat\eta_1,\cdots, \hat\eta_m)'$ is constructed as $\hat\eta_\ell=h_\ell(\E_{np}\hat{e},\dots,\E_{np}\hat{e}^{m})$ with $h_\ell$, $1 \leq \ell \leq m$, introduced in Condition \ref{ass:err_new} and $\mathbb{E}_{np}\hat{e}^k=(np)^{-1}\sum_{i,j} \hat{e}_{ij}^k$ the empirical moments of $\hat{e}_{ij}$. Throughout our theoretical analysis, we consider the regularization parameter fixed at $\lambda = C_0n^{-1/2}\log p$ with $C_0$ some large enough constant for all the Lasso procedures. Therefore, we will drop the dependence of various quantities on $\lambda$ whenever there is no confusion. For example, we will write $\cA^\textsf{aug}(\btheta; \lambda)$ and $ \hat{\bbeta}^\textsf{aug}(\btheta; \lambda)$ as $\cA^\textsf{aug}(\btheta)$ and $ \hat{\bbeta}^\textsf{aug}(\btheta)$, respectively. 

Denote by $\bU(\btheta) := n^{-1}[\bX, \wtilde{\bX}(\btheta)]'[\bX, \wtilde{\bX}(\btheta)]$ and $\bv(\btheta) := n^{-1}[\bX, \wtilde{\bX}(\btheta)]'\by$ and define 
$\bT(\btheta):=\vect(\vech\bU(\btheta), \bv(\btheta))\in\mathbb{R}^{P}$ with $P:=p(2p+3)$. The following lemma states that the statistic $\bT(\btheta)$ plays a crucial role in our procedure. 

\begin{lem}\label{lem:rep}
	The set of variables $\what\cS$ selected by Procedure \ref{proc:1} depends only on $\bT(\btheta)$. 
\end{lem}

For any given $\btheta$, define the active set $\cA^*(\btheta) := \cA_{1}^\textsf{aug}(\btheta)\cup \cA_{2}^\textsf{aug}(\btheta)\subset \{1,\dots,p\}$, where $\cA_{1}^\textsf{aug}(\btheta):=\{j:j\in\{1,\dots,p\}\cap \cA^\textsf{aug}(\btheta) \}$ and $\cA_{2}^\textsf{aug}(\btheta):=\{j-p:j\in\{p+1,\dots,2p\}\cap \cA^\textsf{aug}(\btheta) \}$.  That is, $\cA^*(\btheta)$ is equal to the support of knockoff statistics $(W_1(\btheta), \cdots, W_p(\btheta))'$ if there are no ties on the magnitudes of the augmented Lasso coefficient vector $\hat\bbeta^{\textsf{aug}}(\btheta)$.  

We next focus on the low-dimensional structure of $\bT(\btheta)$ inherited from the augmented Lasso because it will be made clear that this is the key to controlling the FDR \textit{without} sample splitting.
For any subset $\cA\subset\{1,\dots,p\}$,  define a lower-dimensional expression of the vector as $\bT_{\cA}({\btheta}) := \vect(\vech\bU_{\cA}({\btheta}),\bv_{\cA}({\btheta}))$ with $\bU_{\cA}({\btheta})$ the principle submatrix of $\bU(\btheta)$ formed by columns and rows in $\cA$ and $\bv_{\cA}({\btheta})$ the subvector of $\bv(\btheta)$ formed by components in $\cA$. Then it is easy to see that $\bU_{\cA}(\btheta) = n^{-1}[\bX_{\cA},\wtilde{\bX}_{\cA}({\btheta})]'[\bX_{\cA},\wtilde{\bX}_{\cA}({\btheta})]$ and $\bv_{\cA}(\btheta)  = n^{-1}[\bX_{\cA},\wtilde{\bX}_{\cA}({\btheta})]'\by $. Motivated by Lemma \ref{lem:rep}, we define a family of mappings indexed by $\cA$ that describes the \textit{selection algorithm} of Procedure \ref{proc:1} with given data set $([\bX_{\cA},\wtilde{\bX}_{\cA}({\btheta})],\by)$ that forms $\bT_{\cA}(\btheta)$. Formally, define a mapping $S_\cA:\mathbb{R}^{|\cA|(2|\cA|+3)} \to 2^\cA$ as $\bt_\cA \mapsto S_\cA(\bt_\cA)$ for given $\bT_\cA(\btheta)=\bt_\cA$, where $2^\cA$ refers to the power set of $\cA$. That is,   $S_\cA(\bt_\cA)$ represents the outcome of first restricting ourselves to the smaller set of variables $\cA$ and then applying IPAD to $\bT_{\cA}(\btheta)=\bt_{\cA}$ to further select variables from set $\cA$. 

\begin{lem}\label{lem:lowdimrep}
	Under Conditions \ref{ass:regerr}--\ref{ass:err_new}, 
	for any subset $\cA \supset\cA^*(\btheta)$ we have $S_{\{1,\dots,p\}}(\bT(\btheta))=S_\cA(\bT_\cA(\btheta))$.
\end{lem}

When restricting on set $\cA$, we can apply Procedure \ref{proc:1} to a lower-dimensional data set $([\bX_{\cA},\wtilde{\bX}_{\cA}({\btheta})],\by)$ that forms $\bT_{\cA}(\btheta)$ to further select variables from $\cA$.  The previous two lemmas ensure that this gives us a subset of $\cA$ that is identical to $S_{\{1,\dots,p\}}(\bT(\btheta))$. Note that the lower-dimensional problem based on $\bT_{\cA}(\btheta)$  can be easier compared to the original one. We also would like to emphasize that the dimensionality reduction to a smaller model $\cA$ is only for assisting the theoretical analysis and our Procedure  \ref{proc:1} does not need any knowledge of such set $\cA$. 

It is convenient to define $\bt_0 = \E\bT(\btheta^0) \in \mathbb{R}^P$.  Denote by 
\begin{equation}\label{def: setI}
\mathbb{I}:= \left\{\bt \in \mathbb{R}^{P}: \|\bt - \bt_0\|_{\max} \leq a_{np} := C_1 ( k^{1/2}+s^{3/2} )\tilde{c}_{np} \right\},
\end{equation}
where $C_1$ is some positive constant and $\tilde{c}_{np}=p^{-1/2}\log n + n^{-1/2}\log p$. 
For any subset $\cA \subset \{1,\cdots, p\}$, 
let $\mathbb{I}_{\cA}$ be the subspace of $\mathbb{I}$ when taking out the coordinates corresponding to $\E\bT_\cA(\btheta^0)$. Thus $\mathbb{I}_{\cA}\subset \mathbb{R}^{|\cA|(2|\cA|+3)}$. 
In addition to Conditions \ref{ass:regerr}--\ref{ass:eigensepa}, we need an assumption on the algorithmic stability of Procedure \ref{proc:1}. 

\begin{con}[Algorithmic stability]\normalfont \label{ass:stability}
	For any subset $\cA \subset \{1,\dots,p\}$ that satisfies $|\cA|\leq k\leq n\wedge p$, 
	there exists a positive sequence $\rho_{np}\to 0$ as $n\wedge p\to \infty$ such that
	\begin{align*}
	\sup_{|\cA|\leq k}\sup_{\bt_1,\bt_2\in \mathbb{I}_\cA} \frac{\big| S_\cA(\bt_2) \triangle S_\cA(\bt_1) \big|}{|S_\cA(\bt_1)| \wedge |S_\cA(\bt_2)|} =  O(\rho_{np}),
	\end{align*} 
	where $ \triangle$ stands for the symmetric difference between two sets. 
\end{con}

Intuitively the above condition assumes that  the knockoffs procedure is stable with respect to a small perturbation to the input $\bt$ in any lower-dimensional subspace $\mathbb{I}_{\cA}$. Under these regularity  conditions, the asymptotic FDR control of our IPAD procedure can be established. 

\begin{thm}[Robust FDR control]\label{thm:FDRcontrol}
	Assume that Conditions \ref{ass:regerr}--\ref{ass:stability} hold. Fix an arbitrary positive constant $\nu$. If $(s,k,n,p)$ satisfies  $s\vee k \leq n\wedge p$, $c_{np} \leq c/[r^2M^2C(\nu+2)]^{1/2}$, and $( k^{1/2}+s^{3/2})\tilde{c}_{np} \to 0$ 
	as $n\wedge p\rightarrow \infty$ with $c$ and $C$ some positive constants defined in Lemma \ref{lem:ineq} in Appendix, then the set of variables $\what{\cS}$ obtained by Procedure \ref{proc:1} (IPAD) with the LCD knockoff statistics controls the FDR \eqref{def:FDR} 
	to be no larger than $q + O\left( \rho_{np} + n^{-\nu} + p^{-\nu} \right)$. 
\end{thm}

Recall that by definition, the FDR is a function of $\bT(\hat{\btheta})$ and can be written as $\E\FDP(\bT(\hat{\btheta}))$ while the FDR computed with the oracle knockoffs, $\E\FDP(\bT(\btheta^0))$, is perfectly controlled to be no larger than $q$. This observation motivates us to first establish asymptotic equivalence of $\bT(\hat{\btheta})$ and  $\bT(\btheta^0)$ with large probability. Then a natural idea is to show that $\E\FDP(\bT(\hat{\btheta}))$ converges to $\E \FDP(\bT(\btheta^0))$ in probability, which turns out to be highly nontrivial because of  the discontinuity of $\FDP(\cdot)$ (the convergence would be straightforward via the Portmanteau lemma if $\FDP(\cdot)$ were continuous). 
Condition \ref{ass:stability} above provides a remedy to this issue by imposing the algorithmic stability assumption.

\subsection{Power analysis} \label{subsec:Power}
We have established the asymptotic FDR control for our IPAD procedure in Section \ref{subsec:FDR}. We now look at the other side of the coin -- the power \eqref{eq: power}.  Recall that in IPAD, we apply the knockoffs inference procedure to the knockoff statistics LCD, which are constructed using the augmented Lasso in \eqref{ex:lasso}. Therefore the final set of variables selected by IPAD is a subset of variables picked by the augmented Lasso. For this reason, the power of IPAD is always upper bounded by that of Lasso. We will show in this section that there is in fact no power loss relative to the augmented Lasso in the  asymptotic sense. 

\begin{con}[Signal strength I]\normalfont \label{ass:betamin2}
	For any subset $\cA\subset \cS^0$ that satisfies $|\cA|/s>1-\gamma$ for some $\gamma\in(0,1]$, it holds that 
	$\|\bbeta_{\cA}\|_1 > b_{np} sn^{-1/2}\log p$ for some positive sequence $b_{np}\to \infty$. 
\end{con}

\begin{con}[Signal  strength II]\normalfont \label{ass:active}
	There exists some constant $C_2\in (2(qs)^{-1},1)$ such that $|\cS_2|\geq C_2s$ with $\cS_2=\{j:|\beta_j|\gg (s/n)^{1/2}\log p\}$. 
\end{con}

Condition \ref{ass:betamin2} requires that the overall signal is not too weak, but is weaker than the conventional beta-min condition $\min_{j\in\cS^0} |\beta_j| \gg n^{-1/2}\log p$. 
Under Condition \ref{ass:active}, we can show that $|\what\cS|\geq C_2s$ with probability at least $1-O(p^{-\nu}+n^{-\nu})$ using similar techniques to those of Lemma 6 in \cite{RANK}. The intuition is that given $s\rightarrow \infty$, for a variable selection procedure to have high power it should select at least a reasonably large number of variables.  The result $|\what\cS|\geq C_2s$ will be used to derive the asymptotic order of threshold $T$, which is in turn crucial to establish the theorem below on power. 

\begin{thm}[Power guarantee]\label{thm:power}
	Assume that Conditions \ref{ass:regerr}--\ref{ass:eigensepa} and  \ref{ass:betamin2}--\ref{ass:active} hold. Fix an arbitrary positive constant $\nu$. If $(s,k,n,p)$ satisfies $2s\leq k \leq n\wedge p$, $c_{np} \leq c/(r^2M^2C(\nu+2))^{1/2}$, and $sk^{1/2}\tilde{c}_{np} \to 0$ as $n\wedge p\rightarrow \infty$ with $c$ and $C$ some positive constants defined in Lemma \ref{lem:ineq}, then both the Lasso procedure based on $(\bX, \by)$ and our IPAD procedure (Procedure \ref{proc:1}) have power bounded from below by $\gamma-o(1)$ as $n\wedge p\rightarrow \infty$. In particular, if $\gamma=1$ IPAD has no power loss compared to Lasso asymptotically.
\end{thm}

\section{Simulation studies} \label{sec:simu}
We have shown in Section \ref{sec:theory} that IPAD can asymptotically control the FDR in high-dimensional setting and there can be no power loss in applying the procedure. We next move on to numerically investigate the finite-sample performance of IPAD using synthetic data sets. We will compare IPAD with the knockoff filter in \cite{barber2015controlling} (BCKnockoff) and  the high-dimensional knockoff filter in \cite{barber2016knockoff} (HD-BCKnockoff). In what follows, we will first explain in detail the model setups and simulation settings,  then discuss the implementation of the aforementioned methods, and finally summarize the comparison results.    

\subsection{Simulation designs and settings} \label{Sim&set}
In all simulations, the  design matrix $ \bX \in \mathbb{R}^{ n \times p } $ is generated from the factor model
\begin{equation} \label{eq1}
	\bX = \bF^0 (\bLambda^{0})' + \sqrt{r\theta}\bE = \bC^0 + \sqrt{r\theta}\bE,
\end{equation}
where $ \bF^0 = (\bff_1^0,\cdots,\bff_n^0)' \in \mathbb{R}^{n\times r}$ is the matrix of latent factors, $\bLambda^0=(\blambda_1^0,\dots,\blambda_p^0)'\in\mathbb{R}^{p\times r}$ is the matrix of factor loadings, $\bE\in\mathbb{R}^{n\times p}$ is the matrix of model errors, and $ \theta $ is a constant controlling the signal-to-noise ratio. The term $\sqrt{r}$ is used to single out the effect of the number of factors in calculating the signal-to-noise ratio in factor model \eqref{eq1}. We then rescale each column of $ \bX $ to have $ \ell_2 $-norm one and simulate the response vector $ \by = (y_1,\cdots, y_n)' $ from the following model 
\begin{equation}
	y_i = f(\bx_i)+ \sqrt{c} \ep_i, \  i=1,\cdots, n, 
\end{equation}
where $f: \mathbb{R}^{p} \rightarrow \mathbb R$ is the link function which can be linear or nonlinear, $ c > 0 $ is a constant controlling the signal-to-noise ratio, and $\bep = (\ep_1,\cdots,\ep_n)'$ is the vector of model error. We next explain the four different designs of our simulation studies.

\subsubsection{Design 1: linear model with normal factor design matrix} \label{Sec4.1.1} 
The elements of $ \bF^0 $, $ \bLambda^0 $, $ \bE $, and $ \bep $ are drawn independently from $ \mathcal{N}(0,1) $. The link function takes a linear form, that is, 
\begin{equation*}
	\by = \bX \bbeta + \sqrt{c} \bep, 
\end{equation*}
where the coefficient vector $ \bbeta = (\beta_1,\cdots, \beta_p)'\in \mathbb{R}^p$ is generated by first choosing $ s $ random locations for the true signals and then setting $ \beta_j $ at each location to be either $A$ or $-A$ randomly with $A$ some positive value. The remaining $ p-s $ components of $\bbeta$ are set to zero. 

\subsubsection{Design 2: linear model with fat-tail factor matrix and serial dependence} \label{Sec4.1.2} 
The elements of $ \bE $ are generated as
\begin{equation}
	e_{ij} = \left( \frac{\nu-2}{\chi^2_{\nu,j}}\right) u_{ij},
\end{equation}
where $ u_{ij} \sim i.i.d.\ \mathcal{N}(0,1) $ for all $i = 1,\cdots, n$ and $j=1,\cdots, p$, and  $\chi^2_{\nu,j}, \  j=1,\cdots, p$ are i.i.d. random variables from chi-square distribution with $\nu=8$ degrees of freedom. The rest of the design is the same as in Design 1. It is worth mentioning that in this case, the entries of matrix $ \bE $ have fat-tail distribution with serial dependence in each column because of the common factor $\chi^2_{\nu,j}$. This design is used to check the robustness of IPAD method with respect to the serial dependence and the fat-tail distribution of $ \bE $. 

\subsubsection{Design 3: linear model with misspecified design matrix} \label{Sec4.1.3} 
To evaluate the robustness of IPAD procedure to the misspecification of the factor model structure \eqref{eq1},  we set $ \bLambda = \mathbf{0} $, $ r\theta = 1 $ and simulate the rows of matrix $ \bE $ independently from $ \mathcal{N}(\mathbf{0},\bSigma) $ with  $\bSigma = (\sigma_{ij})$, $ \sigma_{ij}= \rho^{|i-j|} $ for $ \leq i, j  \leq p$. The remaining design is the same as in Design 1. It is seen that our assumption on the independence of the entries of $\bE$ is violated.  This design is used to test the robustness of IPAD to misspecification of the factor model structure of $\bX$. 

\subsubsection{Design 4: nonlinear model with normal factor design matrix} \label{Sec4.1.4} 
Our last design is used to evaluate the performance of IPAD method when the link function $f$ is nonlinear. To be more specific, we assume  the following nonlinear model between the response and covariates
\begin{equation*}
	\by = \sin(\bX \bbeta)\exp(\bX \bbeta)+\sqrt{c}\bep,
\end{equation*}
where the coefficient vector $\bbeta$, design matrix $\bX$, and model error $\bep$ are generated similarly as in Design 1. 

\subsubsection{Simulation settings}

The target FDR level is set to be $ q = 0.2 $ in all simulations. 
For Design 1 and Design 2, we set $ n = 2000 $, $ p = 2000 $, $ A = 4 $, $ s = 50 $, $ c = 0.2 $, $ r = 3 $, and $ \theta = 1 $. In order to evaluate the sensitivity of our method to the dimensionality $ p $ and the model sparsity $ s $, we also explore the settings of $ p = 1000,3000$ and $ s = 100,150 $. In Design 3, we set $ r = 0 $ and $ \rho = 0,0.5$. In Design 4, since the model is nonlinear, we use nonparametric method to fit the model and consider lower-dimensional settings of $ p = 50, 250, 500 $. We also decrease the number of observations to $ n = 1000 $ and number of true variables to $ s = 10 $. Moreover, we set $ \theta =  1,2 $ and $ c = 0.1,0.2,0.3 $ to test the effects of signal-to-noise ratio on the performance of IPAD procedure in Design 4.      

\subsection{Estimation procedure}

In implementing the IPAD algorithm suggested in Section \ref{sec:metho}, we use the $ PC_{p1} $ criterion proposed in \cite{bai2002determining} to estimate the number of factors $r$. With an estimated number of factors $\hat r$, we use the principle component method discussed in Section \ref{subsec:FDR} to obtain an estimate $\what\bC$ of matrix $\bC^0$. 
Denote by $\what \bE = (\hat{ e}_{ij}) = \bX - \what\bC$. Recall that in the construction of knockoff variables, the distribution of $\bE$ needs to be estimated. Throughout our simulation studies, we misspecify the model and treat the entries of $\bE$ as i.i.d. Gaussian random variables. Under this working model assumption, the only unknown parameter is the variance which can be estimated by the following maximum likelihood estimator 
\begin{equation*}
	\hat{\sigma}^2 = (np)^{-1}\sum_{i=1}^{n} \sum_{j=1}^{p} \hat{e}_{ij}^2.
\end{equation*}
Then the knockoffs matrix $\what\bX$ is constructed using \eqref{Xk} with the entries of $\bE_{\hat{\bbbeta}}$ drawn independently from $\mathcal N(0, \hat{\sigma}^2)$. For the two comparison methods BCKnockoff and HD-BCKnockoff, we follow the implementation in \cite{barber2015controlling} and \cite{barber2016knockoff}, respectively. Thus it is seen that neither BCKnockoff nor HD-BCKnockoff uses the factor structure in $\bX$ when constructing the knockoff variables.  

In Designs 1--3, with the constructed empirical knockoffs matrix $\what\bX$ we apply the Lasso method to fit the model with $\by$ the response vector and $[\bX, \what{\bX}]$ the augmented design matrix. Then the LCD discussed in Section \ref{subsec:newmethod} is used in the construction of knockoff statistics. In Design 4, we assume the nonlinear relationship between the response and the covariates. In this case, random forest is used for estimation of the model. To construct the knockoff statistics, we use the variable importance measure of mean decrease accuracy (MDA) introduced in \cite{Breiman2001}. This measure is based on the idea that if a variable is unimportant, then rearranging its values should not degrade the prediction accuracy. The MDA for the $j$th variable, denoted as $\widehat{\text{MDA}}_j$, measures the amount of increase in prediction error when the values of the $j$th variable in the out-of-sample prediction are permuted randomly.  Then intuitively,  $\widehat{\text{MDA}}_j$ will be small and around zero if the $j$th variable is unimportant in predicting the response. For each original variable $\bx_j$, we compute $ W_j $ statistic as $ |\widehat{\text{MDA}}_j| - |\widehat{\text{MDA}}_{j+p}|$, $j=1,\cdots, p$.

\subsection{Simulation results}  \label{Sec4.3}  
For each method, we use 100 simulated data sets to calculate its empirical FDR and power, which are the average FDP and TDP (true discovery proportion as in (\ref{eq: power})) over 100 repetitions, respectively. Two different thresholds, knockoff and knockoff$+$  ($T_1$ and $T_2$ in Result \ref{thm:candes}, respectively), are used in the knockoffs inference implementation. It is worth mentioning that as shown in \cite{candes2016panning} and summarized in Result \ref{thm:candes}, knockoff$+$ controls FDR (\ref{def:FDR}) exactly while knockoff controls only the modified FDR (\ref{def:mFDR}).

Tables \ref{tab1} and \ref{tab2} summarize the results from Designs 1 and 2, respectively. As shown in Table \ref{tab1}, all approaches can control empirical FDR at the target level ($ q = 0.2 $) and knockoff+, which is more conservative, reduces power negligibly. It is worth mentioning that even for Design 2, in which the design matrix $ \bX $ is drawn from fat-tail distribution with serial dependence, we still have FDR under control with decent level of power. This suggests that the no serial correlation assumption in our theoretical analysis could just be technical. Compared to the results by BCKnockoff and HD-BCKnockoff, we see that using the extra information from the factor structure  in constructing knockoff variables can help with both FDR and power. Table \ref{tab2} also shows the effects of model sparsity on the  performance of various approaches. It can be seen that when the number of true signals is increased from 50 to 150, the FDR is still under control and the empirical power of IPAD remains steady. 

\begin{table}[htbp]	
	\centering
	{\scriptsize 
		\captionsetup{font=footnotesize}
		\caption{Simulation results for Designs 1 and 2 of Section \ref{Sim&set} with different values of dimensionality $ p $}\label{tab1}%
		\begin{threeparttable}
			\begin{tabular}{lcccccccccc}
				\toprule
				\toprule
				& \multicolumn{5}{c}{Design 1}          & \multicolumn{5}{c}{Design 2} \\
				\cmidrule(lr){2-6}\cmidrule(lr){7-11}
				& FDR   & Power & $ \text{FDR}_+ $  & $ \text{Power}_+ $ & $ R^2 $ & FDR   & Power & $ \text{FDR}_+ $  & $ \text{Power}_+ $ & $ R^2 $ \\
				\midrule
				& \multicolumn{10}{c}{$p = 1000$} \\
				\cmidrule{2-11}    
				IPAD  & 0.195 & 0.991 & 0.180 & 0.990 & 0.659 & 0.199 & 0.961 & 0.180 & 0.960 & 0.652 \\
				BCKnockoff & 0.207 & 0.942 & 0.192 & 0.938 & 0.659 & 0.172 & 0.887 & 0.152 & 0.885 & 0.653 \\
				\cmidrule{2-11}
				& \multicolumn{10}{c}{$p = 2000$} \\
				\cmidrule{2-11}    
				IPAD  & 0.194 & 0.979 & 0.179 & 0.979 & 0.649 & 0.199 & 0.935 & 0.183 & 0.933 & 0.656 \\
				HD-BCKnockoff & 0.142 & 0.706 & 0.127 & 0.691 & 0.649 & 0.136 & 0.607 & 0.113 & 0.581 & 0.644 \\
				\cmidrule{2-11}          
				& \multicolumn{10}{c}{$p = 3000$} \\
				\cmidrule{2-11}    
				IPAD  & 0.191 & 0.964 & 0.176 & 0.963 & 0.652 & 0.188 & 0.913 & 0.171 & 0.911 & 0.658 \\
				HD-BCKnockoff & 0.172 & 0.668 & 0.149 & 0.658 & 0.652 & 0.125 & 0.559 & 0.099 & 0.524 & 0.651 \\
				\bottomrule
			\end{tabular}%
			\begin{tablenotes}
				\item {\footnotesize Note that $\text{FDR}_+ $ and $ \text{Power}_+ $ are the values of  FDR  and Power corresponding to the knockoff+ threshold $ T_2 $.}
			\end{tablenotes}
		\end{threeparttable}
	}
\end{table}%

\begin{table}[htbp]	
	\centering
	{\scriptsize
		\captionsetup{font=footnotesize}
		\caption{Simulation results for Designs 1 and 2 of Section \ref{Sim&set} with different sparsity level $ s $}\label{tab2}%
		\begin{tabular}{lcccccccccc}
			\toprule
			\toprule
			
			& \multicolumn{5}{c}{Design 1}          & \multicolumn{5}{c}{Design 2} \\
			
			\cmidrule(lr){2-6}\cmidrule(lr){7-11} 
			& FDR   & Power & $ \text{FDR}_+ $  & $ \text{Power}_+ $ & $ R^2 $ & FDR   & Power & $ \text{FDR}_+ $  & $ \text{Power}_+ $ & $ R^2 $ \\
			\midrule
			
			& \multicolumn{10}{c}{$s = 50$} \\
			\cmidrule{2-11}    
			IPAD  & 0.194 & 0.979 & 0.179 & 0.979 & 0.649 & 0.199 & 0.935 & 0.183 & 0.933 & 0.656 \\
			
			HD-BCKnockoff & 0.142 & 0.706 & 0.127 & 0.691 & 0.649 & 0.136 & 0.607 & 0.113 & 0.581 & 0.644 \\
			\cmidrule{2-11}          
			
			& \multicolumn{10}{c}{$s = 100$} \\
			
			\cmidrule{2-11}
			IPAD  & 0.191 & 0.978 & 0.183 & 0.977 & 0.783 & 0.181 & 0.937 & 0.174 & 0.936 & 0.789 \\
			
			HD-BCKnockoff & 0.152 & 0.703 & 0.140 & 0.698 & 0.787 & 0.106 & 0.583 & 0.097 & 0.573 & 0.778 \\
			
			\cmidrule{2-11}          
			& \multicolumn{10}{c}{$s = 150$} \\
			
			\cmidrule{2-11}    
			IPAD  & 0.183 & 0.973 & 0.178 & 0.972 & 0.842 & 0.188 & 0.935 & 0.182 & 0.935 & 0.848 \\
			
			HD-BCKnockoff & 0.139 & 0.660 & 0.130 & 0.654 & 0.858 & 0.115 & 0.578 & 0.106 & 0.570 & 0.843 \\
			\bottomrule
	\end{tabular}}%
\end{table}%

Table \ref{tab3} is devoted to the case of Design 3, where the rows of matrix $\bX$ are generated independently from multivariate normal distribution with AR(1) correlation structure. This is a setting where the factor model structure in $\bX$ is misspecified. Since  BCknockoff and HD-BCknockoff make no use of the factor structure in generating knockoff variables, in both low- and high-dimensional examples both methods control FDR exactly at the target level. IPAD based methods have empirical FDR slightly over the target level, which may be caused by the misspecification of the factor structure. On the other hand, IPAD based approaches have much higher empirical power than comparison methods.    

\begin{table}[htbp]
	\centering
	{\scriptsize
		\captionsetup{font=footnotesize}
		\caption{Simulation results for Design 3 of Section \ref{Sim&set}}	\label{tab3}%
		
		\begin{tabular}{lcccccccccc}
			
			\toprule
			
			\toprule
			
			& \multicolumn{5}{c}{$ \rho = 0 $}            & \multicolumn{5}{c}{$ \rho = 0.5 $} \\
			
			\cmidrule(lr){2-6}\cmidrule(lr){7-11} 
			& FDR   & Power & $ \text{FDR}_+ $  & $ \text{Power}_+ $ & $ R^2 $ & FDR   & Power & $ \text{FDR}_+ $  & $ \text{Power}_+ $ & $ R^2 $ \\
			\midrule
			
			& \multicolumn{10}{c}{$p = 1000$} \\
			
			\cmidrule{2-11}    
			IPAD  & 0.204 & 0.995 & 0.189 & 0.995 & 0.444 & 0.226 & 0.984 & 0.216 & 0.984 & 0.446 \\
			
			BCKnockoff & 0.188 & 0.919 & 0.172 & 0.917 & 0.444 & 0.137 & 0.827 & 0.117 & 0.821 & 0.445 \\
			
			\cmidrule{2-11}          
			& \multicolumn{10}{c}{$p = 2000$} \\
			
			\cmidrule{2-11}    
			IPAD  & 0.203 & 0.993 & 0.189 & 0.993 & 0.447 & 0.220 & 0.982 & 0.202 & 0.980 & 0.445 \\
			
			HD-BCKnockoff & 0.151 & 0.630 & 0.126 & 0.603 & 0.449 & 0.115 & 0.522 & 0.090 & 0.467 & 0.442 \\
			
			\cmidrule{2-11}          
			& \multicolumn{10}{c}{$p = 3000$} \\
			
			\cmidrule{2-11}    
			IPAD  & 0.225 & 0.988 & 0.205 & 0.987 & 0.445 & 0.219 & 0.979 & 0.206 & 0.978 & 0.443 \\
			
			HD-BCKnockoff & 0.150 & 0.589 & 0.126 & 0.560 & 0.446 & 0.092 & 0.439 & 0.064 & 0.381 & 0.447 \\
			
			\bottomrule
			
	\end{tabular}}%
\end{table}

Table \ref{tab4} corresponds to Design 4 in which response $ \by $ is related to $ \bX $ nonlinearly. Since BCKnockoff and HD-BCKnockoff are designed for linear models, only the results from IPAD method are reported. It can be seen form Table \ref{tab4} that IPAD approach can control FDR with reasonably high power even in the nonlinear setting. We also observe that in nonlinear setting, the power of IPAD deteriorates faster as dimensionality $ p $ increases compared to the linear setting.

\begin{table}[!htbp] \centering
	\centering
	{\scriptsize
		\captionsetup{font=footnotesize} 
		\caption{Simulation results for Design 4 of Section \ref{Sim&set}}	\label{tab4} 
		\begin{tabular}{@{\extracolsep{5pt}} lcccccccccc} 
			\hline \hline
			& \multicolumn{5}{c}{$\theta  = 1$} & \multicolumn{5}{c}{$\theta  = 2$}\\
			\cmidrule(lr){2-6} \cmidrule(lr){7-11}  
			& FDR   & Power & $ \text{FDR}_+ $  & $ \text{Power}_+ $ & $ R^2 $ & FDR   & Power & $ \text{FDR}_+ $  & $ \text{Power}_+ $ & $ R^2 $ \\
			\hline
			&\multicolumn{10}{c}{$ p = 50$}\\
			\cmidrule(lr){2-11} 
			$ c = 0.1 $ & $0.109$ & $0.839$ & $0.081$ & $0.720$ & $0.707$ & $0.110$ & $0.943$ & $0.061$ & $0.858$ & $0.707$ \\   
			$ c = 0.2 $ & $0.137$ & $0.847$ & $0.068$ & $0.726$ & $0.547$ & $0.097$ & $0.920$ & $0.061$ & $0.837$ & $0.547$ \\  
			$ c = 0.3 $ & $0.137$ & $0.765$ & $0.091$ & $0.582$ & $0.451$ & $0.123$ & $0.907$ & $0.076$ & $0.774$ & $0.451$ \\ 
			\cmidrule(lr){2-11}
			&\multicolumn{10}{c}{$ p = 250$}\\
			\cmidrule(lr){2-11}
			$ c = 0.1 $ & $0.189$ & $0.740$ & $0.104$ & $0.504$ & $0.702$ & $0.174$ & $0.876$ & $0.139$ & $0.788$ & $0.702$ \\
			$ c = 0.2 $ & $0.218$ & $0.666$ & $0.131$ & $0.522$ & $0.552$ & $0.209$ & $0.831$ & $0.118$ & $0.660$ & $0.552$ \\
			$ c = 0.3 $ & $0.200$ & $0.569$ & $0.101$ & $0.361$ & $0.451$ & $0.224$ & $0.766$ & $0.141$ & $0.599$ & $0.451$ \\
			\cmidrule(lr){2-11}
			&\multicolumn{10}{c}{$ p = 500$}\\
			\cmidrule(lr){2-11}
			$ c = 0.1 $ & $0.243$ & $0.661$ & $0.169$ & $0.497$ & $0.702$ & $0.223$ & $0.831$ & $0.173$ & $0.740$ & $0.702$ \\
			$ c = 0.2 $ & $0.204$ & $0.507$ & $0.111$ & $0.266$ & $0.543$ & $0.216$ & $0.749$ & $0.126$ & $0.594$ & $0.543$ \\ 
			$ c = 0.3 $ & $0.247$ & $0.478$ & $0.128$ & $0.299$ & $0.451$ & $0.241$ & $0.691$ & $0.156$ & $0.550$ & $0.451$ \\
			\hline 
	\end{tabular}}
\end{table}

\section{Empirical analysis} \label{sec:rdata} 
Our simulation results in Section \ref{sec:simu} suggest that IPAD is a powerful approach with asymptotic FDR control. We further examine the application of IPAD to the quarterly data on 109 macroeconomic variables from the third quarter of year 1960 (1960Q3) to the fourth quarter of year 2008 (2008Q4) in the United States discussed in \cite{stock2012generalized}. These variables are transformed by taking logarithms and/or differencing following \cite{stock2012generalized}. Our real data analysis consists of two parts. In the first part, we focus on the performance of IPAD method in terms of empirical FDR and power. In the second part, the forecasting performance of IPAD method will be evaluated.

\subsection{Simulation study} \label{Sec5.1}
To evaluate the performance of IPAD approach in terms of empirical FDR and power with real economic data, we set up one additional Monte Carlo simulation study. In this design, we use the transformed macroeconomic variables described above as the design matrix $ \bX $, but simulate response $ \by $ from the model in Design 1 in Section \ref{Sim&set}. We set the number of true signals, the amplitude of signals,  and the target FDR level to $ s = 10 $, $ A = 4 $,  and $ q = 0.2 $, respectively. 

Table \ref{tab5} shows the results for IPAD and HD-BCKnockoff approaches. As expected, HD-BCKnockoff can control FDR but suffers from lack of power. On the other hand, IPAD has empirical FDR slightly higher than the target level ($ q = 0.2 $) while its power is reasonably high. These results are consistent with our theory in Section \ref{sec:theory} because IPAD only controls FDR asymptotically. Additional reason for having slightly higher FDR than the target level can be deviation of the design matrix from our factor model assumption. Overall this simulation study indicates that IPAD can control FDR at around the target level with reasonably high power when we use the macroeconomic data set. In the next section, using the same data set we will compare the forecasting performance of IPAD with that of some commonly used forecasting methods in the literature.

\begin{table}[htbp]
	\centering
	{\scriptsize
		\captionsetup{font=footnotesize} 
		\caption{Real data simulation results with $(n, p ) = (195, 109)$}	\label{tab5}%
		\begin{tabular}{lccccc}
			\toprule
			\toprule
			& FDR   & Power & $ \text{FDR}_+ $  & $ \text{Power}_+ $ & $ R^2 $  \\
			\midrule
			& \multicolumn{5}{c}{$c = 0.2$} \\
			\cmidrule{2-6}    
			IPAD  & 0.278 & 0.812 & 0.223 & 0.796 & 0.747 \\
			HD-BCKnockoff & 0.096 & 0.009 & 0.010 & 0.002 & 0.758 \\
			\cmidrule{2-6}         
			& \multicolumn{5}{c}{$c = 0.3$} \\
			\cmidrule{2-6}    
			IPAD  & 0.280 & 0.757 & 0.221 & 0.723 & 0.665 \\
			HD-BCKnockoff & 0.149 & 0.121 & 0.027 & 0.036 & 0.678 \\
			\cmidrule{2-6}
			& \multicolumn{5}{c}{$c = 0.5$} \\
			\cmidrule{2-6}    
			IPAD  & 0.286 & 0.661 & 0.215 & 0.571 & 0.560 \\
			HD-BCKnockoff & 0.119 & 0.009 & 0.008 & 0.001 & 0.554 \\          
			\bottomrule
	\end{tabular}}%
\end{table}%

\subsection{Forecasting results} \label{Sec5.2}
In this section, we apply the IPAD approach to the real economic data set for  forecasting. One-step ahead prediction is conducted using rolling window of size 120. More specifically, one of the 109 variables is chosen as the response and the remaining 108 variables are treated as predictors. For each quarter between 1990Q3 and 2008Q4, we use the previous 120 periods for model fitting and then one-step ahead prediction is conducted based on the fitted model. We compare the following different methods, where each method is implemented in a same way as IPAD for one-step ahead prediction.

\begin{enumerate}
	\item[1)] Autoregression of order one (AR(1)). Assume that 
	\begin{equation*}
		y_t = \alpha_0 + \rho y_{t-1} + \varepsilon_t,
	\end{equation*}
	where $ y_t $ is regressed on $ y_{t-1} $, and $\alpha_0$ and $\rho$ are the AR(1) coefficients that need to be estimated. With the ordinary least squares estimates $ \hat{\alpha}_0 $ and $ \hat{\rho} $,  the one-step ahead prediction based on this model is $ \hat{y}_{T+1} = \hat{\alpha}_0 + \hat{\rho} y_{T} $.
	
	\item[2)] Factor augmented AR(1) (FAR). We first extract $ m$ factors $\bff_1,\cdots, \bff_m$ form the 109 transformed macroeconomic variables by principal component analysis (PCA). Denote by $\tilde\bff_t\in \mathbb R^{m}$ the factor vector at time $t$ extracted from the rows of matrix $[\bff_1,\cdots, \bff_m] \in \mathbb{R}^{n\times m}$. Then we regress $ y_t $ on $ y_{t-1} $ and $\tilde \bff_{t-1}$ and fit the following model
	\begin{equation*}
		y_t = \alpha_0 + \rho y_{t-1} + \bgamma' \tilde\bff_{t-1} +\varepsilon_t
	\end{equation*}
	with $\bgamma \in \mathbb{R}^m$. The number of factors $m$ is determined using the $ PC_{p1} $ criterion in \cite{bai2002determining}. Similar to AR(1) model, one-step ahead forecast of $ y_t $ at time $ T $ is 
	\begin{equation*}
		\hat{y}_{T+1} = \hat{\alpha}_0 + \hat{\rho} y_{T} + \hat{\bgamma}' \tilde\bff_{T}. 
	\end{equation*}
	
	\item[3)] Lasso method. The $ y_t $ is regressed on $ y_{t-1} $, $ \tilde\bff_{t-1} $, and the 108 transformed macroeconomic variables $ \mathbf{z}_{t-1}\in \mathbb R^{108} $ at time $t-1$
	\begin{equation*}
		y_t = \alpha_0 + \rho y_{t-1} + \bgamma' \tilde\bff_{t-1} + \bdelta' \mathbf{z}_{t-1} +\varepsilon_t,
	\end{equation*}
	where $\tilde{\bff}_t$ is the same as in the FAR(1) model, and $\alpha_0, \rho$, and $\bdelta \in \mathbb R^{108}$ are regression coefficients that need to be estimated. The coefficients are estimated by Lasso method with regularization parameter chosen by the cross-validation. With the estimated Lasso coefficient vector $\hat\bbeta_{\text{Lasso}}$, one-step ahead forecast of $ y_t $ at time $ T $ is 
	\begin{equation*}
		\hat{y}_{T+1} = \hat{\bbeta}_{\text{Lasso}}'\bx_{T},
	\end{equation*}
	where $\bx_T$ is the augmented predictor vector at time $T$.
	
	\item[4)] IPAD method. We  regress $ y_t $  on the augmented vector $(y_{t-1}, \mathbf{z}_{t-1}')' $. The lagged variable $ y_{t-1} $ is assumed to be always in the model. To account for this, we implement IPAD in three steps. First, we regress $y_t$ on $y_{t-1}$ and obtain the residuals $e_{y,t}$. Second, we regress each of the 108 variables in $ \mathbf{z}_{t-1}$ on $y_{t-1}$ and obtain the residual vector $\be_{z,t-1}$. At last, we fit model \eqref{model1}--\eqref{factmodel} using the IPAD approach by treating $e_{y,t}$ as the response and $\be_{z,t-1}$ as predictors, which returns us a set of selected variables (a subset of the 108 macroeconomic variables). With the set of variables $\widehat{\cS}$ selected by IPAD, we fit the following model by the least-squares regression
	\begin{equation}\label{fitted}
		y_{t} = \alpha_0 + \rho y_{t-1} +\bdelta' \mathbf{z}_{t-1, \widehat{\cS}} + \varepsilon_t,
	\end{equation}
	where $\mathbf{z}_{t, \widehat{\cS}}$ stands for the subvector of $\mathbf z_t$ corresponding to the set of variables $\widehat{\cS}$ selected by IPAD at time $t$.  Since  $\widehat{\cS}$ from IPAD is random due to the randomness in generating knockoff variables, we apply the IPAD procedure 100 times and compute  the average of these 100 one-step ahead predictions based on \eqref{fitted} and use the mean value as the final predicted value of $y_{T+1}$.
\end{enumerate}       

Table \ref{tab6} shows the root mean-squared prediction error (RMSE) of these methods. As can be seen, the RMSE of IPAD is very close to those of comparison methods. To statistically compare the relative prediction accuracy of IPAD versus other approaches, we have used the  Diebold--Mariano test \cite{diebold1995comparing}, where the square of one-step ahead prediction error is used as the loss function. Table \ref{tab7} reports the test results. The results indicate that one-step ahead prediction accuracy of IPAD is comparable to other approaches.

\begin{table}[htbp]
	\centering
	{\scriptsize
		\captionsetup{font=footnotesize} 
		\caption{Root mean-squared error of one-period ahead forecast of various macroeconomic variables}	\label{tab6}%
		\begin{tabular}{lcccc}
			\toprule
			\toprule
			& AR    & FAR   & Lasso & IPAD \\
			\midrule
			RGDP  & 2.245 & 1.929 & 2.070 & 2.106 \\
			CPI-ALL & 1.526 & 1.552 & 1.579 & 1.571 \\
			Imports & 7.549 & 5.871 & 6.595 & 6.993 \\
			IP: cons dble & 9.683 & 8.353 & 8.424 & 9.175 \\
			Emp: TTU & 1.112 & 0.989 & 1.167 & 1.100 \\
			U: mean duration & 0.573 & 0.487 & 0.502 & 0.494 \\
			HStarts: South & 0.074 & 0.071 & 0.076 & 0.074 \\
			NAPM new ordrs & 4.800 & 4.378 & 4.659 & 4.673 \\
			PCED-NDUR-ENERGY & 31.927 & 32.121 & 33.546 & 32.164 \\
			Emp. Hours & 2.102 & 1.899 & 2.080 & 1.944 \\
			FedFunds & 0.421 & 0.396 & 0.406 & 0.392 \\
			Cons credit & 2.573 & 2.537 & 2.648 & 2.580 \\
			EX rate: Canada & 10.132 & 10.139 & 10.122 & 10.113 \\
			DJIA  & 23.117 & 23.997 & 24.585 & 23.398 \\
			Consumer expect & 6.496 & 6.888 & 6.681 & 6.661 \\
			\bottomrule
	\end{tabular}}%
\end{table}%

\begin{table}[htbp]
	\centering
	{\scriptsize
		\captionsetup{font=footnotesize} 
		\caption{Diebold--Mariano test for comparing prediction accuracy of IPAD against other procedures}	\label{tab7}%
		\begin{tabular}{lccc}
			\toprule
			\toprule
			& IPAD vs. AR & IPAD vs. FAR & IPAD vs. Lasso \\
			\midrule
			RGDP  & -0.780 & 1.160 & 0.462 \\
			CPI-ALL & 0.521 & 0.394 & -0.218 \\
			Imports & -0.976 & $2.631^{**}$ & 1.464 \\
			IP: cons dble & -1.026 & 1.567 & $2.487^{*}$ \\
			Emp: TTU & -0.140 & 1.692 & -1.845 \\
			U: mean duration & -$3.383^{***}$ & 0.672 & -0.505 \\
			HStarts: South & 0.096 & 0.821 & -0.766 \\
			NAPM new ordrs & -0.517 & 1.814 & 0.076 \\
			PCED-NDUR-ENERGY & 0.753 & 0.049 & -1.759 \\
			Emp. Hours & -1.200 & 0.297 & -$2.063^{*}$ \\
			FedFunds & -0.971 & -0.134 & -0.625 \\
			Cons credit & 0.207 & 0.359 & -0.661 \\
			EX rate: Canada & -0.466 & -0.138 & -0.037 \\
			DJIA  & 0.585 & -0.959 & -1.428 \\
			Consumer expect & 1.212 & -1.038 & -0.277 \\
			\bottomrule
	\end{tabular}}%
\end{table}%

\begin{figure}[!htb]
	\centering
	\subfloat[1990-1999]{\includegraphics[width=0.7\textwidth]{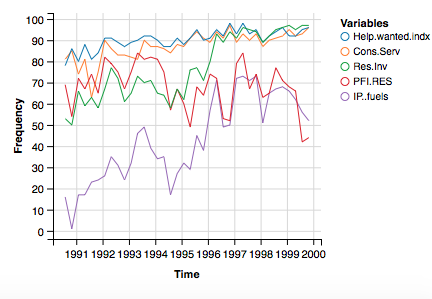}\label{Sub_fig1}}
	\hfill
	\subfloat[2000-2008]{\includegraphics[width=0.7\textwidth]{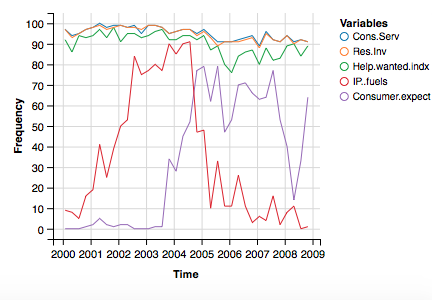}\label{Sub_fig2}}
	\caption{Frequencies of top selected variables in predicting real GDP growth. The set of selected variables are index of help-wanted advertising in newspapers (Help wanted indx), real personal consumption expenditures - services (Cons-Serv), real gross private domestic investment - residential (Res.Inv), residential price index (PFI-RES), industrial production index - fuels (IP:fuels), and  University of Michigan index of consumer expectations (Consumer expect).}
	\label{fig_1} 
\end{figure}

It is worth mentioning that one main advantage of IPAD is its interpretability and stability. Using IPAD for forecasting, we not only enjoy the same level of accuracy as other methods but also obtain the information on variable importance with stability. Recall that for each one-step ahead prediction, we apply IPAD 100 times and obtain 100 sets of selected variables. Thus we can calculate the selection frequency of each variable in each one-step ahead prediction.  Figure \ref{fig_1} depicts the frequencies of top five selected variables in predicting real GDP growth before and after year 2000, where the variable importance is ranked according to the aggregated frequencies over the entire time period before or after 2000. We have experimented with different cutoff years around year 2000, and the top five ranked variables stay the same so only the results corresponding to cutoff year 2000 are reported.  Changes in index of help wanted advertising in newspapers, percentages of changes in real personal consumption of services, and percentage of changes in real gross private domestic investment in residential sector were the top three important variables in predicting real GDP growth during the whole period. It is interesting to see that percentage of changes in residential price index was among top five important variables in predicting GDP growth during the 90s, and then starting from year 2000 it was replaced by changes in index of consumer expectations about stability of economy. Moreover,  it is also seen that the percentage of changes in industrial production of fuels was of great importance for predicting real GDP growth during some periods but not the others.

As a comparison, it is very difficult to interpret the results of FAR. As for Lasso based method, there is no theoretical guarantee on FDR control and in addition, Lasso usually gives us models with much larger size. For instance, in predicting real GDP growth,  IPAD on average selects 5.42 macroeconomic variables while Lasso on average selects 13.32 variables. To summarize, our real data analysis indicates that IPAD is an applicable approach for controlling FDR with competitive prediction power and high interpretability and stability.                 

\section{Discussions} \label{Sec6}
We have suggested in this paper a new procedure IPAD for feature selection in high-dimensional linear models that achieves asymptotic FDR control while retaining high power. Our model setting involves a latent factor model that is motivated by applications in economics and finance. 
Our method falls in the general model-X knockoffs framework in \cite{candes2016panning}, but allows the unknown covariate distribution for the knockoff variable construction. With the LCD knockoff statistics, we have shown that the FDR of IPAD can be asymptotically under control while the power can be asymptotically the same as that of Lasso. Our simulation study and empirical analysis also suggest that IPAD has highly competitively performance compared to many widely used forecasting methods such as Lasso and FAR, but with much higher interpretability and stability. 

Our work has focused on the scenario of static models. It would be interesting to extend the IPAD procedure to high-dimensional dynamic models with time series data. It is also interesting to consider nonlinear models and more flexible machine learning methods for forecasting as well as more refined factor model structures on the covariates for the knockoffs inference with IPAD, and develop theoretical guarantees for
the IPAD framework in these more general model settings. These extensions are beyond the scope of the current paper and are interesting topics for future research.

\bibliographystyle{chicago}
\bibliography{references}


\newpage
\appendix
\setcounter{section}{0}
\renewcommand{\theequation}{A.\arabic{equation}}
\setcounter{equation}{0}

\quad \vspace{0.05in}

\begin{center}
\textbf{\Large Appendix}
\end{center}

\bigskip

\noindent This appendix contains all the proofs and technical details for the theoretical results of the paper. In particular, Section \ref{secA} details the proofs of Lemmas \ref{lem:rep}--\ref{lem:lowdimrep} and Theorems \ref{thm:FDRcontrol}--\ref{thm:power}, Section \ref{secB} presents some key lemmas and their proofs, and Section \ref{sec:additionallem} provides some additional technical lemmas and their proofs. 

To ease the technical presentation, let us introduce some notation. We denote by $\lesssim$ the inequality up to some positive constant factor. Restricting the columns of $\bX$ and $\wtilde{\bX}(\hat{\btheta})$ to the variables in index set $\cA$ such that $|\cA| \leq k$, we obtain the $n\times k$ submatrices $\bX_{\cA}$ and $\wtilde{\bX}_{\cA}(\hat{\btheta})$, respectively.  
Moreover, we define $\bT_{\cA}(\hat{\btheta}) := \vect(\vech\bU_{\cA}(\hat{\btheta}),\bv_{\cA}(\hat{\btheta})) \in \mathbb{R}^{k(2k+3)}$ with $\bU_{\cA}(\hat{\btheta})$ the principle submatrix of $\bU(\hat\btheta)$ formed by columns and rows in set $\cA$, and $\bv_{\cA}(\hat{\btheta})$ the subvector of $\bv(\hat\btheta)$ formed by components in set $\cA$. Then it is easy to see that $\bU_{\cA}(\hat\btheta) = n^{-1}[\bX_{\cA},\wtilde{\bX}_{\cA}(\hat{\btheta})]'[\bX_{\cA},\wtilde{\bX}_{\cA}(\hat{\btheta})]$ and $\bv_{\cA}(\hat\btheta)  = n^{-1}[\bX_{\cA},\wtilde{\bX}_{\cA}(\hat{\btheta})]'\by $.
For the oracle factor loading matrix $\bLambda^0$, with a slight abuse of notation we use $\bLambda_{\cA}^0$ to denote the \textit{row} restricted to the variables in $\cA$ for notational convenience. Recall that $\nu>0$ is a fixed positive number, $c_{np}= (p^{-1}\log n)^{1/2} + (n^{-1}\log p)^{1/2}$, and $\tilde{c}_{np}=p^{-1/2}\log n + n^{-1/2}\log p$. We define $\pi_{np}=n^{-\nu}+p^{-\nu}$. Since $\lambda$ is fixed at $ C_0n^{-1/2}\log p$,  in all the proofs we will drop the dependence of various quantities on $\lambda$ whenever there is no confusion.

\section{Proofs of main results} \label{secA}

\subsection{Proof of Lemma \ref{lem:rep}} \label{secA.1}
For $\lambda$ fixed at $ C_0n^{-1/2}\log p$ and each given $\btheta$, $W_j(\btheta)=w_j([\bX,\wtilde{\bX}(\btheta)],\by)$ depends only on $\hat\bbeta^{\textsf{aug}}(\btheta)$ by the LCD construction. 
Moreover, the Lasso solution $\hat\bbeta^{\textsf{aug}}(\btheta)$ satisfies the Karush--Kuhn--Tucker (KKT) conditions:  
\begin{align}
&\bv(\btheta) - \bU(\btheta)\hat\bbeta^{\textsf{aug}}(\btheta) = n^{-1}\lambda \bz, \label{lem:kkt1}\\
& \text{where } \bz = (z_1,\cdots, z_{2p})^T \text{ with } z_j \in 
\begin{cases}
\{\sgn(\hat{\beta}_j)\} &\mbox{if}~~\hat{\beta}_j \not= 0,\\
[-1,1] &\mbox{if}~~\hat{\beta}_j = 0,
\end{cases}
~~~\mbox{for}~~~j=1,\dots,2p. \label{lem:kkt2}
\end{align}
This means that $\hat\bbeta^{\textsf{aug}}(\btheta)$ depends on the data $([\bX,\wtilde{\bX}(\btheta)],\by))$ only through $\bU(\btheta)$ and $\bv(\btheta)$. 
Thus using notation $\bT(\btheta)=\vect(\vech\bU(\btheta), \bv(\btheta))$ with the fact that $\bU(\btheta)$ is symmetric,  we can reparametrize $w_j([\bX,\wtilde{\bX}(\btheta)],\by)$ as $w_j(\bT(\btheta))$ with a slight abuse of notation. 
Furthermore, note that the thresholds $T_1$ and $T_2$ are both completely determined by $w_j(\bT(\btheta))$. Consequently, by the construction of $\what\cS$ we can see that $\what\cS$ depends only on $\bT(\btheta)$, which completes the proof of Lemma \ref{lem:rep}.

\subsection{Proof of Lemma \ref{lem:lowdimrep}} \label{secA.2} 
We continue to use the same $\lambda$ and $\btheta$ as in Lemma \ref{lem:rep} and its proof. Recall that $S_\cA(\bt_\cA)$ represents the outcome of first restricting ourselves to the smaller set of variables $\cA$ and then applying IPAD to $\bT_{\cA}(\btheta)=\bt_{\cA}$ to further select variables from $\cA$. Also recall that $\cA^*(\btheta)$ is the support of knockoff statistics $W_j(\btheta)$. Thus the knockoff threshold $T_1$ or $T_2$ depends only on $W_j(\btheta)$ with $j\in \cA^*(\btheta)$. 

On the other hand, when we restrict ourselves to $\cA \supset \cA^*(\btheta)$ we solve the following KKT conditions with respect to $\tilde\bbeta := (\tilde\beta_1,\cdots, \tilde\beta_{2|\cA|})^T \in \mathbb{R}^{2|\cA|}$ to get the Lasso solution:
\begin{align}
& \tilde\bbeta 
= (\bU_{\cA}(\btheta)'\bU_{\cA}(\btheta))^{-1}(\bv_{\cA}(\btheta) - n^{-1}\lambda \tilde{\bz}), \label{KKT2a}\\
&\text{where } \tilde{\bz} = (\tilde z_1,\cdots, \tilde z_{2|\cA|})^T \text{ with } \tilde z_j \in 
\begin{cases}
\{\sgn(\tilde{\beta}_j)\} &\mbox{if}~~\tilde{\beta}_j \not= 0,\\
[-1,1] &\mbox{if}~~\tilde{\beta}_j = 0,
\end{cases}
~~~\mbox{for}~~~j=1,\dots,2|\cA|. \label{KKT2b}
\end{align}
Since $\lambda$ is always fixed at the same value $C_0n^{-1/2}\log p$, it is seen that the solution to the above KKT conditions is identical to $\hat\bbeta_{\cA\cA}^{\textsf{aug}}(\btheta)$, where the latter denotes the subvector of $\hat\bbeta^{\textsf{aug}}(\btheta)$ formed by stacking $\hat\beta_{j_1}^{\textsf{aug}}(\btheta)$, $j_1\in \cA$ and 
$\hat\beta_{p+j_2}^{\textsf{aug}}(\btheta)$, $j_2\in \cA$ all together. Therefore, the Lasso solution to \eqref{KKT2a}--\eqref{KKT2b} and the Lasso solution to \eqref{lem:kkt1}--\eqref{lem:kkt2} have the identical support (when viewed in the original $2p$-dimensional space) and in addition, identical values on the support. This guarantees that  $S_{\{1,\dots,p\}}(\bT(\btheta))$ and $S_{\cA}(\bT_\cA(\btheta))$ are identical and thus concludes the proof of Lemma \ref{lem:lowdimrep}.

\subsection{Proof of Theorem \ref{thm:FDRcontrol}} \label{secA.3} 
Recall that for a given $\btheta$, $\cA^*(\btheta)$ is the support of knockoff statistics $(W_1(\btheta),\cdots, W_p(\btheta))'$.   Define set
\begin{align*}
\what{\cA}(\hat\btheta):=\cA^*(\hat{\btheta})\cup \cA^*(\btheta^0).
\end{align*}
It follows from  \eqref{def: k} that the cardinality of $\what\cA(\hat\btheta)$ is bounded by $k$. Hereafter we write $\what{\cA}(\hat\btheta)$ as $\what\cA$ for notational simplicity. 

By Lemmas \ref{lem:rep}--\ref{lem:lowdimrep} and the definition of the FDP, we know that $S_{\{1,\dots,p\}}(\bT(\hat\btheta))= S_{\widehat\cA}(\bT_{\widehat\cA}(\hat\btheta))$ and thus the resulting FDR's are the same. Therefore, we can restrict ourselves to the smaller model $\widehat\cA$ when studying the FDR of IPAD.  The same arguments as above also hold for the oracle knockoffs; that is, the FDR of IPAD applied to $\bT(\btheta^0)$ is the same as that applied to $\bT_{\what{\cA}}(\btheta^0)$. Note that all the FDR's we discuss here are with respect to the full model $\{1,\cdots, p\}$.   
For this reason, in what follows we will abuse the notation and use $\FDR_{\what\cA}(\bT_{\what{\cA}}(\hat{\btheta}))$ and $\FDR_{\what{\cA}}(\bT_{\what{\cA}}(\btheta^0))$   to denote the  FDR of IPAD based on $\bT_{\what{\cA}}(\btheta)$ and $\bT_{\what{\cA}}(\btheta^0)$, respectively. We want to emphasize that although we put a subscript $\what\cA$ in FDR's, their values are still deterministic as argued above. Summarizing the facts, we obtain
\begin{align*}
\FDR_{\what\cA}(\bT_{\what{\cA}}(\hat{\btheta})) & = \FDR_{\{1,\cdots, p\}}(\bT(\hat{\btheta})), \\
\FDR_{\what\cA}(\bT_{\what{\cA}}(\btheta^0)) & = \FDR_{\{1,\cdots, p\}}(\bT(\btheta^0)).
\end{align*} 

Meanwhile, by construction $\wtilde{\bX}(\btheta^0)$ satisfies the two properties in Definition \ref{def:knockoff} and is a valid model-X knockoffs matrix. 
Therefore, for any value of the regularization parameter, the LCD statistics $W_j(\btheta^0)$ based on $([\bX,\wtilde{\bX}(\btheta^0)],\by)$ together with Result \ref{thm:candes} ensure the exact FDR control at some target level $q\in(0,1)$. Summarizing this, we obtain that the FDR of IPAD applied to $\bT(\btheta^0)$ is controlled at target level $q$. 

Combining the arguments in the previous two paragraphs, we deduce
\begin{align*}
\FDR_{\what\cA}(\bT_{\what{\cA}}(\btheta^0))=\FDR_{\{1,\cdots, p\}}(\bT(\btheta^0))\leq q.
\end{align*}
Thus the desired results follow automatically if we can prove that $\FDR_{\what{\cA}}(\bT_{\what{\cA}}(\hat{\btheta}))$ is asymptotically close to $\FDR_{\what{\cA}}(\bT_{\what{\cA}}(\btheta^0))$. We next proceed to prove it. 

Recall the definitions of  $\mathbb{I}$ and $\mathbb{I}_{\cA}$ as in \eqref{def: setI}. Define the event
\begin{align*}
\mathcal{E}_{np}=\left\{\bT_{\what{\cA}}(\hat{\btheta}) \in \mathbb{I}_{\what{\cA}}\right\} \cap \left\{\bT_{\what{\cA}}(\btheta^0) \in \mathbb{I}_{\what{\cA}}\right\}.
\end{align*}
Lemma \ref{lem:rateconv} in Section \ref{secB.1} establishes $\hat\btheta \in \bTheta_{np}$ with probability at least $1-O(\pi_{np})$ and $\btheta^0\in \bTheta_{np}$. 
Hence, Lemma \ref{lem:diff_U-EU} in Section \ref{secB.2} guarantees that 
\begin{align}\label{ineq:prob1}
\Pro \left( \mathcal{E}_{np}^c \right) \leq 2\Pro \left( \sup_{|\cA|\leq k,\, \btheta \in \Theta_{np}}\left\|\bT_{\cA}({\btheta}) - \E[\bT_{\cA}(\btheta^0)] \right\|_{\max} > a_{np} \right) 
= O(\pi_{np}),
\end{align}
where $a_{np}=C_1(k^{1/2}+s^{3/2})\tilde{c}_{np}$ for some constant $C_1>0$. 

For a given deterministic set $\cA \subset\{1,\cdots, p\}$, let $\FDP_\cA(\cdot)$ be the FDP function corresponding to $\FDR_\cA(\cdot)$. By the definition of FDP function, we have for any $\bt_1,\bt_2 \in \mathbb{R}^{|\cA|(2|\cA|+3)}$,
\begin{align*}
& \FDP_\cA(\bt_2)-\FDP_\cA(\bt_1) = \frac{|\cS^1\cap S_\cA(\bt_2)|}{| S_\cA(\bt_2)|} - \frac{|\cS^1\cap S_\cA(\bt_1)|}{| S_\cA(\bt_1)|}   \\
&\quad  = \frac{|\cS^1\cap S_\cA(\bt_2)|\cdot (|S_\cA(\bt_1)| - |S_\cA(\bt_2)|)}{| S_\cA(\bt_1)|\cdot | S_\cA(\bt_2)|} + \frac{|\cS^1\cap S_\cA(\bt_2)| - |\cS^1\cap S_\cA(\bt_1)|}{|S_\cA(\bt_1)|}.
\end{align*}
Further, note that
\begin{align*}
&|\cS^1\cap S_\cA(\bt_2)|/|S_\cA(\bt_2)| \leq 1, ~~~~\left|| S_\cA(\bt_2)| - |S_\cA(\bt_1)|\right| \leq \left| S_\cA(\bt_2) \triangle S_\cA(\bt_1) \right|,\\
&\left||\cS^1\cap S_\cA(\bt_2)| - |\cS^1\cap S_\cA(\bt_1)|\right| \leq \left| \{S_\cA(\bt_2) \triangle S_\cA(\bt_1)\}\cap \cS^1 \right|. 
\end{align*}
Combining the results above yields
\begin{align*}
&|\FDP_\cA(\bt_1)-\FDP_\cA(\bt_2)|\\
&\quad \leq  \frac{\big||S_\cA(\bt_1)| - |S_\cA(\bt_2)|\big|}{| S_\cA(\bt_1)|} + \frac{\big| \{S_\cA(\bt_2) \triangle S_\cA(\bt_1)\}\cap \cS^1 \big|}{|S_\cA(\bt_1)|}
\leq 2 \frac{\big| S_\cA(\bt_2) \triangle S_\cA(\bt_1) \big|}{|S_\cA(\bt_1)|} .
\end{align*}
Similarly we have
\begin{align*}
|\FDP_\cA(\bt_1)-\FDP_\cA(\bt_2)|
\leq 2 \frac{\big| S_\cA(\bt_2) \triangle S_\cA(\bt_1) \big|}{|S_\cA(\bt_2)|} .
\end{align*}
Thus it holds that 
\begin{align}
\sup_{|\cA|\leq k}\sup_{\bt_1,\bt_2\in \mathbb{I}_\cA}|\FDP_\cA(\bt_1)-\FDP_\cA(\bt_2)|
&\leq \sup_{|\cA|\leq k}\sup_{\bt_1,\bt_2\in \mathbb{I}_\cA} \frac{\big| S_\cA(\bt_2) \triangle S_\cA(\bt_1) \big|}{|S_\cA(\bt_1)|\wedge |S_\cA(\bt_2)|} \notag\\
&= O(\rho_{np}),  \label{ineq:FDRcontrol2}
\end{align}
where the last two steps are due to Condition \ref{ass:stability}. 
Therefore, \eqref{ineq:prob1} and \eqref{ineq:FDRcontrol2} together with the fact that $\FDP(\cdot) \in [0,1]$ entail that 
\begin{align*}
&\left|\FDR_{\what{\cA}}(\bT_{\what{\cA}}(\hat{\btheta})) -\FDR_{\what{\cA}}(\bT_{\what{\cA}}(\btheta^0)) \right|
=\left|\E \FDP_{\what{\cA}}(\bT_{\what{\cA}}(\hat{\btheta}))-\E\FDP_{\what{\cA}}(\bT_{\what{\cA}}(\btheta^0)) \right| \\
&\quad \leq \E \left|\FDP_{\what{\cA}}(\bT_{\what{\cA}}(\hat{\btheta}))-\FDP_{\what{\cA}}(\bT_{\what{\cA}}(\btheta^0))\right| \\
& \quad \leq \E\left[ \left|\FDP_{\what{\cA}}(\bT_{\what{\cA}}(\hat{\btheta}))-\FDP_{\what{\cA}}(\bT_{\what{\cA}}(\btheta^0))\right| \mid \mathcal E_{np} \right] \Pro\left( \mathcal{E}_{np} \right) +  2 \Pro \left( \mathcal{E}_{np}^c \right)\\
&\quad \leq \sup_{|\cA|\leq k}\sup_{\bt_1,\bt_2\in \mathbb{I}_\cA}\left|\FDP_\cA(\bt_1)-\FDP_\cA(\bt_2)\right| 
+ O(\pi_{np}) \\
&\quad  = O(\rho_{np}) + O(\pi_{np}).
\end{align*}
This completes the proof of Theorem \ref{thm:FDRcontrol}.

\subsection{Proof of Theorem \ref{thm:power}} \label{secA.4}
By the definition of the LCD statistics, we construct the augmented Lasso estimator for each $\btheta\in\Theta_{np}$, which is defined as 
\begin{align}
\hat{\bbeta}^{\textsf{aug}}(\btheta) 
= \arg\min_{\bb\in\mathbb{R}^{2p}} \left\| \by - [\bX, \widetilde{\bX}(\btheta)] \bb \right\|_2^2 + \lambda \|\bb\|_1. \label{auglasso}
\end{align}
The Lasso estimator of regressing $\by$ on only $\bX$ is also given by
\begin{align}
\hat{\bbeta} 
= \arg\min_{\bb\in\mathbb{R}^{p}} \left\| \by - \bX \bb \right\|_2^2 + \lambda \|\bb\|_1, \label{lasso}
\end{align}
where $\lambda=O(n^{-1/2}\log p)$. According to the true model $\cS^0$, the underlying true parameter vector corresponding to $\hat{\bbeta}^{\textsf{aug}}(\btheta)  $ should be given by $\bbeta^\textsf{aug}:=(\bbeta',\bzero')'\in\mathbb{R}^{2p}$ with $\bbeta=(\bbeta_{\cS^0}',\bzero')'\in\mathbb{R}^{p}$ and $|\cS^0|=s$ for any $\btheta\in\Theta_{np}$. 
By Lemma \ref{lem:lasso} in Section \ref{secB.3}, with probability at least $1-O(\pi_{np})$ the Lasso estimators satisfy 
\begin{align*}
\sup_{\btheta\in \Theta_{np}}\left\|\hat{\bbeta}^\textsf{aug}(\btheta)-\bbeta^\textsf{aug}\right\|_1 &= O(s\lambda),\\
\left\|\hat{\bbeta} -\bbeta \right\|_1 &= O(s\lambda),
\end{align*}
where $\lambda =O(n^{-1/2}\log p)$.  

We now prove that under Condition \ref{ass:betamin2}, the power of the augmented Lasso \eqref{auglasso} is bounded from below by $\gamma\in[0,1]$; that is,
\begin{align}
\E \left| \what{\cS}_\textsf{auglasso} \cap \cS^0 \right|/s \geq \gamma, \label{lasspow}
\end{align}
where $\what{\cS}_\textsf{auglasso}=\{j: \hat{\beta}_j^\textsf{aug}(\btheta)\not=0\}$. 
To this end, we first show that with asymptotic probability one, 
\begin{align}\label{eq: lassopower}
|\what{\cS}_\textsf{auglasso}^c \cap \cS^0|/s\leq 1-\gamma.
\end{align}
The key is to use proof by contradiction. Suppose $|\what{\cS}_\textsf{auglasso}^c \cap \cS^0|/s > 1-\gamma$. Then we can see that
\begin{align*}
\sup_{\btheta\in \Theta_{np}}\left\|\hat{\bbeta}^\textsf{aug}(\btheta)-\bbeta^\textsf{aug}\right\|_1
&\geq \sup_{\btheta\in \Theta_{np}}\left\|\hat{\bbeta}_{\what{\cS}_\textsf{auglasso}^c}^\textsf{aug}(\btheta)-\bbeta_{\what{\cS}_\textsf{auglasso}^c}^\textsf{aug}\right\|_1 \\
&=\left\|\bbeta_{\what{\cS}_\textsf{auglasso}^c}^\textsf{aug}\right\|_1
\geq \left\|\bbeta_{\what{\cS}_\textsf{auglasso}^c \cap \cS^0}^\textsf{aug}\right\|_1
> b_{np} sn^{-1/2}\log p,
\end{align*}
where the last step is by Condition \ref{ass:betamin2}. 
However, by Lemma  \ref{lem:lasso} with probability at least $1-O(\pi_{np})$, the left hand side above  is bounded from above by $O(s\lambda)$ with $\lambda=O(n^{-1/2}\log p)$. These two results contradict with each other since $b_{np}\to \infty$. Hence \eqref{eq: lassopower} is proved. Therefore, the result in \eqref{lasspow} follows immediately since 
$|\what{\cS}_\textsf{auglasso} \cap \cS^0| = s - |\what{\cS}_\textsf{auglasso}^c \cap \cS^0|$ and 
\begin{align*}
&\E \left| \what{\cS}_\textsf{auglasso} \cap \cS^0 \right|/s  \geq \gamma \Pro \left(|\what{\cS}_\textsf{auglasso} \cap \cS^0|/s > \gamma \right) \\
& = \gamma \Pro \left(|\what{\cS}_\textsf{auglasso}^c \cap \cS^0|/s\leq 1-\gamma \right) = \gamma (1-O(\pi_{np})).
\end{align*}

Using the same argument, we can show that the power of the Lasso \eqref{lasso} is also bounded from below by $\gamma(1-O(\pi_{np}))$ under Condition \ref{ass:betamin2}. That is, we have
\begin{align*}
\E \left| \what{\cS}_\textsf{lasso} \cap \cS^0 \right|/s \geq \gamma (1-O(\pi_{np})),
\end{align*}
where $\what{\cS}_\textsf{lasso}=\{j: \hat{\beta}_j\not=0\}$. 

Next we show that our knockoffs procedure has at least the same power as the augmented Lasso and hence the Lasso itself. Namely, we prove
\begin{align}
\E\left| \what{\cS} \cap \cS^0 \right|/s \geq \gamma \label{kopow}
\end{align}
with threshold $T_2$. Note that the same argument is still valid for $ T_1 $. 
Let $|W_{(1)}|\geq \dots \geq |W_{(p)}|$ and define $j^*$ as  $|W_{(j^*)}|=T_2$. 
Then by the definition of $T_2$, it holds that $-T_2 < W_{j^*+1}\leq 0$. Here we have assumed that there are no ties on the magnitudes of $W_j$'s which should be a reasonable assumption considering the continuity of the Lasso solution. 
As in the proof of Theorem 3 in \cite{RANK}, it is sufficient to consider the following two cases.

\textbf{Case 1.} Consider the case of $-T_2<W_{(j^*+1)}<0$. In this case, from the definition of threshold $ T_2 $ we have
\begin{equation*}
\frac{ 2+ |\{ j : W_{(j)} \leq -T_2 \}| }{ | \{j : W_{(j)} \geq T_2 \}| } > q.
\end{equation*}
Using the same argument as in Lemma 6 of \cite{RANK} together with Lemma \ref{lem:lasso}, we can prove from Condition \ref{ass:active} that $|\what\cS|\geq C_2s$ with probability at least $1-O(\pi_{np})$. 
This leads to $ | \{j : W_{(j)} \leq -T_2 \}| >  C_2qs -2 $ with the same probability. Now from the same argument as in A.5 of \cite{RANK}, we can obtain  $T_2=O(\lambda)$. On the other hand, Lemma \ref{lem:lasso} and some algebra establish that 
\begin{align}
O(s\lambda) &= \|\hat{\bbeta}^\textsf{aug}(\hat{\btheta})-\bbeta^\textsf{aug}\|_1 
= \sum_{j=1}^p |\hat{\beta}_j^\textsf{aug}(\hat{\btheta})-\beta_j| + \sum_{j=1}^p|\hat{\beta}_{j+p}^\textsf{aug}(\hat{\btheta})| \notag\\
&= \sum_{j\in \what{\cS}\cap\cS^0} |\hat{\beta}_j^\textsf{aug}(\hat{\btheta})-\beta_j| 
+ \sum_{j\in \cS^1} |\hat{\beta}_j^\textsf{aug}(\hat{\btheta})| \notag\\
&\qquad\qquad+ \sum_{j\in \what{\cS}^c\cap\cS^0} |\hat{\beta}_j^\textsf{aug}(\hat{\btheta})-\beta_j| 
+ \sum_{j=1}^p|\hat{\beta}_{j+p}^\textsf{aug}(\hat{\btheta})|. \label{pow:ineq01}
\end{align}
We then consider the lower bound of the last term in \eqref{pow:ineq01}. 
For any $j\in \what{S}^c$, it holds that $|\hat{\beta}_{j+p}^\textsf{aug}(\hat{\btheta})|>|\hat{\beta}_j^\textsf{aug}(\hat{\btheta})|-T_2$. Hence we obtain
\begin{align}
\sum_{j=1}^p|\hat{\beta}_{j+p}^\textsf{aug}(\hat{\btheta})|
&\geq \sum_{j\in\what{\cS}\cap\cS^0}|\hat{\beta}_{j+p}^\textsf{aug}(\hat{\btheta})|
+\sum_{j\in\what{\cS}^c\cap\cS^0}|\hat{\beta}_{j+p}^\textsf{aug}(\hat{\btheta})| \notag\\
&\geq \sum_{j\in\what{\cS}\cap\cS^0}|\hat{\beta}_{j+p}^\textsf{aug}(\hat{\btheta})|
+\sum_{j\in\what{\cS}^c\cap\cS^0}|\hat{\beta}_{j}^\textsf{aug}(\hat{\btheta})|-T_2|\what{\cS}^c\cap\cS^0|. \label{pow:ineq02}
\end{align}
Plugging \eqref{pow:ineq02} into \eqref{pow:ineq01} and applying the triangle inequality yield
\begin{align*}
O(s\lambda) &\geq \sum_{j\in \what{\cS}\cap\cS^0} |\hat{\beta}_j^\textsf{aug}(\hat{\btheta})-\beta_j| 
+ \sum_{j\in \cS^1} |\hat{\beta}_j^\textsf{aug}(\hat{\btheta})| 
+ \sum_{j\in \what{\cS}^c\cap\cS^0} |\hat{\beta}_j^\textsf{aug}(\hat{\btheta})-\beta_j| \notag\\
&\qquad \qquad + \sum_{j\in\what{\cS}\cap\cS^0}|\hat{\beta}_{j+p}^\textsf{aug}(\hat{\btheta})|
+\sum_{j\in\what{\cS}^c\cap\cS^0}|\hat{\beta}_{j}^\textsf{aug}(\hat{\btheta})|-T_2|\what{\cS}^c\cap\cS^0| \notag\\
&\geq \sum_{j\in \what{\cS}\cap\cS^0} |\hat{\beta}_j^\textsf{aug}(\hat{\btheta})-\beta_j| 
+ \sum_{j\in \cS^1} |\hat{\beta}_j^\textsf{aug}(\hat{\btheta})| \notag\\
&\qquad\qquad +\sum_{j\in\what{\cS}^c\cap\cS^0}|\beta_{j}| + \sum_{j\in\what{\cS}\cap\cS^0}|\hat{\beta}_{j+p}^\textsf{aug}(\hat{\btheta})|
-T_2|\what{\cS}^c\cap\cS^0| \notag\\
&\geq \sum_{j\in\what{\cS}^c\cap\cS^0}|\beta_{j}|
-T_2|\what{\cS}^c\cap\cS^0|
=\|\bbeta_{\what{\cS}^c\cap\cS^0}\|_1
-T_2|\what{\cS}^c\cap\cS^0|. 
\end{align*}
Since $T_2|\what{\cS}^c\cap\cS^0|=O(s\lambda)$ for $\lambda=O(n^{-1/2}\log p)$ due to the discussion above, we consequently obtain
\begin{align}
\|\bbeta_{\what{\cS}^c\cap\cS^0}\|_1 =O(sn^{-1/2}\log p). \label{pow:ineq03}
\end{align}
Suppose $|\what{\cS}^c\cap\cS^0|/s>1-\gamma$. Then Condition \ref{ass:betamin2} gives $\|\bbeta_{\what{\cS}^c\cap\cS^0}\|_1>b_{np} sn^{-1/2}\log p$ for some positive diverging sequence $b_{np}$; this contradicts with \eqref{pow:ineq03}. Thus we obtain $|\what{\cS}^c\cap\cS^0|/s\leq 1-\gamma$ with asymptotic probability one, which leads to \eqref{kopow} by taking expectation. 

\textbf{Case 2.} Consider the case of $W_{(j^*+1)}=0$. In this case, by the definition of threshold $ T_2 $
\begin{equation}
\frac{ 1+ |\{j: W_{(j)} < 0\}| }{ |\{j: W_{(j)} > 0\}| } \leq q.
\end{equation}
If $ | \{j: W_{(j)} < 0\} | > C_3 s $ for some constant $ C_3 >0 $, then from the same argument as in A.5 of \cite{RANK}, we can obtain  $T_2=O(\lambda)$, and the rest of the proof is the same as in Case 1.
On the other hand, if $ | \{j: W_{(j)} < 0\} | \leq o(s) $ we have
\begin{align*}
	| \{j: W_{(j)} \neq 0\} \cap \cS^0 | & =  | \{j: W_{(j)} >0\} \cap \cS^0 | + | \{j: W_{(j)} < 0\} \cap \cS^0 | \\
	&\leq |\what{\cS} \cap\cS^0| + o(s).
\end{align*}
Now note that $ | \{j: W_{(j)} \neq 0 \}| \geq | \{j : |\hat{\beta}_j^{\textsf{aug}}| \neq 0,~ j = 1, \cdots, p\}|$. Then we can see that with asymptotic probability one,
\begin{align*}
	| \{j: W_{(j)} \neq 0\} \cap \cS^0 | 
	& \geq | \{ j : \hat{\beta}_j^{\textsf{aug}} \neq 0,~ j = 1, \cdots, p\}  \cap \cS^0 |\\
	&= |\what{\cS}_\textsf{auglasso} \cap \cS^0| .\\
	& \geq \gamma s (1 - o(1) ).
\end{align*}
Consequently, we obtain $|\what{\cS} \cap\cS^0|/s\geq \gamma (1 - o(1))$, which leads to \eqref{kopow} by taking expectation. Combining these  two cases concludes the proof of Theorem \ref{thm:power}.

\renewcommand{\theequation}{B.\arabic{equation}}
\setcounter{equation}{0}

\section{Some key lemmas and their proofs} \label{secB}

\subsection{Lemma \ref{lem:rateconv} and its proof} \label{secB.1}
\begin{lem}\label{lem:rateconv}
	Assume that Conditions \ref{ass:fac}--\ref{ass:eigensepa} hold. Then with probability at least $1-O(\pi_{np})$, the estimator $\hat{\btheta}=(\vect(\what{\bC})',\hat{\bbbeta}')'$ lies in the shrinking  set given by
	\begin{align*}
	\Theta_{np} &=\left\{\btheta=(\vect(\bC)',\bbbeta')':\left\|\bC-\bC^0\right\|_{\max} + \left\|\bbbeta-\bbbeta^0\right\|_{\max} 
	\leq  O(c_{np}) \right\},
	\end{align*}
	where $c_{np}=(n^{-1}\log p)^{1/2}+(p^{-1}\log n)^{1/2}$ and $\pi_{np}=p^{-\nu}+n^{-\nu}$.
\end{lem}

\noindent\textit{Proof}. We divide the proof into two parts.  
We prove the bound for $\|\what{\bC}-\bC^0\|_{\max}$ in Part 1 and then for $\|\hat{\bbbeta}-\bbbeta^0 \|_{\max}$ in Part 2. 

\textbf{Part 1.}
Note that $\|\what{\bC}-\bC^0\|_{\max}=\max_{i,j}|\hat{c}_{ij}-c_{ij}^0|$, where the maximum is taken over $i\in\{1,\dots,n\}$ and $j\in\{1,\dots,p\}$. 
We write $\bff_i^*=\bH'\bff_i^0$ and $\blambda_j^*=\bH^{-1}\blambda_j^0$ with rotation matrix $\bH$ defined in Lemma \ref{lem:HVbound} in Section \ref{secC.1}. 
From the definition of $c_{ij}$, it holds that 
\begin{align*}
\hat{c}_{ij}-c_{ij}^0 = (\hat{\bff}_i-\bff_i^*)'\blambda_j^* + \hat{\bff}_i'(\hat{\blambda}_j-\blambda_j^*).
\end{align*}
From Lemma \ref{lem:HVbound}, we can assume  $\|\bH\|_2+\|\bH^{-1}\|_2+\|\bV\|_2+\|\bV^{-1}\|_2\lesssim 1$, which occurs with probability at least $1-O(p^{-\nu}) $. 
We also have $\max_{i\in\{1,\dots,n\}}\|\hat{\bff}_i\|_2^2\lesssim 1$ a.s.\ by the assumed restriction $\hat{\bF}'\hat{\bF}/n=\bI_r$ as mentioned on p.213 of \cite{bai2002determining}. Hence, 
the triangle and Cauchy--Schwarz inequalities with Conditions \ref{ass:fac} and \ref{ass:floa} give
\begin{align}
\max_{i,j}|\hat{c}_{ij}-c_{ij}| 
&\leq \max_{i}\|\hat{\bff}_i-\bff_i^*\|_2 \max_{j}\|\blambda_j^*\|_2 
+ \max_{i}\|\hat{\bff}_i\|_2 \max_{j}\|\hat{\blambda}_j-\blambda_j^*\|_2 \notag \\
&\lesssim \max_{i}\|\hat{\bff}_i-\bff_i^*\|_2 
+ \max_{j}\|\hat{\blambda}_j-\blambda_j^*\|_2. \label{eq:lempb2:01}
\end{align}
Then it is sufficient to derive upper bounds for $\max_{i}\|\hat{\bff}_i-\bff_i^*\|_2$ and $\max_{j}\|\hat{\blambda}_j-\blambda_j^*\|_2$ that hold with high probability. 
Using the decomposition of A.1 in \cite{bai2003} along with taking maximum over $i,\ell \in\{1,\dots,n\}$, we can deduce
\begin{align}
&\max_i\|\hat{\bff}_i-\bff_i^*\|_2 \notag\\
&\qquad \leq \|\bV^{-1}\|_2 \max_i \left( 
(\sigma_e^2/n)\|\hat{\bff}_i\|_2
+n^{-1}\sum_{\ell=1}^n \|\hat{\bff}_\ell\|_2 \left|p^{-1}\sum_{j=1}^p\left(e_{\ell j}e_{ij}-\E[e_{\ell j}e_{ij}] \right)\right| \right. \notag \\
&\qquad \qquad + n^{-1} \left. \sum_{\ell=1}^n \|\hat{\bff}_\ell {\bff_\ell^0}'\|_2 \left\|p^{-1}\sum_{j=1}^p \blambda_j^0 e_{ij} \right\|_2
+ n^{-1}\sum_{\ell=1}^n \|\hat{\bff}_\ell {\bff_i^0}'\|_2 \left\|p^{-1}\sum_{j=1}^p \blambda_j^0 e_{\ell j} \right\|_2
\right) \notag \\
&\qquad \lesssim O(n^{-1}) + \max_{i,\ell} \left|p^{-1}\sum_{j=1}^p\left(e_{\ell j}e_{ij}-\E[e_{\ell j}e_{ij}] \right)\right| 
+ \max_{i}\left\|p^{-1}\sum_{j=1}^p \blambda_j^0 e_{ij} \right\|_2 \notag \\
&\qquad \lesssim O(n^{-1}) + R_1 + R_2,  \label{eq:lempb2:02}
\end{align}
where we have used the boundedness of $\|\hat\bff_\ell\|_2$ discussed above and $\|\bff_\ell^0\|_2\leq r^{1/2}\|\bff_\ell^0\|_{\max}\lesssim 1$ in Condition \ref{ass:fac} for the second inequality, and defined $T_1=\max_{i,\ell} \left|p^{-1}\sum_{j=1}^p\left(e_{\ell j}e_{ij}-\E[e_{\ell j}e_{ij}] \right)\right|$ and $T_2=\max_{i,k}\left|p^{-1}\sum_{j=1}^p \lambda_{jk}^0 e_{ij} \right|$. 
Similarly, the expression on p.165 of \cite{bai2003} with 
taking maximum over $i\in\{1,\dots,n\}$ and $j\in\{1,\dots,p\}$ leads to 
\begin{align}
&\max_{j}\|\hat{\blambda}_j-\blambda_j^*\|_2 \notag\\
&\qquad\leq \|\bH\|_2\max_{j}\left\|n^{-1}\sum_{i=1}^n\bff_i^0e_{ij} \right\|_2 
+ \left\| n^{-1}\sum_{i=1}^n\hat\bff_i(\hat{\bff}_i-\bff_i^*)'\right\|_2\left\|\bH^{-1}\right\|_2\max_{j}\left\|\blambda_j^0\right\|_2 \notag\\
&\qquad \qquad+ \max_{j}\left\| n^{-1}\sum_{i=1}^n(\hat{\bff}_i-\bff_i^*)e_{ij}\right\|_2 \notag\\
&\qquad\lesssim \max_{j}\left\|n^{-1}\sum_{i=1}^n \bff_i^0 e_{ij} \right\|_2
+ \max_{i}\left\|\hat{\bff}_i-\bff_i^*\right\|_2 
+ \max_{i} \left\|\hat{\bff}_i-\bff_i^*\right\|_2 \max_{j} \left(n^{-1}\sum_{i=1}^n e_{ij}^2\right)^{1/2} \notag \\
&\qquad = R_3 + \max_{i}\|\hat{\bff}_i-\bff_i^*\|_2 (1+R_4),\label{eq:lempb2:03}
\end{align}
where $R_3=\max_{j,k}\left|n^{-1}\sum_{i=1}^n f_{ik}^0 e_{ij} \right|_2$ 
and $R_4=\max_{j} \left(n^{-1}\sum_{i=1}^n e_{ij}^2\right)^{1/2}$, and the Cauchy--Schwarz inequality has been used to obtain the second inequality. To evaluate $R_4$, we note that
\begin{align*}
R_4^2 \leq \max_{j}\E e_{ij}^2 + \max_{j} \left|n^{-1}\sum_{i=1}^n \left(e_{ij}^2 - \E e_{ij}^2\right) \right|.
\end{align*}
The first term is bounded by $2C_e^2$. For the second term, Lemma \ref{lem:ineq}(a) in Section \ref{secC.2} with $p$ replaced by $n$ and the union bound give
\begin{align*}
&\Pro \left( \max_{j} \left|n^{-1}\sum_{i=1}^n \left(e_{ij}^2 - \E e_{ij}^2\right) \right| > u \right) 
\leq p\max_{j} \Pro \left( \left|n^{-1}\sum_{i=1}^n \left(e_{ij}^2 - \E e_{ij}^2\right) \right| > u \right) \\
&\qquad \leq 2p\exp(-nu^2/C)
\end{align*}
for all $0\leq u\leq c$. Thus putting $u=(C(\nu+1) n^{-1}\log p)^{1/2}$ and using condition $c_{np}\leq c/(r^2M^2C(\nu+2))^{1/2}$, we obtain $R_4^2=O(1)+O((n^{-1}\log p)^{1/2})=O(1)$ with probability at least $1-O(p^{-\nu})$. 
This together with the observation from \eqref{eq:lempb2:01}--\eqref{eq:lempb2:03} yields 
\begin{align*}
\max_{i,j}|\hat{c}_{ij}-c_{ij}^0| 
&\lesssim R_3 + \left\{R_1 + R_2 + O(n^{-1}) \right\}(1+R_4) \\
&\lesssim R_1 + R_2 + R_3 + O(n^{-1}).
\end{align*}
Hence the convergence rate of $\max_{i,j}|\hat{c}_{ij}-c_{ij}^0| $ is determined by the slowest term out of $R_1$, $R_2$, $R_3$, and $O(n^{-1})$. We evaluate these terms by Lemma \ref{lem:ineq} in Section \ref{secC.2} and the union bound with condition $c_{np}\leq c/(r^2M^2C(\nu+2))^{1/2}$ as above. 
First for $R_1$, Lemma \ref{lem:ineq}(a) by letting $u_1=(C(\nu+2) p^{-1}\log n )^{1/2}$ results in 
\begin{align*}
\Pro\left( R_1 > u_1 \right) 
\leq 2 n^2 \exp \left\{ -p(\nu+2) p^{-1} \log n  \right\} 
=O(n^{-\nu}). 
\end{align*}
Next for $R_2$, Lemma \ref{lem:ineq}(c) with $u_2=(2(\nu+1) p^{-1}\log n)^{1/2}$ gives
\begin{align*}
\Pro\left( R_2 > u_2 \right) 
&\leq 2rn \exp \left\{ - p(\nu+1) p^{-1}\log n \right\} = O(n^{-\nu}).
\end{align*}
Finally for $R_3$, Lemma \ref{lem:ineq}(b) with putting $u_3=(C(\nu+1) n^{-1}\log p)^{1/2}$ leads to
\begin{align*}
\Pro\left( R_3 > u_3 \right) 
&\leq 2rp \exp \left\{ - n(\nu+1) n^{-1}\log p \right\} = O(p^{-\nu}). 
\end{align*}
Consequently, we obtain the first result $\|\what{\bC}-\bC^0\|_{\max}=O(c_{np})$, which holds with probability at least $1-O(\pi_{np})$. 

\textbf{Part 2.}
Next we derive the convergence rate of $\hat{\bbbeta}$. It is sufficient to prove only the case 
when $\bbbeta^0$ is a scalar (so that we write $\bbbeta^0=\eta_1^0$) since dimensionality $m$ is fixed and $\eta_k^0$'s share the identical property thanks to Condition \ref{ass:err_new}. 
Recall notation $\mathbb{E}_{np}e^k=(np)^{-1}\sum_{i,j} e_{ij}^k$. Letting $\delta_{ij}=c_{ij}^0-\hat{c}_{ij}$, we have $\hat{e}_{ij}=x_{ij}-\hat{c}_{ij}=e_{ij}+\delta_{ij}$. For an arbitrary fixed $k\in\{1,\dots,m\}$, the binomial expansion entails 
\begin{align}
\left| \mathbb{E}_{np} \hat{e}^k - \E e^k \right|
&= \left| \mathbb{E}_{np} (e+\delta)^k - \E e^k \right| \notag\\
&= \left| \mathbb{E}_{np}(e^k-\E e^k) + \mathbb{E}_{np} \sum_{\ell=0}^{k-1} \binom{k}{\ell} e^\ell\delta^{k-\ell} \right| \notag\\
&\leq \left| \mathbb{E}_{np}(e^k-\E e^k) \right| + \sum_{\ell=0}^{k-1} \binom{k}{\ell}\max_{i,j}|\delta_{ij}|^{k-\ell}\mathbb{E}_{np} |e|^\ell \notag\\
&\lesssim \left| \mathbb{E}_{np}(e^k-\E e^k) \right| + O\left(\max_{i,j}|\delta_{ij}| \right) \sum_{\ell=0}^{k-1} \mathbb{E}_{np} |e|^\ell.\label{ineq:eta2}
\end{align}
For all $k\in\{1,\dots,m\}$, the strong law of large numbers with Theorem 2.5.7 in \cite{Durrett2010} entails 
$|\mathbb{E}_{np}e^{k}- \E e^{k} | = o((np)^{-1/2}\log(np))$ a.s.\ under Condition \ref{ass:err_new}. Furthermore, the second term of \eqref{ineq:eta2} is $O(c_{np})$ with probability at least $1-O(\pi_{np})$ from Part 1 and the same law of large numbers. Consequently, we obtain
\begin{align*}
\left| \mathbb{E}_{np} \hat{e}^k - \E e^k \right|
\lesssim c_{np}.
\end{align*}
Therefore by the construction of $\hat{\eta}_1$ and local Lipschitz continuity of $h_1$ in Condition \ref{ass:err_new}, we see that
\begin{align*}
\left|\hat{\eta}_1 -\eta_1^0 \right| 
&= \left| h_1\left( \mathbb{E}_{np} \hat{e}, \dots, \mathbb{E}_{np} \hat{e}^{m} \right) - h_1\left(\E e, \dots, \E e^{m} \right) \right| \notag\\
&\lesssim \max_{k\in\{1,\dots,m\}} \left| \mathbb{E}_{np} \hat{e}^k - \E e^k \right|
\end{align*}
with probability at least $1-O(\pi_{np})$. This completes the proof of Lemma \ref{lem:rateconv}.

\subsection{Lemma \ref{lem:diff_U-EU} and its proof} \label{secB.2}
\begin{lem}\label{lem:diff_U-EU}
	Assume that Conditions \ref{ass:regerr}--\ref{ass:err_new} hold. Then with probability at least $1-O(\pi_{np})$, the following statements hold
	\begin{align*}
	(a)&~~\sup_{|\cA|\leq k,\, \btheta \in \Theta_{np}}\left\|\bU_{\cA}({\btheta}) - \E[\bU_{\cA}(\btheta^0)] \right\|_{\max} = O\left(k^{1/2}\tilde{c}_{np}\right), \\
	(b)&~~\sup_{|\cA|\leq k,\, \btheta \in \Theta_{np}}\left\|\bv_{\cA}({\btheta}) - \E[\bv_{\cA}(\btheta^0)] \right\|_{\max} = O\left(s^{3/2} \tilde{c}_{np}\right),
	\end{align*}
	where $\Theta_{np}$ was defined in Lemma \ref{lem:rateconv} and $\tilde{c}_{np} = n^{-1/2}\log p + p^{-1/2}\log n$. Consequently, we have 
	\begin{align*}
	\sup_{|\cA|\leq k,\, \btheta \in \Theta_{np}}\left\|\bT_{\cA}({\btheta}) - \E[\bT_{\cA}(\btheta^0)] \right\|_{\max} = O\left( \left( k^{1/2} + s^{3/2}\right)\tilde{c}_{np}\right).
	\end{align*}
\end{lem}

\noindent\textit{Proof}. To complete the proof of $(a)$, we verify the following
\begin{align*}
(a\textendash i)&~~\sup_{|\cA|\leq k,\, \btheta \in \Theta_{np}}\left\|\bU_{\cA}({\btheta}) - \bU_{\cA}(\btheta^0) \right\|_{\max} 
\lesssim k^{1/2}\tilde{c}_{np}, \\
(a\textendash ii)&~~\left\|\bU(\btheta^0) - \E[\bU(\btheta^0)] \right\|_{\max} \lesssim (n^{-1}\log p)^{1/2}.
\end{align*}
From $(a\textendash i)$ and $(a\textendash ii)$, we can conclude that
\begin{align*}
&\sup_{|\cA|\leq k,\, \btheta \in \Theta_{np}}\left\|\bU_{\cA}({\btheta}) - \E[\bU_{\cA}(\btheta^0)] \right\|_{\max} \\
&\qquad\leq \sup_{|\cA|\leq k,\, \btheta \in \Theta_{np}}\left\|\bU_{\cA}({\btheta}) - \bU_{\cA}(\btheta^0) \right\|_{\max}
+ \sup_{|\cA|\leq k}\left\|\bU_{\cA}(\btheta^0) - \E[\bU_{\cA}(\btheta^0)] \right\|_{\max} \\
&\qquad\leq \sup_{|\cA|\leq k,\, \btheta \in \Theta_{np}}\left\|\bU_{\cA}({\btheta}) - \bU_{\cA}(\btheta^0) \right\|_{\max}
+ \left\|\bU(\btheta^0) - \E[\bU(\btheta^0)] \right\|_{\max} \\
&\qquad\lesssim k^{1/2}\tilde{c}_{np},
\end{align*}
which yields result (a). 

We begin with showing $(a\textendash i)$; this is the uniform extension of Lemma \ref{lem:diffsUV_max}(a) in Section \ref{secC.3} over $|\cA|\leq k$. In fact, the proof is almost the same, with the only difference that bound \eqref{bound:E1} should be replaced with the bound derived in Lemma \ref{lem:maxeigv}(c); that is,  
\begin{align}\label{bound:E2}
\max_{|\cA|\leq k} \left\|n^{-1/2}\bE_\cA \right\|_2 \lesssim 1\vee \left(kn^{-1}\log p \right)^{1/2},
\end{align}
which holds with probability at least $1-O(p^{-\nu})$. Notice that $\left(kn^{-1}\log p \right)^{1/2}\leq \log^{1/2} p$. Therefore, even if we use \eqref{bound:E2} instead of \eqref{bound:E1} in the proof of Lemma \ref{lem:diffsUV_max}(a) we can still derive the same convergence rate $k^{1/2}\tilde{c}_{np}$ as in Lemma \ref{lem:diffsUV_max}(a), and hence $(a\textendash i)$ holds with probability at least $1-O(\pi_{np})$. 

For $(a\textendash ii)$, we see that
\begin{align}
&\left\|\bU(\btheta^0) - \E[\bU(\btheta^0)] \right\|_{\max} 
\leq \left\|n^{-1}\bX'\bX - \E[ n^{-1}\bX'\bX ] \right\|_{\max} \notag\\
&\qquad + \left\|n^{-1}\wtilde{\bX}(\btheta^0)'\wtilde{\bX}(\btheta^0) - \E[ n^{-1}\wtilde{\bX}(\btheta^0)'\wtilde{\bX}(\btheta^0)] \right\|_{\max} \notag\\
&\qquad + 2\left\|n^{-1}\bX'\wtilde{\bX}(\btheta^0) - \E[ n^{-1}\bX'\wtilde{\bX}(\btheta^0)] \right\|_{\max} 
=: W_1 + W_2 + 2W_3.\label{ineq:004}
\end{align}
We derive the bounds for each of these terms. First, $W_1$ is bounded as
\begin{align*}
W_{1} &\leq \left\|n^{-1}{\bC^0}'\bC^0 - \E[ n^{-1}{\bC^0}'\bC^0 ] \right\|_{\max}
+ \left\|n^{-1}\bE'\bE - \E n^{-1}\bE'\bE\right\|_{\max} 
+2 \left\|n^{-1}\bE'{\bC^0} \right\|_{\max} \\
&=:W_{1,1}+W_{1,2}+W_{1,3}. 
\end{align*}
Under Condition \ref{ass:floa}, we deduce
\begin{align*}
W_{1,1} &=  \max_{j,\ell\in\{1,\dots,p\}}\left|\sum_{k,m=1}^r\lambda_{jk}^0\lambda_{\ell m}^0 n^{-1}\sum_{i=1}^n\left(f_{ik}^0f_{im}^0-\E f_{ik}^0f_{im}^0\right)\right| \\
&\leq rM^2\max_{j,\ell\in\{1,\dots,p\}}\left| n^{-1}\sum_{i=1}^n\left(f_{ik}^0f_{im}^0-\E f_{ik}^0f_{im}^0\right)\right| .
\end{align*}
From Lemma \ref{lem:ineq}(d) with Condition \ref{ass:fac} and the union bound, we have 
\begin{align*}
&\Pro \left(\max_{j,\ell\in\{1,\dots,p\}}\left| n^{-1}\sum_{i=1}^n\left(f_{ik}^0f_{im}^0-\E f_{ik}^0f_{im}^0\right)\right| > u\right) \\
&\qquad \leq p^2 \max_{j,\ell\in\{1,\dots,p\}}\Pro \left(\left| n^{-1}\sum_{i=1}^n\left(f_{ik}^0f_{im}^0-\E f_{ik}^0f_{im}^0\right)\right| > u\right) 
\leq 2p^2 \exp\left( -nu^2/C \right).
\end{align*}
Hence, letting $u=(C(\nu+2)n^{-1}\log p)^{1/2}$ above yields the bound $W_{1,1}\lesssim (n^{-1}\log p)^{1/2}$ with probability at least $1-O(p^{-\nu})$. 
Next for $W_{1,2}$, we can find from Lemma \ref{lem:ineq}(a) with $p$ replaced by $n$ and the union bound that 
\begin{align*}
&\Pro \left( \left\|n^{-1}\bE'\bE - \E n^{-1}\bE'\bE\right\|_{\max} > u \right)
\leq p^2\max_{j,\ell}\Pro \left( \left|n^{-1}\sum_{i=1}^n\left( e_{ij}e_{i\ell} - \E e_{ij}e_{i\ell} \right)\right| > u \right) \\
&\qquad \leq 2p^2 \exp \left( -nu^2/C \right).
\end{align*}
Letting $u=(C(\nu+2)n^{-1}\log p)^{1/2}$ and using $n^{-1}\log p \leq c^2/(C(\nu+2))$, we obtain $W_{1,2}\lesssim (n^{-1}\log p)^{1/2}$ with probability at least $1-O(p^{-\nu})$. Next for $W_{1,3}$, the union bound gives
\begin{align*}
&\Pro \left( \left\| n^{-1}\bE'\bF^0 {\bLambda^0}'\right\|_{\max} > u \right)
=\Pro \left( \max_{j,\ell \in \{1,\dots,p\}}\left|n^{-1}\sum_{k=1}^r\sum_{i=1}^{n} e_{ij}f_{ik}^0\lambda_{\ell k}^0 \right| > u \right) \\
&\qquad\leq \Pro \left( r\max_{j,\ell \in \{1,\dots,p\}} \max_{k \in \{1,\dots,r\}} \left|n^{-1}\sum_{i=1}^{n} e_{ij}f_{ik}^0\right| \left|\lambda_{\ell k}^0 \right| > u \right) \\
&\qquad\leq rp\max_{k\in \{1,\dots,r\}}\max_{j\in \{1,\dots,p\}} \Pro \left(  \left|n^{-1}\sum_{i=1}^{n} e_{ij}f_{ik}^0 \right| > u/(rM) \right).
\end{align*}
Lemma \ref{lem:ineq}(b) states that for all $0\leq u/(rM)\leq c/(rM)$ it holds that 
\begin{align*}
\Pro \left( \left|n^{-1}\sum_{i=1}^{n} e_{ij}f_{ik}^0 \right| > u/(rM) \right)
\leq 2\exp\left\{ -nu^2/(Cr^2M^2) \right\}.
\end{align*}
Therefore, if we put $u= rM(C(\nu+1)n^{-1}\log p)^{1/2}$ using $n^{-1}\log p \leq c^2/(r^2M^2C(\nu+1))$, the upper bound of the probability is further bounded by $2rp^{-\nu}$. Thus we obtain $W_{13}\lesssim (n^{-1}\log p)^{1/2}$ with probability at least  $1- O(p^{-\nu})$.
Consequently, the bound of $W_1$ is 
\begin{align*}
W_1 \leq W_{1,1}+W_{1,2}+W_{1,3} \lesssim (n^{-1}\log p)^{1/2}
\end{align*}
with probability at least  $1- O(p^{-\nu})$. Note that we have the same result for $W_2$ since it has the same distribution as $W_1$. Finally, $W_3$ is bounded as
\begin{align*}
&W_3 
\leq  \left\|n^{-1}{\bC^0}'\bC^0 - \E[ n^{-1}{\bC^0}'\bC^0] \right\|_{\max} 
+ \left\|n^{-1}\bE'\bE_{\bbbeta^0}\right\|_{\max} \\
&\qquad+ \left\|n^{-1}\bE'\bC^0 \right\|_{\max} 
+ \left\|n^{-1}\bE_{\bbbeta^0}'\bC^0 \right\|_{\max} 
=: W_{1,1} + W_{3,1} + W_{1,3} + W_{3,2}.
\end{align*}
The upper bound of $W_{3,1}$ turns out to be $O((n^{-1}\log p)^{1/2})$ that holds with probability at least $1-O(p^{-\nu})$. 
We check this claim. Using the union bound and the inequality of Lemma \ref{lem:ineq}(a) with $p$ replaced by $n$ and putting 
$u=(C(\nu+2) n^{-1}\log p )^{1/2}$ yield
\begin{align*}
\Pro \left( \left\|n^{-1}\bE'\bE_{\bbbeta^0}\right\|_{\max} > u \right) 
\leq p^2 \max_{j,\ell }\Pro \left( \left|n^{-1}\sum_{i=1}^{n}\left(e_{ij}e_{\bbbeta^0,i\ell} \right) \right| > u \right) 
\leq 2 p^{-\nu}.
\end{align*}
Finally, $W_{3,2}$ is found to have the same bound as $W_{1,3}$ because $\bE_{\bbbeta^0}$ is an independent copy of $\bE$. 
Consequently, with probability at least $1-O(p^{-\nu})$, we obtain
\begin{align*}
\left\| \bU(\btheta^0) - \E [\bU(\btheta^0)] \right\|_{\max} 
\lesssim (n^{-1} \log p)^{1/2}.
\end{align*}
This completes the proof of $(a)$ since $p^{-\nu}/\pi_{np}=O(1)$. 

Next we show $(b)$ by verifying the following
\begin{align*}
(b\textendash i)&~~\sup_{|\cA|\leq k,\, \btheta \in \Theta_{np}}\left\|\bv_{\cA}({\btheta}) - \bv_{\cA}(\btheta^0) \right\|_{\max} 
\lesssim s^{3/2}\tilde{c}_{np}, \\
(b\textendash ii)&~~\left\|\bv(\btheta^0) - \E[\bv(\btheta^0)] \right\|_{\max} \lesssim s(n^{-1} \log p)^{1/2}.
\end{align*}
Similar to the proof of $(a)$, we need to modify the proof of Lemma \ref{lem:diffsUV_max}(b) in Section \ref{secC.3} for obtaining the uniform bound with respect to $\cA$, but the obtained result is already uniform over the choice of $\cA$. Thus the same upper bound holds and $(b\textendash i)$ follows. 
Next we show $(b\textendash ii)$. It holds that 
\begin{align*}
&\left\|\bv(\btheta^0) - \E \bv(\btheta^0) \right\|_{\max} \\
&\quad \leq \left\| n^{-1}\bX'\by - \E n^{-1}\bX'\by \right\|_{\max} + \left\| n^{-1}\wtilde{\bX}(\btheta^0)'\by - \E n^{-1}\wtilde{\bX}(\btheta^0)'\by \right\|_{\max} \\
&\quad \leq \left\| \left( n^{-1}\bX'\bX - \E n^{-1}\bX'\bX \right) \bbeta \right\|_{\max} 
+\left\| n^{-1}\bX'\bep - \E n^{-1}\bX'\bep \right\|_{\max} \\
&\quad+ \left\| \left(n^{-1}\wtilde{\bX}(\btheta^0)'\bX - \E [n^{-1}\wtilde{\bX}(\btheta^0)'\bX] \right) \bbeta \right\|_{\max}
+\left\| n^{-1}\wtilde{\bX}(\btheta^0)'\bep - \E [n^{-1}\wtilde{\bX}(\btheta^0)'\bep] \right\|_{\max}  \\
&\quad =: Z_1 + Z_2 + Z_3 + Z_4.
\end{align*}
These terms can be bounded by the results obtained in the proof of $(a\textendash ii)$. 
We see that
\begin{align*}
Z_1 \leq s^{1/2}\left\| n^{-1}\bX'\bX_{\cS^0} - \E n^{-1}\bX'\bX_{\cS^0}\right\|_{\max} \left\| \bbeta_{\cS^0} \right\|_{2}
\lesssim s W_1 \lesssim s(n^{-1}\log p)^{1/2}
\end{align*}
with probability at least  $1- O(p^{-\nu})$.
Next we deduce
\begin{align*}
Z_2 \leq \left\| n^{-1}\bLambda^0{\bF^0}'\bep \right\|_{\max}+\left\| n^{-1}\bE'\bep \right\|_{\max}.
\end{align*}
The first and second terms can be bounded by the same ways as $W_{1,3}$ and $W_{3,1}$ in the proof of $(a)$ above with $\bE$ and $\bE_{\bbbeta^0}$ replaced by $\bep$, respectively. Then the first term dominates the second and hence 
$Z_2 \lesssim (n^{-1}\log p)^{1/2}$ with probability at least  $1- O(p^{-\nu})$. Similarly, we can obtain 
\begin{align*}
Z_3 \leq s^{1/2}\left\| n^{-1}\wtilde{\bX}(\btheta^0)'\bX_{\cS^0} - \E n^{-1}\wtilde{\bX}(\btheta^0)'\bX_{\cS^0} \right\|_{\max} \left\| \bbeta_{\cS^0} \right\|_{2}
\lesssim s W_3\lesssim s(n^{-1}\log p)^{1/2}
\end{align*}
with probability at least  $1- O(p^{-\nu})$.
Note that $Z_4$ has the same bound as $Z_2$. Consequently, collecting terms leads to the result, $Z_1+\dots+Z_4 \lesssim s(n^{-1}\log p)^{1/2}$ with probability at least  $1- O(p^{-\nu})$. This proves ($b\textendash ii$) and concludes the proof of Lemma \ref{lem:diff_U-EU}.

\subsection{Lemma \ref{lem:lasso} and its proof} \label{secB.3}
\begin{lem}\label{lem:lasso}
	Assume that all the conditions of Theorem \ref{thm:power} hold. Then with probability at least $1-O(\pi_{np})$, the Lasso solution in \eqref{auglasso} satisfies
	\begin{align*}
	\sup_{\btheta\in \Theta_{np}}\left\|\hat{\bbeta}^\textsf{aug}(\btheta)-\bbeta^\textsf{aug}\right\|_2 &= O(s^{1/2}\lambda), \\
	\sup_{\btheta\in \Theta_{np}}\left\|\hat{\bbeta}^\textsf{aug}(\btheta)-\bbeta^\textsf{aug}\right\|_1 &= O(s\lambda),
	\end{align*}
	where $\lambda =c_1n^{1/2}\log p$ with $c_1$ some positive constant. 
\end{lem}

\noindent\textit{Proof}. Let $\bdelta(:=\bdelta(\btheta)):=\hat{\bbeta}^\textsf{aug}(\btheta)-\bbeta^\textsf{aug}$. 
We start with introducing two inequalities
\begin{align}
&\sup_{\btheta\in\Theta_{np}}\left\|n^{-1}[\bX, \widetilde{\bX}(\btheta)]'\bep \right\|_{\max} \leq 2^{-1}\lambda, \label{ineq:event1} \\
&\inf_{\btheta\in\Theta_{np},\,\bdelta \in \mathbb{V}} \bdelta'\bU(\btheta)\bdelta /\|\bdelta\|_2^2 \geq \sigma_e^2(1+o(1)), \label{ineq:event2}
\end{align}
where $\lambda =c_1n^{-1/2}\log p$ for some positive constant $c_1$ and 
\begin{align}
\mathbb{V} = \left\{ \bdelta \in \mathbb{R}^{2p}: \|\bdelta_{\cS^1}\|_1 \leq 3\|\bdelta_{\cS^0}\|_1,\, \|\bdelta\|_0\leq k \right\}. \label{ineq:sset}
\end{align}
It is well known that the rate of convergence of the Lasso estimator can be obtained provided that \eqref{ineq:event1} and \eqref{ineq:event2} hold. Thus we show that these two inequalities actually hold with high probability in Step 1, and then derive the convergence rate using \eqref{ineq:event1} and \eqref{ineq:event2} in Step 2. 

\textbf{Step 1.} We check whether \eqref{ineq:event1} and \eqref{ineq:event2} actually hold with high probability. We first verify \eqref{ineq:event1}. 
By the proofs of Lemmas \ref{lem:diffsUV_max} and \ref{lem:diff_U-EU}, we have
\begin{align*}
&\sup_{\btheta\in \Theta_{np}}\left\|n^{-1}[\bX, \widetilde{\bX}(\btheta)]'\bep \right\|_{\max} \\
&\qquad\leq \left\| n^{-1}\bX'\bep \right\|_{\max} 
+ \sup_{\btheta\in \Theta_{np}}\left\| n^{-1}\widetilde{\bX}(\btheta)'\bep-n^{-1}\widetilde{\bX}(\btheta^0)'\bep \right\|_{\max} 
+ \left\|n^{-1}\widetilde{\bX}(\btheta^0)'\bep \right\|_{\max}.
\end{align*}
The first and third terms can both be upper bounded by $O(n^{-1/2}\log p)$ with probability at least $1-O(p^{-\nu})$, following the same lines for deriving bound for $Z_2$ in the proof of Lemma \ref{lem:diff_U-EU}. 
To evaluate the second term, we can use the argument about $V_2$ and its upper bound \eqref{ineq:V2} in the proof of Lemma \ref{lem:diffsUV_max}. That bound still holds with the same rate $O(n^{-1/2}\log p)$ even if we take $\cA=\{1,\dots,p\}$. Thus we conclude  that \eqref{ineq:event1} is true for the given $\lambda$ by choosing an appropriate positive large constant $c_1$, with probability  at least $1-O(\pi_{np})$.

Next to verify \eqref{ineq:event2}, we derive the population lower bound first and then show that the difference is negligible.  
From the construction, we have
\begin{align*}
&\E[n^{-1}{\wtilde{\bX}(\btheta^0)}'\wtilde{\bX}(\btheta^0)] = \E[n^{-1}{\bX}'\bX] = \bLambda^0 \bSigma_f{\bLambda^0}' + \sigma^2_e \bI_{p}, \\
&\E [n^{-1}{\wtilde{\bX}(\btheta^0)}'\bX] = \E [n^{-1}{\bX}'\wtilde{\bX}(\btheta^0)] = \bLambda^0 \bSigma_f{\bLambda^0}'.
\end{align*}
Using these equations, we obtain the lower bound
\begin{align}
\inf_{\bdelta \in \mathbb{V}} \bdelta' \E\left[ \bU(\btheta^0)\right] \bdelta/\|\bdelta\|_2^2 
&= \inf_{\bdelta \in \mathbb{V}} \bdelta' 
\begin{pmatrix}
\bLambda^0 \bSigma_f{\bLambda^0}'+\sigma^2_e \bI_{p} & \bLambda^0 \bSigma_f{\bLambda^0}' \\
\bLambda^0 \bSigma_f{\bLambda^0}' & \bLambda^0 \bSigma_f{\bLambda^0}'+\sigma^2_e \bI_{p}
\end{pmatrix}
\bdelta/\|\bdelta\|_2^2 \notag\\
&= 
\inf_{\bdelta \in \mathbb{V}} \bdelta' \left\{
\begin{pmatrix}
1 & 1 \\
1 & 1
\end{pmatrix}
\otimes \bLambda^0 \bSigma_f{\bLambda^0}' + \sigma_e^2 \bI_{2p} \right\}
\bdelta/\|\bdelta\|_2^2 \notag\\
&\geq \sigma_e^2. \label{ineq:popu}
\end{align}
Because $\bdelta\in\mathbb{V}$ is sparse and satisfies $|\cB|\leq k$ for $\cB:=\supp(\bdelta)$, 
it holds that $\bdelta' \bU(\btheta^0) \bdelta = \bdelta_{\cB}' \bU_{\cB}(\btheta^0)\bdelta_{\cB}$ and $\bdelta' \E\left[ \bU(\btheta^0)\right] \bdelta = \bdelta_{\cB}' \E\left[ \bU_{\cB}(\btheta^0)\right] \bdelta_{\cB}$. Hence from Lemma \ref{lem:diff_U-EU} together with the  condition on dimensionality, we obtain 
\begin{align}
\sup_{|{\cB}|\leq k,\, \btheta\in \Theta_{np}}\left\| \bU_{\cB}(\btheta) - \E [\bU_{\cB}(\btheta^0)] \right\|_{\max}  
&=O\left(k^{1/2}\tilde{c}_{np}\right) \notag\\
&=o(s^{-1})\label{ineq:diffrate}
\end{align}
with probability at least $1-O(\pi_{np})$. 
Thus using \eqref{ineq:diffrate},  we have for any $\bdelta\in \mathbb{V}$,
\begin{align*}
&\bdelta' \E[\bU(\btheta^0)]\bdelta - \bdelta' \bU(\btheta)\bdelta =\bdelta_{\cB}' \left\{ \E[\bU_{\cB}(\btheta^0)] - \bU_{\cB}(\btheta)\right\}\bdelta_{\cB} \\
&\qquad \leq \|\bdelta\|_1^2 \sup_{|{\cB}|\leq k,\, \btheta\in \Theta_{np}}\left\|\bU_{\cB}(\btheta) - \E [\bU_{\cB}(\btheta^0)] \right\|_{\max} 
= \left(\|\bdelta_{\cS^0}\|_1 + \|\bdelta_{\cS^1}\|_1 \right)^2 o(s^{-1}) \\
&\qquad \lesssim \|\bdelta_{\cS^0}\|_1^2 o(s^{-1}) 
\leq \|\bdelta_{\cS^0}\|_2^2 o(1)
\leq \|\bdelta\|_2^2 o(1).
\end{align*}
Rearranging the terms with \eqref{ineq:popu} yields
\begin{align*}
\inf_{\btheta\in\Theta_{np},\,\bdelta \in \mathbb{V}}\bdelta' \bU(\btheta) \bdelta/\|\bdelta\|_2^2 
&\geq \inf_{\bdelta \in \mathbb{V}}\bdelta' \E[\bU(\btheta^0)]\bdelta/\|\bdelta\|_2^2 
- |o(1)|
\geq \sigma_e^2 - |o(1)|,
\end{align*}
resulting in \eqref{ineq:event2}. 
In consequence, two inequalities \eqref{ineq:event1} and \eqref{ineq:event2} hold with probability at least $1-O(\pi_{np})$. 

\textbf{Step 2.} This part is well known in the literature (e.g., \cite{Negahban2012}) so we briefly give the proof omitting the details. 
Because the objective function is given by
\begin{align*}
\hat{\bbeta}^\textsf{aug} (\btheta)
&= \arg\min_{\bb\in\mathbb{R}^{2p}} n^{-1}\left\| \by - [\bX, \widetilde{\bX}(\btheta)] \bb \right\|_2^2 + \lambda \|\bb\|_1, 
\end{align*}
the global optimality of the Lasso estimator implies
\begin{align*}
&(2n)^{-1}\left\|\by - [\bX, \widetilde{\bX}(\btheta)]\hat{\bbeta}^\textsf{aug}(\btheta) \right\|_2^2 + \lambda \left\|\hat{\bbeta}^\textsf{aug}(\btheta) \right\|_1 \\
&\qquad\qquad\leq (2n)^{-1}\left\|\by - [\bX, \widetilde{\bX}(\btheta)]{\bbeta}^\textsf{aug} \right\|_2^2 + \lambda \left\|{\bbeta}^\textsf{aug} \right\|_1,
\end{align*}
where the true parameter vector $\bbeta^\textsf{aug}$ was defined in the proof of Theorem \ref{thm:power}. 
Note that $\sup_{\btheta \in \Theta_{np}}\|\bdelta(\btheta)\|_0\leq k$ by the assumption. 
Expanding the inequality and collecting terms with \eqref{ineq:event1} yield
\begin{align}
2^{-1}\bdelta' \bU(\btheta) \bdelta 
\leq \left\|n^{-1}\bep'[\bX, \widetilde{\bX}(\btheta)] \right\|_{\max} \|\bdelta\|_1 +\lambda \|\bdelta\|_1 
\leq (3/2)\lambda \|\bdelta\|_1. \label{ineq:optlasso}
\end{align}
On the other hand, applying Lemma 1 of \cite{Negahban2012} to our model reveals that $\bdelta \in \mathbb{V}$.  
Thus  we can use \eqref{ineq:event2}, \eqref{ineq:optlasso}, and \eqref{ineq:sset} to get  
\begin{align*}
\|\bdelta\|_2^2(\sigma_e^2+o(1))
\leq 3\lambda \|\bdelta\|_1 = 3\lambda \left( \|\bdelta_{\cS^1}\|_1 + \|\bdelta_{\cS^0}\|_1 \right) \leq 12 \lambda\|\bdelta_{\cS^0}\|_1. 
\end{align*}
Since $|\cS^0|=s$ and $\|\bdelta_{\cS^0}\|_1\leq s^{1/2} \|\bdelta_{\cS^0}\|_2$, it holds that 
$\|\bdelta\|_2\leq 12 s^{1/2}\lambda/(\sigma_e^2+o(1))$. 
Since $\|\bdelta_{\cS^0}\|_2 \leq \|\bdelta\|_2$, we obtain the desired bound $\|\bdelta\|_1 \leq 48 s\lambda/(\sigma_e^2+o(1))$. 
This bound holds uniformly over $\btheta\in\Theta_{np}$, which completes the proof of Lemma \ref{lem:lasso}.

\section{Additional technical lemmas and their proofs} \label{sec:additionallem}

\subsection{Lemma \ref{lem:HVbound} and its proof} \label{secC.1}
\begin{lem}\label{lem:HVbound}
	Denote by $\bV\in\mathbb{R}^{r\times r}$ a diagonal matrix with its entries the $r$ largest eigenvalues of $(np)^{-1}\bX\bX'$ and define $\bH= ({\bLambda^0}'\bLambda^0/p)({\bF^0}'\hat{\bF}/n) \bV^{-1}$. 
	Assume that Conditions \ref{ass:fac}--\ref{ass:eigensepa} hold. Then $\|\bH\|_2+\|\bH^{-1}\|_2+\|\bV\|_2+\|\bV^{-1}\|_2$ is bounded from above by some constant with probability at least $1-O(p^{-\nu})$. 
\end{lem}

\noindent\textit{Proof}. 
Let $\lambda^{k}[\bA]$ denote the $k$th largest eigenvalue of square matrix $\bA$ throughout the proof. 
Because $\|{\bLambda^0}'\bLambda^0/p\|_2\leq M$ and 
\begin{align*}
\|{\bF^0}'\hat{\bF}/n\|_2 &\leq \|n^{-1/2}\bF^0\|_2\|n^{-1/2}\hat{\bF}\|_2 \\
&\leq (rn)^{1/2}\|n^{-1/2}\bF^0\|_{\max} \left(\lambda^1[n^{-1}\hat{\bF}'\hat{\bF}]\right)^{1/2}
\leq r^{1/2}M
\end{align*}
by Conditions \ref{ass:fac}--\ref{ass:floa}, and $\hat{\bF}'\hat{\bF}/n=\bI_r$, we have 
\begin{align*}
\|\bH\|_2 
\leq \left\|{\bLambda^0}'\bLambda^0/p\right\|_2 \left\|{\bF^0}'\hat{\bF}/n\right\|_2 \left\|\bV^{-1}\right\|_2
\lesssim \left\|\bV^{-1}\right\|_2,
\end{align*}
where $\|\bV^{-1}\|_2$ is equal to the reciprocal of the $r$th largest eigenvalue of $(np)^{-1}\bX\bX'$. 
Similarly, under Conditions \ref{ass:fac}--\ref{ass:floa} we also have 
\begin{align*}
\left\|\bH^{-1}\right\|_2\leq \|\bV\|_2\left\|({\bF^0}'\hat{\bF}/n)^{-1}\right\|_2\left\|({\bLambda^0}'\bLambda^0/p)^{-1}\right\|_2
\lesssim \left\|\bV\right\|_2\left\|({\bF^0}'\hat{\bF}/n)^{-1}\right\|_2,
\end{align*}
where $\|\bV\|_2$ is equal to the largest eigenvalue of $(np)^{-1}\bX\bX'$ and the inverse matrix in the upper bound is well defined by \cite{bai2003}. 
To see if $\|({\bF^0}'\hat{\bF}/n)^{-1}\|_2$ is bounded from above, it suffices to bound the minimum eigenvalue of ${\bF^0}'\hat{\bF}\hat{\bF}'{\bF^0}/n^2$ away from zero uniformly in $n$. 
Regarding $r$ eigenvalues of the matrix, Sylvester's law of inertia (e.g., \cite{HornJohnson2012}, Theorem 4.5.8) entails that all the $r$ eigenvalues are positive for all $n$. 
Moreover, by Proposition 1 of \cite{bai2003} we know that the limiting matrix of $\hat{\bF}'{\bF^0}/n$ is nonsingular under Conditions \ref{ass:fac} and \ref{ass:eigensepa}.
Therefore, we can conclude that $\liminf_{n\to\infty}\lambda^r[{\bF^0}'\hat{\bF}\hat{\bF}'{\bF^0}/n^2]>0$ a.s., and hence $\|\bH^{-1}\|_2\lesssim \|\bV\|_2$ follows. 

To complete the proof, it is sufficient to show that the maximum and $r$th largest eigenvalues of $(np)^{-1}\bX\bX'$ are bounded from 
above and away from zero, respectively, for all large $n$ and $p$. 
By the definition of the spectral norm and triangle inequality, we have
\begin{align*}
\left\{\lambda^{1}\left[(np)^{-1}\bX\bX') \right]\right\}^{1/2} 
&=\left\|(np)^{-1/2}\bX\right\|_2 
\leq \left\|(np)^{-1/2} \bF^0{\bLambda^0}'\right\|_2 + \left\|(np)^{-1/2} \bE\right\|_2 \\
&\leq \left\|n^{-1/2} \bF^0\right\|_2 \left\|p^{-1/2}\bLambda^0\right\|_2 + \left\|(np)^{-1/2} \bE\right\|_2. 
\end{align*}
By Conditions \ref{ass:fac} and \ref{ass:floa}, the first term is a.s.\ bounded by a constant as discussed above. The second term is $O((n\wedge p)^{-1/2})=o(1)$ with probability at least $1-2\exp(-|O(n\vee p)|)$ by Lemma \ref{lem:maxeigv}(a) under Condition \ref{ass:err_new}. 
Therefore, the largest eigenvalue of $(np)^{-1}\bX\bX'$ is bounded from above by some constant with probability at least $1-2\exp(-|O(n\vee p)|)$.

Next we bound the $r$th largest eigenvalue of $(np)^{-1}\bX\bX'$ away from zero. Since the matrix is symmetric, Weyl's inequality (e.g., \cite{HornJohnson2012}, Theorem 4.3.1) yields 
\begin{align}
&\lambda^{r}\left[ (np)^{-1}\bX\bX' \right] = \lambda^{r}\left[ (np)^{-1} \left\{ \bF^0{\bLambda^0}'\bLambda^0 {\bF^0}' + \left(\bE\bLambda^0 {\bF^0}' + \bF^0{\bLambda^0}'\bE' \right) + \bE\bE' \right\} \right] \notag\\
&\quad \geq \lambda^{r}\left[ (np)^{-1} \bF^0{\bLambda^0}'\bLambda^0 {\bF^0}'\right] + \lambda^{n}\left[ (np)^{-1} \left(\bE\bLambda^0 {\bF^0}' + \bF^0{\bLambda^0}'\bE' \right) \right] + \lambda^{n}\left[ (np)^{-1} \bE\bE' \right]. \label{ineq:lwbd}
\end{align}
The third term of lower bound \eqref{ineq:lwbd} is obviously nonnegative. 
For the first term of lower bound \eqref{ineq:lwbd}, let $\mathcal{V}$ denote a subspace of $\mathbb{R}^n$. Because $\bF^0{\bLambda^0}'\bLambda^0 {\bF^0}'$ is symmetric, the Courant--Fischer min-max Theorem (e.g., \cite{HornJohnson2012}, Theorem 4.2.6) yields 
\begin{align*}
&\lambda^{r}\left[ (np)^{-1} \bF^0{\bLambda^0}'\bLambda^0 {\bF^0}'\right] 
= \max_{\mathcal{V}:\dim(\mathcal{V})=r} \min_{\bv\in \mathcal{V}\backslash\{\bzero\}} \left\{ (np)^{-1} \frac{\bv' \bF^0{\bLambda^0}'\bLambda^0 {\bF^0}' \bv}{\bv'\bv} \right\} \\
&\qquad \geq \max_{\mathcal{V}:\dim(\mathcal{V})=r} \min_{\bv\in \mathcal{V}\backslash\{\bzero\}} \left( n^{-1} \frac{\bv'\bF^0{\bF^0}'\bv}{\bv'\bv} \right) \min_{ {\bF^0}'\bv \in\mathbb{R}^r\backslash\{\bzero\}} \left( p^{-1}\frac{\bv' \bF^0{\bLambda^0}'\bLambda^0{\bF^0}' \bv}{\bv'\bF^0{\bF^0}'\bv} \right) \\
&\qquad = \lambda^r \left[ n^{-1}\bF^0{\bF^0}' \right] \lambda^r \left[ p^{-1} {\bLambda^0}'\bLambda^0 \right]
= \lambda^r \left[ n^{-1}{\bF^0}'{\bF^0} \right] \lambda^r \left[ p^{-1} {\bLambda^0}'\bLambda^0 \right] \\
&\qquad \geq \lambda^r \left[\bSigma_f \right] \lambda^r \left[ p^{-1} {\bLambda^0}'\bLambda^0 \right] - \left\|n^{-1}{\bF^0}'\bF^0 -\bSigma_f \right\|_2 \\
&\qquad \geq \lambda^r \left[\bSigma_f \right] \lambda^r \left[ p^{-1} {\bLambda^0}'\bLambda^0 \right] - r\left\|n^{-1}{\bF^0}'\bF^0 -\bSigma_f \right\|_{\max}.
\end{align*}
In this lower bound, the first term is bounded away from zero by Conditions \ref{ass:fac}--\ref{ass:floa}. Meanwhile, to evaluate the second term we use Lemma \ref{lem:ineq}(d) in Section \ref{secC.2}, which together with the union bound establishes
\begin{align*}
&\Pro \left(\left\|n^{-1}{\bF^0}'\bF^0 -\bSigma_f \right\|_{\max} > u \right) 
\leq r^2 \max_{k,\ell\in\{1,\dots,r\}} \Pro  \left(\left|n^{-1}\sum_{i=1}^n\left(f_{ik}^0f_{i\ell}^0 -\E f_{ik}^0f_{i\ell}^0 \right) \right| > u \right) \\
&\qquad \leq 2r^2 \exp(-nu^2/C) 
\end{align*}
for any $0\leq u\leq c$. Thus the second one turns out to be $O((n^{-1}\log p)^{1/2})=o(1)$ with probability at least $1-O(p^{-\nu})$ once we set $u=(C\nu n^{-1}\log p)^{1/2}$ and assume $n^{-1}\log p\leq c^2/(C\nu)$ without loss of generality. Therefore, the first term of lower bound \eqref{ineq:lwbd} is bounded away from zero eventually. 
For the second term of \eqref{ineq:lwbd}, since the spectral norm gives the upper bound of the spectral radius we have
\begin{align*}
&\left|\lambda^{n}\left[ (np)^{-1} \left(\bE\bLambda^0 {\bF^0}' + \bF^0{\bLambda^0}'\bE' \right) \right] \right| 
\leq \left\| (np)^{-1} \left(\bE\bLambda^0 {\bF^0}' + \bF^0{\bLambda^0}'\bE' \right) \right\|_2\\
&\qquad \leq 2\left\| (np)^{-1/2} \bE \right\|_2 \left\|p^{-1/2}\bLambda^0\right\|_2 \left\|n^{-1/2}\bF^0\right\|_2 \\
&\qquad =O\left((n\wedge p)^{-1/2}\right)O(1)O(1)=o(1),
\end{align*}
which holds with probability at least $1-2\exp(-|O(n\vee p)|)$ by Lemma \ref{lem:maxeigv}(a) in Section \ref{secC.4}. As a consequence, the desired result holds with probability at least $1-O(p^{-\nu})$ and this concludes the proof of Lemma \ref{lem:HVbound}.

\subsection{Lemma \ref{lem:ineq}  and its proof} \label{secC.2}
\begin{lem}\label{lem:ineq} 
	Assume that Conditions \ref{ass:fac}--\ref{ass:err_new} hold. Then there exist some positive constants $c$ and $C$ such that the following inequalities hold 
	\begin{enumerate}
		\item [(a)] For all $\ell,i \in \{1,\dots,n\}$ and $0\leq u\leq c$, we have
		\begin{align*}
		\Pro \left( \left|p^{-1}\sum_{j=1}^p\left(e_{\ell j}e_{ij}-\E[e_{\ell j}e_{ij}] \right) \right| > u \right)
		\leq 2 \exp\left( - pu^2 /C \right).
		\end{align*}
		\item [(b)] For all $k\in\{1,\dots,r\}$, $j \in \{1,\dots,p\}$, and $0\leq u\leq c$, we have 
		\begin{align*}
		\Pro \left( \left|n^{-1}\sum_{i=1}^n f_{ik}^0e_{ij}\right|  > u \right)
		\leq 2\exp\left( - nu^2/C \right).
		\end{align*}
		\item [(c)] For all $k\in\{1,\dots,r\}$, $i \in \{1,\dots,n\}$, and $u\geq 0$, we have
		\begin{align*}
		\Pro \left( \left|p^{-1}\sum_{j=1}^p \lambda_{jk}^0 e_{ij} \right| > u \right)
		\leq 2 \exp \left(-pu^2/C \right).
		\end{align*}
		\item [(d)] For all $k,\ell \in \{1,\dots,r\}$ and $0\leq u\leq c$, we have
		\begin{align*}
		\Pro \left( \left|n^{-1}\sum_{i=1}^n\left(f_{i k}^0f_{i\ell}^0-\E[f_{i k}^0f_{i\ell}^0] \right) \right| > u \right)
		\leq 2 \exp\left( - nu^2/C \right).
		\end{align*}
	\end{enumerate}
\end{lem}

\noindent\textit{Proof}. 
(a) To obtain the first result, we rely on the Hanson--Wright inequality. Let $\bxi=(\xi_1,\dots,\xi_m)'\in \mathbb{R}^m$ denote a random vector whose components are independent copies of $e\sim \mbox{subG}(C_e^2)$. Then the inequality states that for any (nonrandom) matrix $\bA \in \mathbb{R}^{m\times m}$, 
\begin{align}
\Pro \left( \left| \bxi'\bA\bxi - \E \bxi'\bA\bxi  \right| > u \right) 
\leq 2 \exp\left\{ -\wtilde{C}_H \min \left( \frac{u^2}{K^4\|\bA\|_F^2}, \frac{u}{K^2\|\bA\|_2} \right) \right\}, \label{ineq:HW}
\end{align}
where $K$ is a positive constant such that $\sup_{k\geq 1}k^{-1/2}(\E|e|^k)^{1/k}\leq K$ and $\wtilde{C}_H$ is a positive constant. 
In our setting, we can take $K=3C_e^2$ (e.g., Lemma 1.4 of \cite{RigolletHutter2017}). 
Using this inequality, we first prove the case when $\ell=i$. If we set $m=p$ and $\bA=\diag(p^{-1},\dots,p^{-1})$, then we have 
\begin{align*}
\left|\bxi'\bA\bxi - \E \bxi'\bA\bxi \right| 
= \left| p^{-1}\sum_{j=1}^p (\xi_j^2- \E \xi_j^2) \right|
\overset{d}{=} \left|p^{-1}\sum_{j=1}^p\left(e_{ij}^2 - \E[e_{ij}^2] \right) \right|
\end{align*}
for all $i$. Moreover, we obtain $\|\bA\|_F^2=p^{-1}$ and $\|\bA\|_2=p^{-1}$ in this case. The assumed condition $0<u\leq 9C_e^2=K^2$ entails that  $u^2/K^4\leq u/K^2$ so the result follows from \eqref{ineq:HW} with $\wtilde{C}_H$ replaced by  $C_H=81C_e^4/\wtilde{C}_H$. 

Similarly, we prove the case when $\ell\not=i$. 
We set $m=p+1$ and $\bA=(\ba_1,\dots,\ba_{p+1})$, where $\ba_1=(0, p^{-1},\dots,p^{-1})'$ and 
$\ba_j=\bzero$ for $j=2,\dots,p+1$. That is, the entries of $\bA$ are all zero except that the second to $(p+1)$th components in the first column vector are $p^{-1}$. Under this setting, we observe that
\begin{align*}
\left|\bxi'\bA\bxi - \E \bxi'\bA\bxi \right| 
= \left| p^{-1}\sum_{j=2}^{p+1} \xi_1\xi_j \right|
\overset{d}{=} \left|p^{-1}\sum_{j=1}^p e_{\ell j}e_{ij}  \right|
\end{align*}
for all $\ell\not=i$. Moreover, we obtain $\|\bA\|_F^2=\|\bA\|_2=p^{-1}$ in this case. Therefore, the same bound holds as in the case of $\ell=i$ from \eqref{ineq:HW} again. Consequently, for any $0\leq u\leq 9C_e^2$ we have
\begin{align*}
\Pro \left( \left|p^{-1}\sum_{j=1}^p\left(e_{\ell j}e_{ij}-\E[e_{\ell j}e_{ij}] \right) \right| > u \right)
\leq 2 \exp\left( - pu^2 /C_H\right).
\end{align*}

(b) We prove the second assertion by Bernstein's inequality for the sum of a martingale difference sequence (e.g., Theorem 3.14 in \cite{Bercu2015}). Fix $k=1$ and $j=1$. Define $\mathcal{F}_{i-1}$ as the $\sigma$-field generated from $\{f_{\ell 1}^0: \ell=i,i-1,\dots \}$. 
Then $(f_{i1}^0e_{i1},\mathcal{F}_{i})$ forms a martingale difference sequence because $\E|f_{i1}^0e_{i1}|<\infty$ and $\E[f_{i1}^0e_{i1}|\mathcal{F}_{i-1}]=0$ under Conditions \ref{ass:fac} and \ref{ass:err_new}. 
Since the sub-Gaussianity of $e_{i1}$ implies $\E e_{i1}^2\leq 4C_e^2$ (e.g., Lemma 1.4 of \cite{RigolletHutter2017}), we have 
$V_i := \E \left[ f_{ik}^{0\,2}e_{ij}^2 \mid \mathcal{F}_{i-1} \right] \leq 4 C_e^2M^2$, and hence  
$\sum_{i=1}^n V_i \leq 4n C_e^2 M^2$ a.s.\ due to boundedness $|f_{i1}^0|\leq M$ a.s. On the other hand, by the sub-Gaussianity of $e_{ij}$ and boundedness of $|f_{i1}^0|$ again we observe that for all $p\geq 3$ and $i\in\{1,\dots,n\}$,
\begin{align*}
\E \left[ (0\vee f_{i1}^0e_{i1})^p \mid \mathcal{F}_{i-1} \right] 
\leq M^{p} (2C_e^2)^{p/2} p \Gamma (p/2) 
\leq p! (2C_eM)^{p-2} V_i /2,
\end{align*}
where $\Gamma$ denotes the Gamma function and we have used the estimates $p\Gamma(p/2)\leq p!$ and $2^{p/2-2}\leq 2^{p-2}/2$ for $p\geq 3$ in the last inequality. 
Then an application of Theorem 3.14 in \cite{Bercu2015} by putting $x=u$, $y=4M^2C_e^2$, and $c=2MC_e$ in their notation gives the one-sided result. Making twice the bound yields
\begin{align*}
\Pro \left( \left|n^{-1}\sum_{i=1}^n f_{ik}^0e_{ij}\right|  > u \right)
\leq 2\exp\left( - \frac{nu^2}{8M^2C_e^2+4MC_eu} \right).
\end{align*}
For all $0\leq u\leq MC_e^2$, the upper bound is further bounded by $2\exp(nu^2/(12M^2C_e^2))$. We set $C_I=12M^2C_e^2$. Consequently, for any $0\leq u\leq MC_e^2$ we have 
\begin{align*}
\Pro \left( \left|n^{-1}\sum_{i=1}^n f_{ik}^0e_{ij}\right|  > u \right)
\leq 2\exp\left( - nu^2/C_I \right).
\end{align*}

(c) We prove the third inequality. Note that 
\begin{align*}
\Pro \left( \left|\lambda_{jk}^0 e_{ij}\right| > u \right) 
\leq 2 \exp\left\{-\frac{u^2}{2\lambda_{jk}^{02}C_e^2}\right\} 
\leq 2 \exp\left\{-\frac{u^2}{2M^2C_e^2}\right\}.
\end{align*}
This implies that $\lambda_{jk}^0 e_{ij}$ is a sequence of i.i.d.\ subG($M^2C_e^2$). 
Thus the result is obtained directly by Bernstein's inequality for the sum of independent sub-Gaussian random variables. Consequently, for any $u\geq 0$ putting $C_J=M^2C_e^2$ leads to 
\begin{align*}
\Pro \left( \left|p^{-1}\sum_{j=1}^p \lambda_{jk}^0 e_{ij} \right| > u \right)
\leq 2 \exp \left(-pu^2/C_J \right).
\end{align*}

(d) We show the last inequality. Note that for each $k$, $(f_{ik})_i \sim \mbox{i.i.d.\ subG}(M^2)$ since $|f_{ik}^0|\leq M$ a.s.\ by Lemma 1.8 of \cite{RigolletHutter2017} under Condition \ref{ass:fac}. Thus the remaining is the same as (a). Set $C_K=81M^4/\wtilde{C}_H$ here. Then for any $0\leq u\leq 9M^2$, we have
\begin{align*}
\Pro \left( \left|n^{-1}\sum_{i=1}^n\left(f_{i k}^0f_{i\ell}^0-\E[f_{i k}^0f_{i\ell}^0] \right) \right| > u \right)
\leq 2 \exp\left( - nu^2/C_K \right).
\end{align*}

Finally the obtained inequalities hold even if the constant in the upper bound is replaced with arbitrary fixed constant $C$ such that $C\geq \max\{C_H,C_I,C_J,C_K\}$. Similarly, we can also restrict the range of $u$ for each inequality to be $0\leq u \leq c$ for arbitrary fixed constant $c$ that satisfies $0< c\leq \min(9C_e^2,MC_e^2,9M^2)$. This completes the proof of Lemma \ref{lem:ineq}.

\subsection{Lemma \ref{lem:diffsUV_max} and its proof} \label{secC.3}
\begin{lem}\label{lem:diffsUV_max}
	Assume that Conditions \ref{ass:regerr}--\ref{ass:err_new} hold. Then for any set $\cA$ satisfying $|\cA|\leq k$, the following statements hold with probability at least $1-O(\pi_{np})$
	\begin{align*}
	(a)&~~\sup_{\btheta \in \Theta_{np}}\left\|\bU_{\cA}({\btheta}) - \bU_{\cA}(\btheta^0) \right\|_{\max} = O\left(k^{1/2}\tilde{c}_{np}\right),\\
	(b)&~~\sup_{\btheta \in \Theta_{np}}\left\|\bv_{\cA}({\btheta}) - \bv_{\cA}(\btheta^0) \right\|_{\max} = O\left(s^{3/2}\tilde{c}_{np}\right),
	\end{align*}
	where $\Theta_{np}$ was defined in Lemma \ref{lem:rateconv} and $\tilde{c}_{np} = n^{-1/2}\log p + p^{-1/2}\log n$. Consequently, we have
	\begin{align*}
	\sup_{\btheta \in \Theta_{np}}\left\|\bT_{\cA}({\btheta}) - \bT_{\cA}(\btheta^0) \right\|_{\max} 
	= O\left( \left( k^{1/2} + s^{3/2}\right)\tilde{c}_{np}\right).
	\end{align*}
\end{lem}

\noindent\textit{Proof}. 
We first state some results that are useful in the proof. Since $\|n^{-1/2}\bF^{0}\|_2=O(1)$ a.s.\ by Condition \ref{ass:fac} and $\|k^{-1/2}\bLambda_\cA^{0} \|_2=O(1)$ for any $\cA$ such that $|\cA|\leq k$ under Condition \ref{ass:floa}, we first have
\begin{align*}
\left\| n^{-1/2}{\bC_{\cA}^0} \right\|_2
\leq \left\|n^{-1/2}\bF^{0}\right\|_2 k^{1/2}\left\|k^{-1/2}\bLambda_\cA^0 \right\|_2\lesssim k^{1/2}.
\end{align*}
Next Lemma \ref{lem:maxeigv}(b) in Section \ref{secC.4} gives directly 
\begin{align}\label{bound:E1}
\left\|n^{-1/2}\bE_{\bbbeta^0 \cA} \right\|_2 
\lesssim 1
\end{align}
with probability at least $1-O(p^{-\nu})$. {By Condition \ref{ass:err_new}, we also deduce
	\begin{align*}
	\Pro\left( \sup_{\bbbeta\in\bTheta_{np}}\left\|\bE_{\bbbeta}-\bE_{\bbbeta^0} \right\|_{\max} > u\right) 
	&\leq np \max_{i,j} \Pro \left( \sup_{\bbbeta\in\bTheta_{np}}\left|e_{\bbbeta ij}-e_{\bbbeta^0 ij} \right| > u\right) \\
	&\leq np \max_{i,j} \Pro \left(  |Z| > u/(M^{1/2}c_{np}^{1/2}) \right) \\
	&\leq 2 np\exp\left( -u^2/\left(c_e^2Mc_{np}\right)\right)
	\end{align*}
	for any $u\geq 0$. Thus setting $u=2c_eM^{1/2}c_{np}^{1/2}\log^{1/2}(np)$ with some large enough positive constant $M$,  we obtain that with probability at least $1-O((np)^{-\nu})$, 
	\begin{align*}
	\sup_{\bbbeta\in\bTheta_{np}}\left\|\bE_{\bbbeta}-\bE_{\bbbeta^0} \right\|_{\max} 
	\lesssim c_{np}\log^{1/2}(np) = O(\tilde{c}_{np}).
	\end{align*}}
We will use these results and Lemma \ref{lem:matnorm} in Section \ref{secC.5} in the proofs below. 

To prove (a), we have
\begin{align*}
&\left\|\bU_{\cA}(\btheta) - \bU_{\cA}(\btheta^0) \right\|_{\max} 
\leq \left\|n^{-1}\wtilde{\bX}_{\cA}(\btheta)'\wtilde{\bX}_{\cA}(\btheta) 
- n^{-1}\wtilde{\bX}_{\cA}(\btheta^0)'\wtilde{\bX}_{\cA}(\btheta^0) \right\|_{\max} \\
&\qquad + 2\left\|n^{-1}\bX_{\cA}'\wtilde{\bX}_{\cA}(\btheta) 
- n^{-1}\bX_{\cA}'\wtilde{\bX}_{\cA}(\btheta^0) \right\|_{\max}=:U_1+U_2.
\end{align*}
Observe that $U_1$ is further bounded as
\begin{align*}
&U_1 
\leq \left\|n^{-1}\bC_{\cA}'\bC_{\cA} - n^{-1}{\bC_{\cA}^0}'\bC_{\cA}^0 \right\|_{\max} 
+ \left\|n^{-1}\bE_{\bbbeta \cA}'\bE_{\bbbeta \cA} - n^{-1}\bE_{\bbbeta^0\cA}'\bE_{\bbbeta^0\cA} \right\|_{\max} \\
&\qquad + 2\left\|n^{-1}\bE_{\bbbeta \cA}'\bC_{\cA} -n^{-1}\bE_{\bbbeta^0 \cA}'\bC_{\cA}^0 \right\|_{\max}=:U_{11} + U_{12} + U_{13}.
\end{align*}
By Lemma \ref{lem:matnorm}, it is easy to see that 
\begin{align*}
U_{11}
&\leq \left\|n^{-1}\left(\bC_{\cA}-\bC_{\cA}^0 \right)' \left( \bC_{\cA}-\bC_{\cA}^0 \right) \right\|_{\max} + 2\left\| n^{-1}{\bC_{\cA}^0}' \left(\bC_{\cA}- {\bC_{\cA}^0}\right)\right\|_{\max} \\
&\leq n^{-1/2}\left\| \bC_\cA-\bC_\cA^0 \right\|_{\max} \left\|\bC_\cA-\bC_\cA^0 \right\|_2 + 2\left\|n^{-1/2}\bC_{\cA}^0 \right\|_2\left\|\bC_{\cA} - \bC_{\cA}^{0} \right\|_{\max}   \\
&\lesssim k^{1/2}\left\|{\bC}-\bC^{0} \right\|_{\max}^2 + k^{1/2}\left\|{\bC}-\bC^{0} \right\|_{\max} \\
&=O\left( k^{1/2}c_{np}^2 + k^{1/2}c_{np} \right) = O\left( k^{1/2}c_{np} \right), 
\end{align*}
where the last estimate follows from Lemma \ref{lem:rateconv}. 
Similarly, we deduce
\begin{align*}
U_{12}
&\leq \left\|n^{-1}\left(\bE_{\bbbeta \cA}-\bE_{\bbbeta^0 \cA} \right)' \left( \bE_{\bbbeta \cA}-\bE_{\bbbeta^0 \cA} \right) \right\|_{\max} + 2\left\| n^{-1}\bE_{\bbbeta^0 \cA}'\left(\bE_{\bbbeta \cA}- \bE_{\bbbeta^0 \cA}\right) \right\|_{\max} \\
&\leq n^{-1/2}\left\| \bE_{\bbbeta \cA}-\bE_{\bbbeta^0 \cA} \right\|_{\max} \left\|\bE_{\bbbeta \cA}-\bE_{\bbbeta^0 \cA} \right\|_2 + 2\left\|n^{-1/2}\bE_{\bbbeta^0 \cA} \right\|_2\left\|\bE_{\bbbeta \cA}- \bE_{\bbbeta^0 \cA} \right\|_{\max}   \\
&\lesssim k^{1/2}\left\|\bE_{\bbbeta}-\bE_{\bbbeta^0} \right\|_{\max}^2 
+  \left\|\bE_{\bbbeta}-\bE_{\bbbeta^0} \right\|_{\max} \\
&= O \left( k^{1/2}\tilde{c}_{np}^{2} + \tilde{c}_{np} \right) 
\end{align*}
and 
\begin{align*}
U_{13}
&\leq \left\|n^{-1}\left(\bE_{\bbbeta \cA}-\bE_{\bbbeta^0 \cA} \right)' \left( \bC_{\cA}-\bC_{\cA}^0 \right) \right\|_{\max} \\
&\qquad + \left\| n^{-1}\bE_{\bbbeta^0 \cA}'\left(\bC_{\cA}- \bC_{\cA}^0\right) \right\|_{\max} 
+ \left\| n^{-1}{\bC_{\cA}^0}'\left(\bE_{\bbbeta \cA}- \bE_{\bbbeta^0 \cA}^0\right) \right\|_{\max}\\
&\leq k^{1/2}\left\|\bE_{\bbbeta}-\bE_{\bbbeta^0} \right\|_{\max} \left\| \bC-\bC^0 \right\|_{\max} \\
&\qquad + \left\|n^{-1/2}\bE_{\bbbeta^0 \cA} \right\|_2 \left\| \bC- \bC^0 \right\|_{\max} 
+ \left\| n^{-1/2}{\bC_{\cA}^0} \right\|_2 \left\| \bE_{\bbbeta}- \bE_{\bbbeta^0}^0\right\|_{\max} \\
&= O \left(k^{1/2}\tilde{c}_{np}c_{np} + c_{np} + k^{1/2}\tilde{c}_{np} \right) 
= O \left(k^{1/2}\tilde{c}_{np}\right).
\end{align*}
Combining these bounds of $U_{11}$--$U_{13}$, we have
\begin{align*}
U_1\leq U_{11}+U_{12}+U_{13} \lesssim k^{1/2}\tilde{c}_{np}.
\end{align*}
This holds uniformly in $\btheta\in\Theta_{np}$ with probability at least $1-O(\pi_{np})$ by Lemma \ref{lem:rateconv} and the discussion above. 
Next we obtain 
\begin{align*}
U_2
&\leq \left\|n^{-1}{\bC_{\cA}^0}'(\bC_{\cA}-\bC_{\cA}^0) \right\|_{\max} 
+ \left\|n^{-1}{\bC_{\cA}^0}'(\bE_{\bbbeta \cA}-\bE_{\bbbeta^0 \cA}) \right\|_{\max} \\
&\qquad + \left\|n^{-1}\bE_{\cA}'(\bC_{\cA}-\bC_{\cA}^0) \right\|_{\max} 
+ \left\|n^{-1}\bE_{\cA}'(\bE_{\bbbeta \cA}-\bE_{\bbbeta^0 \cA}) \right\|_{\max} \\
&\leq \left\|n^{-1/2}{\bC_{\cA}^0}\right\|_2 \left\|\bC_{\cA}-\bC_{\cA}^0 \right\|_{\max} 
+ \left\|n^{-1/2}{\bC_{\cA}^0}\right\|_2 \left\|\bE_{\bbbeta \cA}-\bE_{\bbbeta^0 \cA} \right\|_{\max} \\
&\qquad + \left\|n^{-1/2}\bE_{\bbbeta^0 \cA} \right\|_2  \left\|\bC_{\cA}-\bC_{\cA}^0 \right\|_{\max} 
+ \left\|n^{-1/2}\bE_{\bbbeta^0 \cA} \right\|_2 \left\|\bE_{\bbbeta}-\bE_{\bbbeta^0} \right\|_{\max} \\
&=O\left( k^{1/2}c_{np} + k^{1/2}\tilde{c}_{np} + c_{np} + \tilde{c}_{np} \right) \\
&= O\left(k^{1/2}\tilde{c}_{np} \right). 
\end{align*}
This also holds uniformly in $\btheta\in\Theta_{np}$ with probability at least $1-O(\pi_{np})$ by Lemma \ref{lem:rateconv} and the discussion above. Consequently, it holds that 
\begin{align*}
\sup_{\btheta \in \Theta_{np}}\left\|\bU_{\cA}(\btheta) - \bU_{\cA}(\btheta^0) \right\|_{\max} 
\lesssim k^{1/2}\tilde{c}_{np}
\end{align*}
with probability at least $1-O(\pi_{np})$. 

To prove (b), we have
\begin{align*}
&\left\|\bv_{\cA}(\btheta) - \bv_{\cA}(\btheta^0) \right\|_{\max} 
\leq \left\|n^{-1}\wtilde{\bX}_{\cA}(\btheta)'\by 
- n^{-1}\wtilde{\bX}_{\cA}(\btheta^0)'\by \right\|_{\max} \\
&\quad\leq \left\|n^{-1}\wtilde{\bX}_{\cA}(\btheta)'\bX \bbeta
- n^{-1}\wtilde{\bX}_{\cA}(\btheta^0)'\bX \bbeta\right\|_{\max} 
+ \left\|n^{-1}\wtilde{\bX}_{\cA}(\btheta)'\bep
- n^{-1}\wtilde{\bX}_{\cA}(\btheta^0)'\bep \right\|_{\max} \\
&\quad =:V_1+V_2.
\end{align*}
First, because $\bX\bbeta = \bX_{\cS^0}\bbeta_{\cS^0}$ we see that
\begin{align*}
V_1&\leq s^{1/2}\left\|n^{-1}\wtilde{\bX}_{\cA}(\btheta)'\bX_{\cS^0}
- n^{-1}\wtilde{\bX}_{\cA}(\btheta^0)'\bX_{\cS^0}\right\|_{\max} \left\|\bbeta_{\cS^0}\right\|_2 \\
&\lesssim s\left\|n^{-1}\wtilde{\bX}_{\cA}(\btheta)'\bX_{\cS^0}
- n^{-1}\wtilde{\bX}_{\cA}(\btheta^0)'\bX_{\cS^0}\right\|_{\max}.
\end{align*}
Recall that $|\cS^0|=s$ and $s\leq n\wedge p$. By a similar bound of $U_2$, the norm just above can be bounded further as
\begin{align*}
&\left\|n^{-1/2}{\bC_{\cS^0}^0}\right\|_2 \left\|\bC_{\cA}-\bC_{\cA}^0 \right\|_{\max} 
+ \left\|n^{-1/2}{\bC_{\cS^0}^0}\right\|_2 \left\|\bE_{\bbbeta \cA}-\bE_{\bbbeta^0 \cA} \right\|_{\max} \\
&\qquad + \left\|n^{-1/2}\bE_{\bbbeta^0 {\cS^0}} \right\|_2  \left\|\bC_{\cA}-\bC_{\cA}^0 \right\|_{\max} 
+ \left\|n^{-1/2}\bE_{\bbbeta^0 {\cS^0}} \right\|_2 \left\|\bE_{\bbbeta \cA}-\bE_{\bbbeta^0 \cA} \right\|_{\max} \\
&\lesssim s^{1/2} \left\|\bC-\bC^0 \right\|_{\max} 
+ s^{1/2} \left\|\bE_{\bbbeta}-\bE_{\bbbeta^0} \right\|_{\max}
+ \left\|\bC-\bC^0 \right\|_{\max} 
+ \left\|\bE_{\bbbeta}-\bE_{\bbbeta^0} \right\|_{\max} \\
&=O\left(s^{1/2}c_{np} + s^{1/2}\tilde{c}_{np} + c_{np} + \tilde{c}_{np} \right)
=O\left(s^{1/2}\tilde{c}_{np} \right).
\end{align*}
Thus we have 
\begin{align*}
V_1\lesssim ss^{1/2}\tilde{c}_{np} = s^{3/2}\tilde{c}_{np}
\end{align*}
with probability at least $1-O(\pi_{np})$. 
Next the same procedure yields
\begin{align}
V_2 &\leq \left\|\wtilde{\bX}_{\cA}(\btheta)
- \wtilde{\bX}_{\cA}(\btheta^0) \right\|_{\max} \left\|n^{1/2} \bep \right\|_2 \notag\\
&\lesssim \left\|\bC-\bC^0 \right\|_{\max} 
+ \left\|\bE_{\bbbeta}-\bE_{\bbbeta^0} \right\|_{\max}
\lesssim \tilde{c}_{np},\label{ineq:V2}
\end{align}
where $\|n^{1/2} \bep \|_2=(\E\ep^2)^{1/2}+o(1)$ a.s.\ by the law of large numbers for independent random variables. 
Since the results hold uniformly in $\btheta\in\Theta_{np}$, combining them leads to 
\begin{align*}
\sup_{\btheta \in \Theta_{np}}\left\|\bv_{\cA}({\btheta}) - \bv_{\cA}(\btheta^0) \right\|_{\max} \lesssim s^{3/2}\tilde{c}_{np}
\end{align*}
with probability at least $1-O(\pi_{np})$. This concludes the proof of Lemma \ref{lem:diffsUV_max}.

\subsection{Lemma \ref{lem:maxeigv} and its proof} \label{secC.4}
\begin{lem}\label{lem:maxeigv}
	Assume that Condition \ref{ass:err_new} holds. Then the following statements hold
	\begin{enumerate}
		\item [(a)] We have
		\begin{align*}
		\Pro \left( \left\|(n\vee p)^{-1/2}\bE\right\|_2 \lesssim 1 \right) 
		\geq 1-2\exp(-|O(n\vee p)|);
		\end{align*}
		\item [(b)] For any fixed set $\cA$ with $|\cA|\leq k\leq n$, we have
		\begin{align*}
		\Pro \left( \left\|n^{-1/2}\bE_\cA \right\|_2 \lesssim 1 \right) 
		\geq 1-2p^{-\nu};
		\end{align*}
		\item [(c)] For all $k\leq n$, we have 
		\begin{align*}
		\Pro \left( \max_{|\cA|\leq k} \left\|n^{-1/2}\bE_\cA \right\|_2 \lesssim 1\vee \left(n^{-1}k\log p \right)^{1/2} \right) 
		\geq 1-2p^{-\nu},
		\end{align*}
	\end{enumerate}
	where $\nu>0$ is a predetermined constant. 
\end{lem}

\noindent\textit{Proof}. 
Result (a) is obtained by Theorem 5.39 of \cite{vershynin12}. Moreover, by the same theorem there exist some positive constants $c$ and $C$ such that for any $\cA$ with $|\cA|\leq k\leq n$ and every $t\geq 0$,
\begin{align*}
\Pro \left( \sigma_e^{-1} \|n^{-1/2}\bE_\cA\|_2 > 1+C + n^{-1/2}t \right) \leq 2 \exp \left(-ct^2 \right),
\end{align*}
where $\sigma_e^2=\E e^2$. Therefore, result (b) is immediately obtained by putting $t^2=c^{-1}\nu\log p$ since $n^{-1/2}t=o(1)$ and 
$\exp \left(-ct^2 \right) =p^{-\nu} $ in this case. 

For (c), taking the union bound leads to 
\begin{align*}
&\Pro \left( \sigma_e^{-1} \max_{|\cA|\leq k} \|n^{-1/2}\bE_\cA\|_2 > 1+C + n^{-1/2}t \right) \\
&\quad\leq \binom{p}{k} \max_{|\cA|\leq k} \Pro \left( \sigma_e^{-1} \|n^{-1/2}\bE_\cA\|_2 > 1+C + n^{-1/2}t \right)
\leq 2 p^k \exp \left(-ct^2 \right).
\end{align*}
Set $t^2=c^{-1}(\nu+k)\log p$ in this inequality. 
Then we have $n^{-1/2}t= O\left((n^{-1}k\log p)^{1/2}\right)$ and 
\begin{align*}
2 p^k \exp \left(-ct^2 \right) \leq 2 p^k \exp \left(-(\nu+k)\log p \right) = 2p^{-\nu},
\end{align*}
which gives result (c) and completes the proof of Lemma \ref{lem:maxeigv}.

\subsection{Lemma \ref{lem:matnorm} and its proof} \label{secC.5}
\begin{lem}\label{lem:matnorm}
	For matrices $\bA\in \mathbb{R}^{k_1\times n}$ and $\bB\in \mathbb{R}^{n\times k_2}$, we have $\|\bA\bB\|_{\max}\leq n^{1/2}\|\bA\|_2\|\bB\|_{\max}$ and $\|\bA\bB\|_{\max}\leq n^{1/2}\|\bA\|_{\max}\|\bB\|_2$. 
\end{lem}

\noindent\textit{Proof}. 
For any matrix $\bM=(m_{ij})\in\mathbb{R}^{k\times n}$, let $\|\bM\|_{\infty,\infty}$ denote the induced $\ell_\infty$-norm. 
First, we have
\begin{align*}
\|\bM\|_{\infty,\infty} := \sup_{\bv\in\mathbb{R}^n\backslash\{\bzero\}}\frac{\|\bM\bv\|_{\max}}{\|\bv\|_{\max}}
\leq \sup_{\bv\in\mathbb{R}^n\backslash\{\bzero\}}\frac{\|\bM\bv\|_2}{\|\bv\|_2} \frac{\|\bv\|_2}{\|\bv\|_{\max}}
\leq n^{1/2} \|\bM\|_2.
\end{align*}
Therefore, by a simple calculation we see that 
\begin{align*}
\|\bA\bB\|_{\max} 
&= \|\vect(\bA\bB)\|_{\max} = \|(\bI_{k_2}\otimes \bA)\vect(\bB)\|_{\max} \\
&= \frac{\|(\bI_{k_2}\otimes \bA)\vect(\bB)\|_{\max}}{\|\vect(\bB)\|_{\max}}\|\vect(\bB)\|_{\max} \\
&\leq \|\bI_{k_2} \otimes \bA\|_{\infty,\infty} \|\vect(\bB)\|_{\max} 
= \|\bA\|_{\infty,\infty} \|\bB\|_{\max} 
\leq n^{1/2}\|\bA\|_2 \|\bB\|_{\max}.
\end{align*}
The second assertion follows from applying this inequality to $\bB'\bA'$. 
This concludes the proof of Lemma \ref{lem:matnorm}.

\end{document}